\pdfoutput=1
\RequirePackage{ifpdf}
\ifpdf % We are running pdfTeX in pdf mode
\documentclass[pdftex]{sigma}
\else
\documentclass{sigma}
\fi

\numberwithin{equation}{section}

\newtheorem{Theorem}{Theorem}[section]
\newtheorem{Corollary}[Theorem]{Corollary}
\newtheorem{Conjecture}[Theorem]{Conjecture}
\newtheorem{rhp}[Theorem]{Riemann--Hilbert Problem}
 { \theoremstyle{definition}
\newtheorem{Remark}[Theorem]{Remark} }

\newcommand{\ii}{\mathrm{i}}
\newcommand{\ee}{\mathrm{e}}
\newcommand{\dd}{\mathrm{d}}

\begin{document}
%\allowdisplaybreaks

\newcommand{\arXivNumber}{1804.03173}

\renewcommand{\thefootnote}{}

\renewcommand{\PaperNumber}{125}

\FirstPageHeading

\ShortArticleName{On the Increasing Tritronqu\'ee Solutions of the Painlev\'e-II Equation}

\ArticleName{On the Increasing Tritronqu\'ee Solutions\\ of the Painlev\'e-II Equation\footnote{This paper is a~contribution to the Special Issue on Painlev\'e Equations and Applications in Memory of Andrei Kapaev. The full collection is available at \href{https://www.emis.de/journals/SIGMA/Kapaev.html}{https://www.emis.de/journals/SIGMA/Kapaev.html}}}

\Author{Peter D.~MILLER}

\AuthorNameForHeading{P.D.~Miller}

\Address{Department of Mathematics, University of Michigan,\\ East Hall, 530 Church St., Ann Arbor, MI 48109, USA}
\Email{\href{mailto:millerpd@umich.edu}{millerpd@umich.edu}}
\URLaddress{\url{http://www.math.lsa.umich.edu/~millerpd/}}

\ArticleDates{Received April 11, 2018, in final form November 12, 2018; Published online November 15, 2018}

\Abstract{The increasing tritronqu\'ee solutions of the Painlev\'e-II equation with parameter $\alpha$ exhibit square-root asymptotics in the maximally-large sector $|\arg(x)|<\tfrac{2}{3}\pi$ and have recently appeared in applications where it is necessary to understand the behavior of these solutions for complex values of $\alpha$. Here these solutions are investigated from the point of view of a Riemann--Hilbert representation related to the Lax pair of Jimbo and Miwa, which naturally arises in the analysis of rogue waves of infinite order. We show that for generic complex $\alpha$, all such solutions are asymptotically pole-free along the bisecting ray of the complementary sector $|\arg(-x)|<\tfrac{1}{3}\pi$ that contains the poles far from the origin. This allows the definition of a total integral of the solution along the axis containing the bisecting ray, in which certain algebraic terms are subtracted at infinity and the poles are dealt with in the principal-value sense. We compute the value of this integral for all such solutions. We also prove that if the Painlev\'e-II parameter $\alpha$ is of the form $\alpha=\pm\tfrac{1}{2}+\ii p$, $p\in\mathbb{R}\setminus\{0\}$, one of the increasing tritronqu\'ee solutions has no poles or zeros whatsoever along the bisecting axis.}

\Keywords{Painlev\'e-II equation; tronqu\'ee solutions}

\Classification{33E17; 34M40; 34M55; 35Q15}

\renewcommand{\thefootnote}{\arabic{footnote}}
\setcounter{footnote}{0}

\section{Introduction}
The Painlev\'e-II equation with parameter $\alpha\in\mathbb{C}$
\begin{gather}
\frac{\dd^2u}{\dd x^2}=xu+2u^3-\alpha,\qquad u=u(x;\alpha)\label{eq:PII-intro}
\end{gather}
has been the object of intense study ever since it was identified by Painlev\'e more than a century ago as one of only 6 formerly unknown second-order ordinary differential equations of the form $u''=F(x,u)$ with $F$ rational in $u$ and meromorphic in $x$ having what is now called the \emph{Painlev\'e property}: the only singularities of a solution $u$ whose location in the $x$-plane depends on initial conditions are poles. In fact, it is now known that every solution of \eqref{eq:PII-intro} is a meromorphic function of~$x$ all of whose poles are simple and of residue $\pm 1$. Since its original discovery in the context of the solution of an abstract classification problem for differential equations, the equation~\eqref{eq:PII-intro} and its particular solutions have become very important in numerous applications. For instance, similarity solutions of the modified Korteweg--de Vries equation satisfy~\eqref{eq:PII-intro}~\cite{FlaschkaN80}. The oscillations appearing near the leading edge of the dispersive shock wave generated from a wide class of initial data in the weakly-dispersive Korteweg--de Vries equation have a universal profile corresponding to the Hastings--McLeod solution of~\eqref{eq:PII-intro} with $\alpha=0$ \cite{ClaeysG10}. The real graphs of the rational solutions of \eqref{eq:PII-intro} for integer values of $\alpha$ determine the locations of kinks in space-time near a point where generic initial data for the semiclassical sine-Gordon equation crosses the separatrix in the phase portrait of the simple pendulum in a transversal manner~\cite{BuckinghamM12}. In mathematical physics there are also many applications of solutions of \eqref{eq:PII-intro}. Perhaps the most famous one concerns the distribution functions for the largest eigenvalue of random matrices from certain ensembles, which in the scaling limit can in some cases can be written in terms of again the Hastings--McLeod solution \cite{TracyW94}.

The Hastings--McLeod solution of \eqref{eq:PII-intro} is especially important in applications because it is a global solution for real values of $x$, i.e., it has no singularities for $x\in\mathbb{R}$. It is an example of a so-called \emph{tronqu\'ee} solution, namely one having no poles near $x=\infty$ in one or more sectors of opening angle $\tfrac{2}{3}\pi$ of the complex plane. In fact, the Hastings--McLeod solution has no poles near infinity in two disjoint sectors: $|\arg(x)|<\tfrac{1}{3}\pi$ and $|\arg(-x)|<\tfrac{1}{3}\pi$. The intervening sectors symmetric about the imaginary axis are filled with poles for large $|x|$, and the poles form a locally regular lattice. In general, there are two distinct types of tronqu\'ee solutions of the Painlev\'e-II equation \eqref{eq:PII-intro}. The \emph{increasing tronqu\'ee} solutions satisfy $u(x;\alpha)\sim \pm\big({-}\tfrac{1}{2}x\big)^{1/2}$ as $|x|\to\infty$ with $|\arg(-x)|<\tfrac{1}{3}\pi$ (arising from a dominant balance between the terms $xu$ and $2u^3$ in \eqref{eq:PII-intro}), and the \emph{decreasing tronqu\'ee} solutions satisfy $u(x;\alpha)\sim \alpha x^{-1}$ as $|x|\to\infty$ with $|\arg(x)|<\tfrac{1}{3}\pi$ (arising from a dominant balance between the terms $xu$ and $-\alpha$ in \eqref{eq:PII-intro}). The Hastings--McLeod solution is uniquely determined as being both an increasing and a decreasing tronqu\'ee solution. Otherwise, the tronqu\'ee solutions are not generally determined by their leading asymptotics. However, for each choice of sign there is a unique solution of \eqref{eq:PII-intro} for which $u(x;\alpha)\sim\pm \big({-}\tfrac{1}{2}x\big)^{1/2}$ holds as $|x|\to\infty$ with $-\tfrac{1}{3}\pi<\arg(-x)<\pi$ and another unique solution for which the same asymptotic holds for $-\pi<\arg(-x)<\tfrac{1}{3}\pi$. These solutions are related by the rotation symmetry of the Painlev\'e-II equation in which whenever $u(x)$ is a solution, then so is $\ee^{2\pi\ii/3}u\big(\ee^{2\pi\ii/3}x\big)$. Thus there is also for each choice of sign a unique solution for which $u(x;\alpha)\sim \pm\ii\big(\tfrac{1}{2}x\big)^{1/2}$ as $|x|\to\infty$ with $|\arg(x)|<\tfrac{2}{3}\pi$. These six solutions are called the \emph{increasing tritronqu\'ee} solutions of \eqref{eq:PII-intro} in the nomenclature of \cite[Chapter~11]{FokasIKN06} that can be traced back to the terminology introduced by Boutroux~\cite{Boutroux13,Boutroux14}.

The Painlev\'e-I equation also has tritronqu\'ee solutions, and these have been conjectured and/or shown to describe critical phenomena in several different situations \cite{BertolaT13,BuckinghamM15,DubrovinGK09,LuM18}. The tritronqu\'ee solutions of other Painlev\'e equations seem to not arise as frequently; however the increasing tritronqu\'ee solutions of \eqref{eq:PII-intro} have recently appeared in two quite different applications. The first application is the asymptotic description of rational solutions $w=w_n(z;m)$, $n\in\mathbb{Z}$, $m\in\mathbb{C}$, of the Painlev\'e-III equation
\begin{gather}
\frac{\dd^2w}{\dd z^2}=\frac{1}{w}\left(\frac{\dd w}{\dd z}\right)^2-\frac{1}{z}\frac{\dd w}{\dd z}+\frac{4(n+m)w^2+4(n-m)}{z}+4w^3-\frac{4}{w}\label{eq:PIII}
\end{gather}
in the limit of large integer parameter $n$. Given such $n$, the rational solution $w=w_n(z;m)$ is uniquely determined by the rationality condition and the asymptotic property $w_n(z;m)\to 1$ as $z\to\infty$. As $n\to+\infty$, the poles and zeros of $w_n(z;m)$ accumulate within a dilation $nE$ of a~fixed eye-shaped domain $E$ centered at the origin and with corners at the points $\pm\tfrac{1}{2}\ii$.
In the interior of the eye-shaped domain there are accurate asymptotic formul\ae\ for $w_n(z;m)$ in terms of modulated elliptic functions \cite{BothnerM18}. If one examines the function $w_n(z;m)$ near the corner point $z=\pm\tfrac{1}{2}\ii n$, it turns out to be natural to zoom in on the corner point by centering and rescaling the independent variable as $z=\pm\ii\big(\tfrac{1}{2}n+\big(\tfrac{1}{32}n\big)^{1/3}x\big)$ and also to introduce a new dependent variable by $w_n=\pm\ii \big(1-\big(\tfrac{1}{4}n\big)^{-1/3}u^\pm\big)$. Substituting these scalings into the Painlev\'e-III equation \eqref{eq:PIII} and formally considering $n$ large one finds that $u=u^\pm(x)$ is a solution of an $O\big(n^{-1/3}\big)$ perturbation of the Painlev\'e-II equation \eqref{eq:PII-intro} with parameter $\alpha=m$. Noting that the interior angle of the eye-shaped domain $E$ at its corners is exactly $2\pi/3$ and comparing with the known asymptotic behavior of $w_n(z;m)$ in the exterior of $nE$, in \cite{BothnerMS18} the following conjecture is formulated.
\begin{Conjecture}[Bothner, Miller, and Sheng]\label{conj:PIII}
Let $m\in\mathbb{C}$ be fixed. Then
\begin{gather*}
\lim_{n\to+\infty}\big(\tfrac{1}{4}n\big)^{1/3}(1\pm\ii w_n\big({\pm}\ii\big(\tfrac{1}{2}n+\big(\tfrac{1}{32}n\big)^{1/3}x\big);m\big)=u^\pm_\mathrm{TT}(x;m)
\end{gather*}
where $u=u^\pm_\mathrm{TT}(x;m)$ is the unique increasing tritronqu\'ee solution of~\eqref{eq:PII-intro} with parameter $\alpha=m$ and determined by the asymptotic behavior $u^\pm_\mathrm{TT}(x;m)\sim \pm\ii\big(\tfrac{1}{2}x\big)^{1/2}$ as $x\to\infty$ with $|\arg(x)|<\tfrac{2}{3}\pi$.
\end{Conjecture}
Thus the complementary sector $|\arg(-x)|\le\tfrac{1}{3}\pi$ in which the poles of the tritronqu\'ee solution reside for large $|x|$ corresponds to the interior of the eye-shaped domain $E$ in an overlap region where the Painlev\'e-II asymptotics and the modulated elliptic function asymptotics are simultaneously valid. Now the elliptic function asymptotics valid within $E$ are associated with an elliptic curve associated to each point of $E$, and in \cite{BothnerM18} it is shown that the elliptic curve degenerates along the vertical segment of the imaginary axis connecting the two corner points. This degeneration makes one of the periods of the elliptic function blow up and results in a~local dilution of the pole/zero distribution of $w_n(z;m)$ when $z$ is near the vertical segment. In the overlap domain this vertical segment corresponds to the negative real axis in the variable $x$ of the Painlev\'e-II increasing tritronqu\'ee solution. Of course the negative real axis is precisely the ray that bisects the asymptotic pole sector $|\arg(-x)|<\tfrac{1}{3}\pi$ for the tritronqu\'ee solution. This is one of the \emph{critical rays}\footnote{This is terminology used in \cite{Novokshenov12}. In other references the same rays are also called \emph{canonical rays} \cite[Chapters~9--10]{FokasIKN06} or \emph{Stokes rays} \cite[Remark~7.1]{FokasIKN06}.} for the Painlev\'e-II equation \eqref{eq:PII-intro}, and it is well known in the case $\alpha=0$ that solutions behave differently for large $x$ near such rays than elsewhere in the complex plane. See \cite[Chapters~9--10]{FokasIKN06}, where the large-$x$ asymptotics are worked out in detail for a class of solutions for $\alpha=0$ that includes the tritronqu\'ees $u^\pm_\mathrm{TT}(x;0)$ as a particular case. However, since $m\in\mathbb{C}$ is an arbitrary parameter in the sequence of rational Painlev\'e-III solutions $\{w_n(z;m)\}_{n=0}^\infty$, to fully explain the matching between the corner asymptotics suggested by Conjecture~\ref{conj:PIII} and the large-$n$ asymptotics of $w_n(z;m)$ along the central axis of the eye domain $z\in nE$, it is desirable to generalize known asymptotic results for the tritronqu\'ee solutions of~\eqref{eq:PII-intro} with $\alpha=0$ to the setting of arbitrary $\alpha=m\in\mathbb{C}$. This is the aim of Theorem~\ref{thm:centerline} below.

The second application concerns a new solution $\Psi(X,T)$ of the focusing nonlinear Schr\"odinger equation
\begin{gather}
\ii\frac{\partial\Psi}{\partial T} +\frac{1}{2}\frac{\partial^2\Psi}{\partial X^2} + |\Psi|^2\Psi=0,\label{eq:NLS}
\end{gather}
which was recently identified \cite{BilmanLM18} as a scaling limit of a sequence of particular solutions of the same equation modeling so-called rogue waves of increasingly higher amplitude. The special solution $\Psi(X,T)$, the \emph{rogue wave of infinite order}, has many remarkable additional properties. In particular, the function $\Psi(X,0)$ is related to a special transcendental solution of the Painlev\'e-III equation~\eqref{eq:PIII} with $n=-m=\tfrac{1}{2}$ for which all Stokes multipliers of the direct monodromy problem for \eqref{eq:PIII} vanish. In the regime of large variables $(X,T)\in\mathbb{R}^2$, it turns out that $\Psi(X,T)$ exhibits transitional asymptotic behavior when $v:=T|X|^{-3/2}$ is in the neighborhood of the critical value $v_\mathrm{c}:=54^{-1/2}$. In \cite{BilmanLM18} it is proved that as $|X|\to\infty$ with $v-v_\mathrm{c}=O\big(|X|^{-1/3}\big)$, the rogue wave of infinite order $\Psi(X,T)$ can be expressed explicitly in terms of a function $\mathcal{V}(y)$ extracted from a certain model Riemann--Hilbert problem, namely Riemann--Hilbert Problem~\ref{rhp:PII-generalized} in Section~\ref{sec:Representation} below in the case of parameters $p=\ln(2)/(2\pi)$ and $\tau=1$. It is then of some practical interest to obtain some alternative and possibly more effective characterization of $\mathcal{V}(y)$, and to determine its essential properties. One of the goals of this paper is to relate $\mathcal{V}(y)$ to the Painlev\'e-II equation \eqref{eq:PII-intro}, and to identify the particular solution needed to construct $\mathcal{V}(y)$ as the increasing tritronqu\'ee\footnote{We stress that the subscript in the notation $u_\mathrm{TT}^\pm(x;\alpha)$ is mnemonic for ``tritronqu\'ee'' and has nothing to do with the independent variable $T$ in \eqref{eq:NLS}.} $u^-_\mathrm{TT}(x;\alpha)$ for $\alpha=\tfrac{1}{2}+\ii\ln(2)/(2\pi)$. It is known \cite[p.~297]{FokasIKN06} that when $\alpha=0$, $u_\mathrm{TT}^\pm(x;0)$ is well-defined for all real $x$, however for $\mathcal{V}(y)$ to be a meaningful asymptotic description of the rogue wave of infinite order $\Psi(X,T)$, it would be necessary that $u^-_\mathrm{TT}(x;\alpha)$ be a global solution of \eqref{eq:PII-intro} also for the complex value $\alpha=\tfrac{1}{2}+\ii\ln(2)/(2\pi)$. The global nature of tritronqu\'ee solutions of~\eqref{eq:PII-intro} for certain complex $\alpha$ including this particular value is the subject of Theorem~\ref{thm:global} below.

One of the earliest and most important scientific contributions of Andrei Kapaev was a~systematic description of the large-$x$ asymptotics of general solutions of the Painlev\'e-II equation~\eqref{eq:PII-intro}; see for example \cite{Kapaev88,Kapaev92}. These works were based on the isomonodromy method in the setting of the Flaschka--Newell Lax pair representation \cite{FlaschkaN80} of~\eqref{eq:PII-intro}. Kapaev's results were made fully rigorous for $\alpha=0$ with the development by Deift and Zhou~\cite{DeiftZ95} of a suitable analogue of the steepest descent method adapted to matrix Riemann--Hilbert problems such as that arising in the inverse monodromy theory for~\eqref{eq:PII-intro}. The more general results of Kapaev were later also put on rigorous footing by this method; in particular the asymptotic analysis of the increasing tronqu\'ee solutions for general $\alpha$ is described fully for general $x$ avoiding the critical rays in \cite[Chapter~11, Section~5]{FokasIKN06}. Tritronqu\'ee solutions of higher-order equations in the Painlev\'e-II hierarchy have also been studied by Joshi and Mazzocco~\cite{JoshiM03}.

Despite these developments, the results needed for the applications described above do not appear to be in the literature. In this note, we try to fill this gap by proving the following results, in which we assume that
$\alpha\in\mathbb{C}\setminus\big(\mathbb{Z}+\tfrac{1}{2}\big)$, and let $u=u^\pm_\mathrm{TT}(x;\alpha)$ denote the corresponding increasing tritronqu\'ee solution of the Painlev\'e-II equation \eqref{eq:PII-intro} characterized uniquely by the asymptotics
\begin{gather}
u^\pm_\mathrm{TT}(x;\alpha)=\pm\ii\left(\frac{x}{2}\right)^{1/2}+O\big(x^{-1}\big),\qquad x\to\infty,\qquad |\arg(x)|<\frac{2}{3}\pi.\label{eq:u-x-wide-asymp}
\end{gather}
The uniqueness of $u^\pm_\mathrm{TT}(x;\alpha)$ given \eqref{eq:u-x-wide-asymp} combines with elementary Schwarz and odd-reflection symmetries of \eqref{eq:PII-intro} to imply the following identities:
\begin{gather}
u_\mathrm{TT}^+(x;\alpha)=-u_\mathrm{TT}^-(x;-\alpha)=u_\mathrm{TT}^-(x^*;\alpha^*)^*,\qquad (\alpha,x)\in\mathbb{C}^2.\label{eq:TT-symmetries}
\end{gather}
The remaining four increasing tritronqu\'ee solutions are obtained from $u^\pm_\mathrm{TT}(x;\alpha)$ via the cyclic symmetry group generated by $u(x)\mapsto \ee^{2\pi\ii/3}u\big(\ee^{2\pi\ii/3}x\big)$.

The first result concerns the behavior of increasing tritronqu\'ee solutions as $x\to\infty$ along the critical bisecting ray of the complementary sector $|\arg(-x)|<\tfrac{1}{3}\pi$ containing the poles near infinity.
Given $\alpha\in\mathbb{C}\setminus \big(\mathbb{Z}+\tfrac{1}{2}\big)$, let $q$ and $\tau$ be defined
as follows. Firstly, set
\begin{gather}
q_0(\alpha):=-\ii\alpha-\frac{1}{2\pi}\log\big(1+\ee^{-2\pi\ii\alpha}\big),\label{eq:q0}
\end{gather}
where to be precise the principal branch of the logarithm is meant, with imaginary part in $(-\pi,\pi]$. Then denoting by $[x]$ the nearest integer to $x\in\mathbb{R}$ with half-integers rounded down, i.e., $\big[n+\tfrac{1}{2}\big]=n$ for $n\in\mathbb{Z}$,
\begin{gather}
q=q(\alpha):=q_0(\alpha)-\ii [\mathrm{Im}(q_0(\alpha))],\qquad \tau=\tau(\alpha):=\ii (-1)^{[\mathrm{Im}(q_0(\alpha))]}\ee^{\log(1+\ee^{-2\pi\ii\alpha})/2}.\label{eq:q-tau}
\end{gather}
Note that $-\tfrac{1}{2}<\mathrm{Re}(\ii q(\alpha))\le\tfrac{1}{2}$.
\begin{Theorem}\label{thm:centerline}Suppose that $\mathrm{Re}(\ii q(\alpha))\neq\tfrac{1}{2}$. Then $u^-_\mathrm{TT}(x;\alpha)$ is pole free for sufficiently large negative $x$, and the following asymptotic formul\ae\ hold:
\begin{gather*}
u^-_\mathrm{TT}(x;\alpha)=\frac{q\Gamma(\ii q)\ee^{-\pi q/2}}{2\tau\sqrt{\pi}}\ee^{3\pi\ii/4}8^{-\ii q}\ee^{2\ii(-x)^{3/2}/3}(-x)^{-1/4-3\ii q/2}\big(1+O\big(|x|^{M^-(q)}\big)\big),\\
x\to -\infty,\qquad -\frac{1}{2}<\mathrm{Re}(\ii q)<0,\qquad
M^-(q):=\max\left\{-\frac{3}{4}-\frac{3}{2}\mathrm{Re}(\ii q),3\mathrm{Re}(\ii q)\right\}<0,
\\
u^-_\mathrm{TT}(x;\alpha)=\frac{\tau\sqrt{\pi}\ee^{\pi q/2}}{\Gamma(\ii q)}\ee^{-3\pi\ii/4}8^{\ii q}\ee^{-2\ii(-x)^{3/2}/3}(-x)^{-1/4+3\ii q/2}\big(1+O\big(|x|^{M^+(q)}\big)\big),\\
x\to -\infty,\quad 0<\mathrm{Re}(\ii q)<\frac{1}{2},\qquad M^+(q):=\max\left\{-\frac{3}{4}+\frac{3}{2}\mathrm{Re}(\ii q),-3\mathrm{Re}(\ii q)\right\}<0,
\end{gather*}
and
\begin{gather}
u^-_\mathrm{TT}(x;\alpha)=\frac{\tau\sqrt{2q(\ee^{2\pi q}-1)}}{(-x)^{1/4}}\sin(\theta(x))+O\big(|x|^{-1}\big),\qquad x\to -\infty,\nonumber\\
 q\in\mathbb{R},\qquad \theta(x):=-\frac{2}{3}(-x)^{3/2}+\frac{3}{2}q\ln(-x)-\frac{1}{4}\pi+3q\ln(2)-\arg(\Gamma(\ii q)).
\label{eq:uTT-x-negative-oscillatory}
\end{gather}
\end{Theorem}
Corresponding formul\ae\ for $u_\mathrm{TT}^+(x;\alpha)$ can be obtained from Theorem~\ref{thm:centerline} using the symmet\-ries~\eqref{eq:TT-symmetries} provided that $\mathrm{Re}(\ii q(-\alpha))\neq\tfrac{1}{2}$ or $\mathrm{Re}(\ii q(\alpha^*))\neq\tfrac{1}{2}$. In fact, the proof of Theorem~\ref{thm:centerline} will show that it is also possible to obtain an asymptotic description of $u^-_\mathrm{TT}(x;\alpha)$ when $\mathrm{Re}(\ii q(\alpha))=\tfrac{1}{2}$, but it may be necessary to exclude certain neighborhoods of infinitely many values of $x$ where poles may exist, and the resulting formula must include more terms.

\begin{figure}[t]\centering
\includegraphics[scale=0.75]{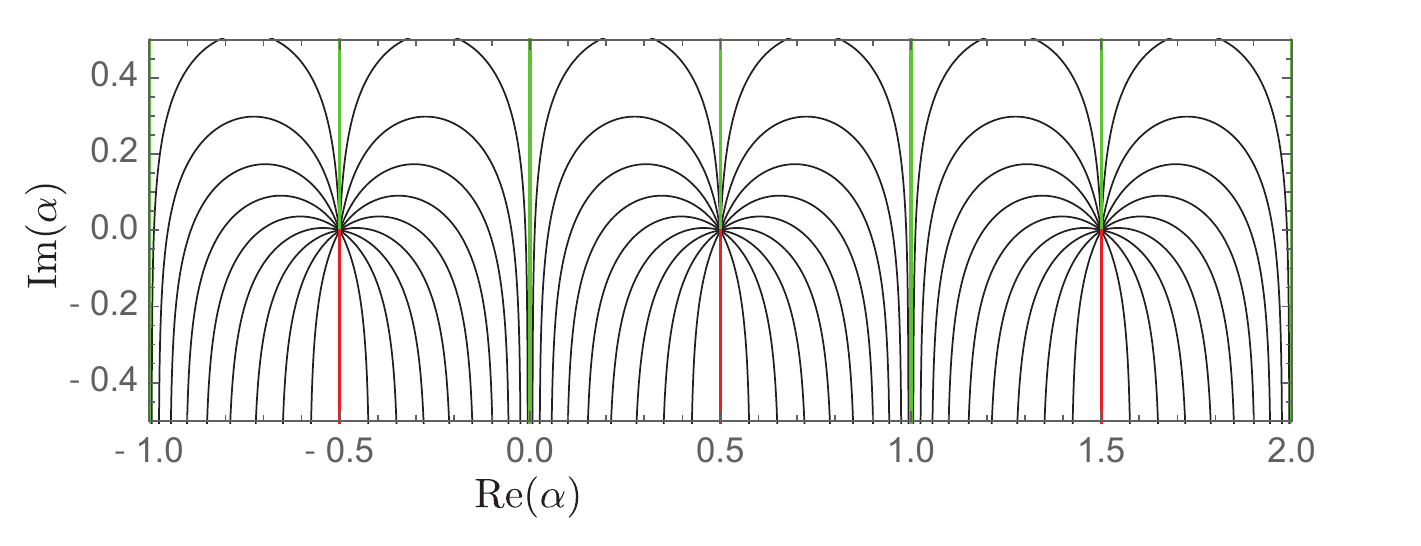}
\caption{Curves of constant $\mathrm{Re}(\ii q(\alpha))$ in the $\alpha$-plane. The red lines correspond to $\mathrm{Re}(\ii q(\alpha))=\tfrac{1}{2}$, i.e., everywhere but on these lines $u^-_\mathrm{TT}(x;\alpha)$ is pole-free for large real $x$ and has simple asymptotics given by Theorem~\ref{thm:centerline}. The green lines correspond to $\mathrm{Re}(\ii q(\alpha))=0$ and here $u^-_\mathrm{TT}(x;\alpha)$ has oscillatory asymptotics as $x\to -\infty$. The vertical strips $\mathrm{Re}(\alpha)\in \big(0,\tfrac{1}{2}\big)\pmod{\mathbb{Z}}$ correspond to $\mathrm{Re}(\ii q(\alpha))\in \big(0,\tfrac{1}{2}\big)$ while the vertical strips $\mathrm{Re}(\alpha)\in \big({-}\tfrac{1}{2},0\big)\pmod{\mathbb{Z}}$ correspond to $\mathrm{Re}(\ii q(\alpha))\in \big({-}\tfrac{1}{2},0\big)$. Note that $\mathrm{Re}(\ii q(\alpha))$ takes on every value in the range $\big({-}\tfrac{1}{2},\tfrac{1}{2}\big]$ in the neighborhood of each half-integer $\alpha=n+\tfrac{1}{2}$, $n\in\mathbb{Z}$.}
\label{fig:alpha-plane}
\end{figure}

\begin{figure}[t]\centering
\includegraphics[scale=0.875]{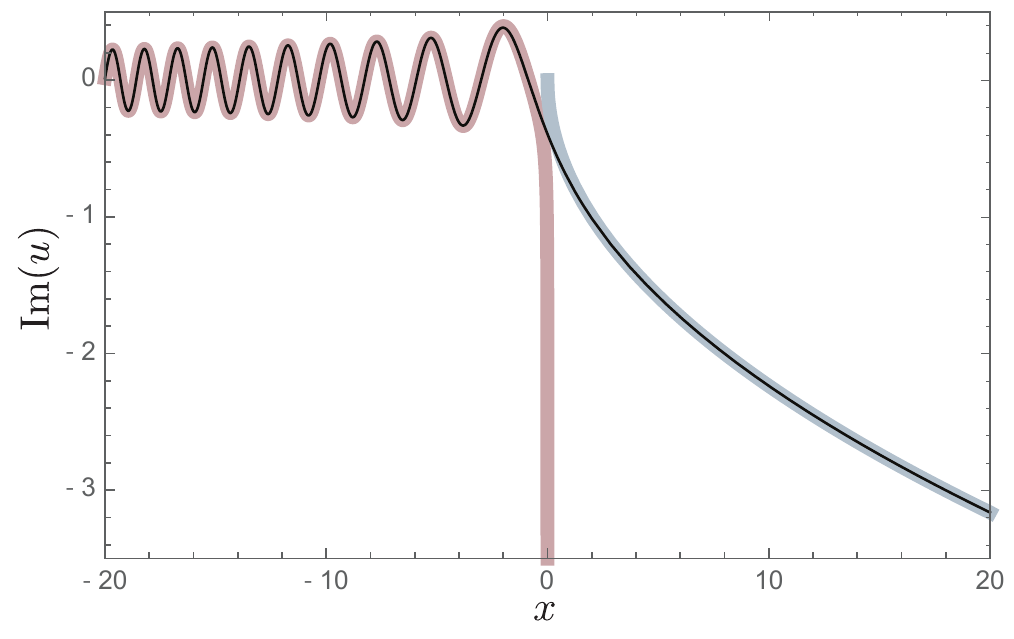}
\caption{Comparison of the leading terms in the asymptotic formul\ae\ \eqref{eq:u-x-wide-asymp} (light blue thick curve plotted for $x>0$) and \eqref{eq:alpha-zero-x-negative} (light red thick curve plotted for $x<0$) with Olver's numerical approximation of the purely imaginary solution $u_\mathrm{TT}^-(x;0)$ (thin black curve).}\label{fig:alpha-zero-compare}
\end{figure}

The only values of $\alpha$ that are not covered by Theorem~\ref{thm:centerline} are those for which $\mathrm{Re}(\ii q(\alpha))=\tfrac{1}{2}$. This means that $\mathrm{Re}(\ii q_0(\alpha))\in\mathbb{Z}+\tfrac{1}{2}$, i.e., $\log\big(1+\ee^{-2\pi\ii\alpha}\big)$ has the form $-2\pi\ii\alpha + r+ 2\pi\ii \big(n+\tfrac{1}{2}\big)$ where $r\in\mathbb{R}$ and $n\in\mathbb{Z}$ are arbitrary parameters. No information is therefore lost by exponentiating, which leads to $\ee^{-2\pi\ii\alpha}=-1/(1+\ee^r)$. It is then straightforward to determine that the excluded values of $\alpha$ have the form $\alpha=\tfrac{1}{2}+n-\ii p$ for $n\in\mathbb{Z}$ and $p>0$. Similarly, the values of $\alpha$ for which $\mathrm{Re}(\ii q(\alpha))=0$ (so the $x\to -\infty$ asymptotics are oscillatory) correspond to $\mathrm{Re}(\ii q_0(\alpha))\in\mathbb{Z}$, i.e., $\log\big(1+\ee^{-2\pi\ii\alpha}\big)$ has the form $-2\pi\ii\alpha+r+2\pi\ii n$ where $r\in\mathbb{R}$ and $n\in\mathbb{Z}$. Therefore $\ee^{-2\pi\ii\alpha}=1/(\ee^r-1)$. The case $r>0$ then corresponds to $\alpha\in\mathbb{Z}+\ii\mathbb{R}$ while $r<0$ corresponds instead to $\alpha=\tfrac{1}{2}+n+\ii p$ for $n\in\mathbb{Z}$ and $p>0$. See Fig.~\ref{fig:alpha-plane}.
Specific examples that are especially relevant include:
\begin{gather*}
q_0(0)=-\frac{\ln(2)}{2\pi}\quad\implies\quad q(0)=-\frac{\ln(2)}{2\pi}\qquad\text{and}\qquad \tau(0)=\ii\sqrt{2}
\end{gather*}
and, for applications to rogue waves of infinite order,
\begin{gather*}
q_0\left(\frac{1}{2}+\ii\frac{\ln(2)}{2\pi}\right)=
-\ii+\frac{\ln(2)}{2\pi}\quad\implies\\
q\left(\frac{1}{2}+\ii\frac{\ln(2)}{2\pi}\right)=\frac{\ln(2)}{2\pi}\qquad\text{and}\qquad \tau\left(\frac{1}{2}+\ii\frac{\ln(2)}{2\pi}\right)=1.
\end{gather*}
In both cases, the relevant asymptotic formula for $u_\mathrm{TT}^-(x;\alpha)$ in the limit $x\to-\infty$ is~\eqref{eq:uTT-x-negative-oscillatory}. In particular, we have
\begin{gather}
u_\mathrm{TT}^-(x;0)=\ii\sqrt{\frac{\ln(2)}{\pi}}\frac{1}{(-x)^{1/4}}\sin\left(-\frac{2}{3}(-x)^{3/2}-\frac{3\ln(2)}{4\pi}\ln(-x)-\frac{1}{4}\pi\right.\nonumber\\
\left. \hphantom{u_\mathrm{TT}^-(x;0)=}{} -\frac{3\ln(2)^2}{2\pi} +\arg\left(\Gamma\left(\ii\frac{\ln(2)}{2\pi}\right)\right)\right)
 +O\big(|x|^{-1}\big),\qquad x\to -\infty,\label{eq:alpha-zero-x-negative}
\end{gather}
a connection formula that is well-known in the literature, cf.~\cite[Theorem~9.1]{FokasIKN06} and \cite{Novokshenov12}.
In the Flaschka--Newell theory \cite{FlaschkaN80}, the solution $u^-_\mathrm{TT}(x;0)$ is associated with Stokes multipliers $s_1=s_2=s_3=\ii$; see \cite{Novokshenov12}. It is striking to compare the (purely imaginary) exact solution $u^-_\mathrm{TT}(x;0)$ with its asymptotic approximations for large positive and negative~$x$ (see~\eqref{eq:u-x-wide-asymp} and~\eqref{eq:alpha-zero-x-negative}, respectively).
To do this, we used the \textit{Mathematica} package \texttt{RHPackage} of Olver~\cite{RHPackage} with the command \texttt{PainleveII[\{I,I,I\},x]} to obtain the numerical approximation of the solution, which we compare with the two asymptotic formul\ae\ in Fig.~\ref{fig:alpha-zero-compare}.

\looseness=-1 The fact that the solution $u^\pm_\mathrm{TT}(x;\alpha)$ has simple asymptotics as $x\to\pm\infty$ suggests that, after subtracting off certain explicit terms, $u^\pm_\mathrm{TT}(x;\alpha)$ may be integrable over $x\in\mathbb{R}$ in a suitable sense. To make the definition of such an integral precise, we may recall that the only possible singularities of each solution of \eqref{eq:PII-intro} are simple poles of residue $\pm 1$, which may in principle occur along the real axis. Hence they may be taken into account by a suitable regularization of the integral. We choose the Hadamard principal value, and then we can establish the following result.
\begin{Theorem}\label{thm:integral}
Suppose that $\mathrm{Re}(\ii q(\alpha))\neq \tfrac{1}{2}$, and let $A<0<B$ be such that all real poles of $u^-_\mathrm{TT}(x;\alpha)$ lie in the interval $(A,B)$. Then the total integral formula
\begin{gather}
\exp\left(\int_{-\infty}^{A}\left[u^-_\mathrm{TT}(x;\alpha)-\frac{\alpha}{x}\right]\,\dd x +
\mathrm{P.V.} \int_{A}^{B}u^-_\mathrm{TT}(x;\alpha)\,\dd x \right.\nonumber\\
\left. \qquad \quad {}+
\int_{B}^{+\infty}\left[u^-_\mathrm{TT}(x;\alpha)+\ii\sqrt{\frac{x}{2}}+\frac{\alpha}{2x}\right]\,\dd x\right)\nonumber\\
\qquad {} =(-1)^{N_+-N_-}\tau(\alpha)\frac{\ee^{\ii\pi(\alpha-1)/2}\Gamma(\alpha+\tfrac{1}{2})}{\sqrt{2\pi}}
\left(-A\sqrt{\frac{B}{2}}\right)^{-\alpha}\ee^{-\ii B^{3/2}\sqrt{2}/3}\label{eq:u-total-integral}
\end{gather}
holds, where $N_\pm$ denotes the number of real poles of $u^-_\mathrm{TT}(x;\alpha)$ of residue $\pm 1$, and in which the integral over $(-\infty,A)$ is a convergent improper integral while that over $(B,+\infty)$ is absolutely convergent.
\end{Theorem}
Again, a corresponding formula for $u_\mathrm{TT}^+(x;\alpha)$ can be obtained from \eqref{eq:u-total-integral} using the symmetries \eqref{eq:TT-symmetries} provided that $\mathrm{Re}(\ii q(-\alpha))\neq\tfrac{1}{2}$ or $\mathrm{Re}(\ii q(\alpha^*))\neq\tfrac{1}{2}$. Similar results have been obtained for other well-known solutions of the Painlev\'e-II equation such as the Hastings--McLeod and Ablowitz-Segur solutions for $\alpha=0$, see, e.g.,~\cite{BaikBDI09}, although usually such integration formul\ae\ have been considered only for global (i.e., pole-free on the integration axis) solutions.

Next, we recall that in the setting of the rogue wave solution $\Psi(X,T)$ of infinite order, we need to consider the special complex value of $\alpha=\tfrac{1}{2}+\ii\ln(2)/(2\pi)$ and know that the solution $u^-_\mathrm{TT}(x;\alpha)$ has no poles at all on the real line. In fact, we can show more.
\begin{Theorem}\label{thm:global}
Suppose that $p>0$. Then $u^-_\mathrm{TT}(x;\tfrac{1}{2}+\ii p)$, $u^-_\mathrm{TT}\big(x;-\tfrac{1}{2}+\ii p\big)$, $u^+_\mathrm{TT}\big(x;\tfrac{1}{2}-\ii p\big)$, and $u^+_\mathrm{TT}\big(x;-\tfrac{1}{2}-\ii p\big)$ are
global solutions for $x\in\mathbb{R}$, i.e., they are analytic for $x\in\mathbb{R}$. Moreover, they have no real zeros.
\end{Theorem}
The reader may observe that, in the cases covered by Theorem~\ref{thm:global}, the relevant solution always has oscillatory asymptotics as $x\to -\infty$ according to \eqref{eq:uTT-x-negative-oscillatory} in Theorem~\ref{thm:centerline}, and that the leading term has infinitely many zeros in the limit. Indeed, for $u_\mathrm{TT}^-(x;\alpha)$ the indicated values of $\alpha$ in Theorem~\ref{thm:global} lie on the green half-lines emerging vertically from $\alpha=\pm\tfrac{1}{2}$ in Fig.~\ref{fig:alpha-plane}. But whereas the solutions $u^\pm_\mathrm{TT}(x;0)$ are purely imaginary and hence the zeros of the leading term are perturbations of actual real zeros of the solution, the tritronqu\'ee solutions that are the subject of Theorem~\ref{thm:global} are essentially complex-valued. Thus, the fact that they have no real zeros simply means that the error term in~\eqref{eq:uTT-x-negative-oscillatory} is nonzero and has a phase with a component orthogonal to that of $\tau\sqrt{2q(\ee^{2\pi q}-1)}$ in neighborhoods of $x$-values satisfying $\theta(x)=\pi n$ for $n\in\mathbb{Z}$ large.

The reason for excluding the half-integral values of $\alpha$ in the above results is that the Riemann--Hilbert representation of the increasing tritronqu\'ee solution that we study below fails to yield a~determinate expression for the solution in such cases\footnote{As will be seen early in Section~\ref{sec:Representation}, if $\alpha-\tfrac{1}{2}\in\mathbb{Z}\setminus\{0\}$ then Riemann--Hilbert Problem~\ref{rhp:PII-generalized} below has no solution at all, while if $\alpha=\tfrac{1}{2}$ it has a trivial solution through which $u_\mathrm{TT}^\pm(x;\tfrac{1}{2})$ is represented as an indeterminate fraction: $0/0$.}. On the other hand, it was discovered by Gambier \cite{Gambier10} and is now well-known that for $\alpha-\tfrac{1}{2}\in\mathbb{Z}$, the Painlev\'e-II equation \eqref{eq:PII-intro} is solvable via B\"acklund transformations and the general solution $a\mathrm{Ai}(x)+b\mathrm{Bi}(x)$ of the Airy equation (see \cite[Chapter~11, Section~4]{FokasIKN06}); the tronqu\'ee solutions of these special cases were recently studied in detail by Clarkson \cite{Clarkson16}. However, it is clear from Theorem~\ref{thm:centerline} and Fig.~\ref{fig:alpha-plane} that the asymptotic behavior of $u^\pm_\mathrm{TT}(x;\alpha)$ is very sensitive to the value of $\alpha$ near half-integers $\mathbb{Z}+\tfrac{1}{2}$. This raises the interesting question of double-scaling asymptotics, i.e., consideration of the limits $x\to -\infty$ and $\alpha\to n+\tfrac{1}{2}$ simultaneously at appropriate related rates, a problem for the future. A different double-scaling limit related to solutions of~\eqref{eq:PII-intro} has been recently addressed by Bothner~\cite{Bothner17}, and the joint asymptotic behavior of $u(x;\alpha)$ when $x$ and $\alpha$ are both large has also been studied~\cite{BuckinghamM14,BuckinghamM15,Kapaev97}.

Finally, we formulate a corollary of the above results and their connection with Riemann--Hilbert Problem~\ref{rhp:PII-generalized} formulated in Section~\ref{sec:Representation} below that is needed in the application to rogue waves of infinite order \cite{BilmanLM18}.
\begin{Corollary}Let $p:=\ln(2)/(2\pi)$ and $\tau=1$. Then Riemann--Hilbert Problem~{\rm \ref{rhp:PII-generalized}} has a unique solution for every $y\in\mathbb{R}$, and the function $\mathcal{V}(y)$ extracted from it via the formula \eqref{eq:UV-define} is zero-free and critical point-free for real $y$ and satisfies $\mathcal{V}'(y)/\mathcal{V}(y)=-\big(\tfrac{2}{3}\big)^{1/3}u^-_\mathrm{TT}\big({-}\big(\tfrac{2}{3}\big)^{1/3}y;\tfrac{1}{2}+\ii
p\big)$. Moreover, $\mathcal{V}(y)$ has the following asymptotic behavior:
\begin{gather}
\mathcal{V}(y)=\frac{\tau p\Gamma(\ii p)}{2\sqrt{\pi}}\ee^{-3\pi\ii/4}\ee^{-\pi p/2}2^{-\ii p}\ee^{-2\ii (-y/3)^{3/2}}\nonumber\\
\hphantom{\mathcal{V}(y)=}{}\times (-3y)^{-\ii p/2-1/4}\big[1+O\big(|y|^{-5/4}\big)\big],\qquad y\to -\infty,\label{eq:V-asymp-y-Minus}
\end{gather}
and
\begin{gather}
\mathcal{V}(y)=-\left(\frac{y}{6}\right)^{\ii p+1/2}\big(1+O\big(y^{-3/4}\big)\big),\qquad y\to +\infty.\label{eq:V-asymp-y-Plus-specific}
\end{gather}
\end{Corollary}

The rest of this paper is organized as follows. In Section~\ref{sec:Representation} we present a Riemann--Hilbert problem connected with the Jimbo--Miwa theory of the Painlev\'e-II equation \eqref{eq:PII-intro} and use the Deift-Zhou steepest descent method to study its asymptotic behavior and therefore establish precisely which solution of the latter equation it encodes, namely the increasing tritronqu\'ee solution $u=u_\mathrm{TT}^-(x;\alpha)$. Then in Sections~\ref{sec:centerline}, \ref{sec:integral}, and \ref{sec:global} we use the Riemann--Hilbert representation to prove Theorems~\ref{thm:centerline}, \ref{thm:integral}, and \ref{thm:global}, respectively. A certain parabolic cylinder parametrix needed in Sections~\ref{sec:Representation} and \ref{sec:centerline} is described in Appendix~\ref{app:PC}.

\subsection*{Notation} We use subscripts $+$ and $-$ to denote the boundary values taken by a sectionally analytic function on an oriented jump contour from the left and right sides respectively. We also make frequent use of the Pauli matrices
\begin{gather*}
\sigma_1:=\begin{bmatrix}0&1\\1&0\end{bmatrix},\qquad\sigma_2:=\begin{bmatrix}0 & -\ii\\\ii & 0\end{bmatrix},\qquad\text{and}\qquad\sigma_3:=\begin{bmatrix}1&0\\0 & -1\end{bmatrix}.
\end{gather*}

\section[A Jimbo--Miwa representation of $u^-_\mathrm{TT}(x;\alpha)$]{A Jimbo--Miwa representation of $\boldsymbol{u^-_\mathrm{TT}(x;\alpha)}$}\label{sec:Representation}
\begin{figure}[t]\centering
\includegraphics[scale=0.9]{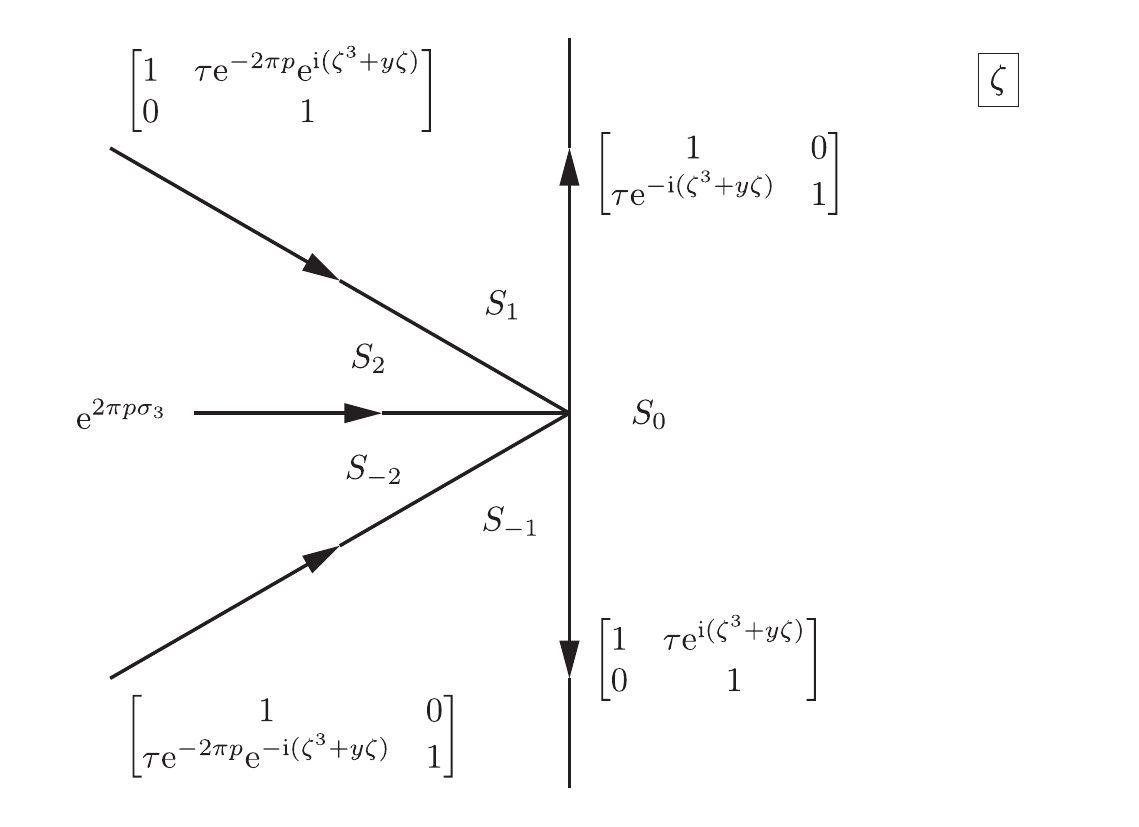}
\caption{The jump contour in the $\zeta$-plane and jump matrix $\mathbf{V}^{\mathrm{PII}}$.}\label{fig:PII-jumps-generalized}
\end{figure}

In \cite{BilmanLM18}, the analysis of the rogue wave of infinite order $\Psi(X,T)$ in the transitional regime where $T|X|^{-3/2}\approx 54^{-1/2}$ leads to consideration of a certain local model Riemann--Hilbert problem which coincides with the following in the special case of parameters $p=\ln(2)/(2\pi)$ and $\tau=1$. Consider the jump contour and jump matrix shown in Fig.~\ref{fig:PII-jumps-generalized}.
\begin{rhp}[Jimbo--Miwa Painlev\'e-II problem]\label{rhp:PII-generalized} Let $y,p,\tau\in\mathbb{C}$ be related by $\tau^2=\ee^{2\pi p}-1$. Seek a $2\times 2$ matrix-valued function $\mathbf{W}(\zeta;y)=\mathbf{W}(\zeta;y,p,\tau)$ with the following properties.
\begin{itemize}\itemsep=0pt
\item[]\textbf{Analyticity:} $\mathbf{W}(\zeta;y)$ is analytic for $\zeta$ in the five sectors $S_0$: $|\arg(\zeta)|<\tfrac{1}{2}\pi$, $S_1$: $\tfrac{1}{2}\pi<\arg(\zeta)<\tfrac{5}{6}\pi$, $S_{-1}$: $-\tfrac{5}{6}\pi<\arg(\zeta)<-\tfrac{1}{2}\pi$, $S_2$: $\tfrac{5}{6}\pi<\arg(\zeta)<\pi$, and $S_{-2}$: $-\pi<\arg(\zeta)<-\tfrac{5}{6}\pi$. It takes continuous boundary values on the excluded rays and at the origin from each sector.
\item[]\textbf{Jump conditions:} $\mathbf{W}_+(\zeta;y)=\mathbf{W}_-(\zeta;y)\mathbf{V}^\mathrm{PII}(\zeta;y)$, where $\mathbf{V}^\mathrm{PII}(\zeta;y)$ is the matrix defined on the jump contour shown in Fig.~{\rm \ref{fig:PII-jumps-generalized}}.
\item[]\textbf{Normalization:} $\mathbf{W}(\zeta;y)\zeta^{\ii p\sigma_3}\to\mathbb{I}$ as $\zeta\to\infty$ uniformly in all directions.
\end{itemize}
\end{rhp}
Note that the condition $\tau^2=\ee^{2\pi p}-1$ ensures that the cyclic product of the jump matrices at the origin is the identity, which is a necessary condition for the continuity of the boundary values from each sector at $\zeta=0$. Noting also that all jump matrices have unit determinant, it is therefore a simple consequence of Liouville's theorem that this problem has at most one solution, and if it exists it must have unit determinant. If $p\in\ii\mathbb{Z}$ then all jump matrices become the identity so given the continuity of the boundary values at the origin from each sector the solution $\mathbf{W}(\zeta;y)$ would need to be entire; but this yields a contradiction with the normalization condition unless also $p=0$. Hence there is no solution for $p\in\ii\mathbb{Z}\setminus\{0\}$. If $p=0$ it is easy to check that $\mathbf{W}(\zeta;y)\equiv\mathbb{I}$ is
the unique solution. For all other values of $p\in\mathbb{C}$, the solution will exist for generic values of $y\in\mathbb{C}$ that avoid certain poles.

\subsection{Differential equations}\label{sec:DE}
Given parameters $p$ and $\tau$ with $\tau^2=\ee^{2\pi p}-1$, and assuming solvability of Riemann--Hilbert Problem~\ref{rhp:PII-generalized} in the neighborhood of some value of $y\in\mathbb{C}$, we can derive from the solution certain differential equations via the dressing construction. It is a consequence of the exponential decay to the identity of the jump matrix $\mathbf{V}^\mathrm{PII}(\zeta;y)$ as $\zeta\to\infty$ that the normalization condition on~$\mathbf{W}(\zeta;y)$ holds in the stronger sense that:
\begin{gather}\allowdisplaybreaks
\mathbf{W}(\zeta;y)\sim\left(\mathbb{I}+\sum_{j=1}^\infty \mathbf{W}^j(y)\zeta^{-j}\right)\zeta^{-\ii p\sigma_3},\qquad \zeta\to\infty,\nonumber\\
\frac{\partial\mathbf{W}}{\partial y}(\zeta;y) \sim\left(\sum_{j=1}^\infty\frac{\dd\mathbf{W}}{\dd y}^{j}(y)\zeta^{-j}\right)\zeta^{-\ii p\sigma_3},\qquad\zeta\to\infty,\label{eq:large-zeta-PII}\\
\frac{\partial\mathbf{W}}{\partial\zeta}(\zeta;y) \sim-\left(\ii p\sigma_3\zeta^{-1}+\!\sum_{j=2}^\infty\! \big[(j-1)\mathbf{W}^{j-1}(y)+\ii p\mathbf{W}^{j-1}(y)\sigma_3\big]\zeta^{-j}\right)\zeta^{-\ii p\sigma_3},\qquad\zeta\to\infty.\nonumber
\end{gather}
It is easy to see that the ``undressed'' matrix $\mathbf{Z}(\zeta;y):=\mathbf{W}(\zeta;y)\ee^{\ii(\zeta^3+y\zeta)\sigma_3/2}$ is analytic in the same domain that $\mathbf{W}(\zeta;y)$ is, and satisfies analogous jump conditions except that the factors $\ee^{\pm\ii (\zeta^3+y\zeta)}$ in the jump matrix $\mathbf{V}^\mathrm{PII}(\zeta;y)$ have been replaced in all cases by $1$. It follows easily that
$\mathbf{Z}_\zeta(\zeta;y)\mathbf{Z}(\zeta;y)^{-1}$ and $\mathbf{Z}_y(\zeta;y)\mathbf{Z}(\zeta;y)^{-1}$ are both entire functions of $\zeta$ whose asymptotic expansions as $\zeta\to\infty$ are easily computed from \eqref{eq:large-zeta-PII}. Applying Liouville's theorem shows that these entire functions are polynomials:
\begin{gather}
\mathbf{Z}_\zeta\mathbf{Z}^{-1}=\frac{3}{2}\ii\sigma_3\zeta^2 + \frac{3}{2}\ii \big[\mathbf{W}^1(y),\sigma_3\big]\zeta\nonumber\\
\hphantom{\mathbf{Z}_\zeta\mathbf{Z}^{-1}=}{} +\frac{1}{2}\ii y\sigma_3 + \frac{3}{2}\ii \big[\mathbf{W}^2(y),\sigma_3\big]-\frac{3}{2}\ii\big[\mathbf{W}^1(y),\sigma_3\big]\mathbf{W}^1(y),
\label{eq:Z-zeta-Z-inverse}\\
\mathbf{Z}_y\mathbf{Z}^{-1}=\frac{1}{2}\ii\sigma_3\zeta +\frac{1}{2}\ii\big[\mathbf{W}^1(y),\sigma_3\big],\nonumber
\end{gather}
and from the necessarily vanishing coefficient of $\zeta^{-1}$ in the expansion of $\mathbf{Z}_\zeta(\zeta;y)\mathbf{Z}(\zeta;y)^{-1}$ one finds that
\begin{gather}
\frac{3}{2}\big[\mathbf{W}^3(y),\sigma_3\big]-\frac{3}{2}\big[\mathbf{W}^{2}(y),\sigma_3\big]\mathbf{W}^{1}(y) +\frac{3}{2}\big[\mathbf{W}^1(y),\sigma_3\big]\big(\mathbf{W}^1(y)^2-\mathbf{W}^{2}(y)\big)\nonumber\\
\qquad{} +\frac{1}{2} y\big[\mathbf{W}^1(y),\sigma_3\big]= p\sigma_3,\label{eq:Z-zeta-Z-inverse-zeta-m-1}
\end{gather}
while setting to zero the coefficient of $\zeta^{-1}$ in the expansion of $\mathbf{Z}_y(\zeta;y)\mathbf{Z}(\zeta;y)^{-1}$ gives
\begin{gather}
\frac{\dd\mathbf{W}^1}{\dd y}(y)+\frac{1}{2}\ii\big[\mathbf{W}^2(y),\sigma_3\big]-\frac{1}{2}\ii\big[\mathbf{W}^1(y),\sigma_3\big]\mathbf{W}^1(y)=\mathbf{0}.
\label{eq:Z-y-Z-inverse-zeta-m-1}
\end{gather}
The latter allows \eqref{eq:Z-zeta-Z-inverse} to be rewritten in terms of the matrix coefficient $\mathbf{W}^1(y)$ alone:
\begin{gather}
\mathbf{Z}_\zeta\mathbf{Z}^{-1}=\frac{3}{2}\ii\sigma_3\zeta^2 + \frac{3}{2}\ii\big[\mathbf{W}^1(y),\sigma_3\big]\zeta
+ \frac{1}{2}\ii y\sigma_3 - 3\frac{\dd\mathbf{W}^1}{\dd y}(y).\label{eq:Z-zeta-Z-inverse-OD}
\end{gather}
It is convenient to use the representation \eqref{eq:Z-zeta-Z-inverse} for the diagonal terms and \eqref{eq:Z-zeta-Z-inverse-OD} for the off-diagonal terms. Thus, if
\begin{gather}
\mathcal{U}(y):=W^1_{12}(y)\qquad\text{and}\qquad\mathcal{V}(y):=W^1_{21}(y),\label{eq:UV-define}
\end{gather}
then $\mathbf{Z}(\zeta;y)$ is a simultaneous fundamental solution matrix of the two Lax pair equations:
\begin{gather*}
\frac{\partial\mathbf{Z}}{\partial\zeta}=\mathbf{A}\mathbf{Z},\qquad\mathbf{A}:=\frac{3}{2}\ii\sigma_3\zeta^2 + 3\ii\begin{bmatrix}0 & -\mathcal{U}\\\mathcal{V}&0\end{bmatrix}\zeta + \frac{1}{2}\ii\begin{bmatrix}y+6\mathcal{U}\mathcal{V} & 6\ii \mathcal{U}'\\6\ii\mathcal{V}' & -y-6\mathcal{U}\mathcal{V}\end{bmatrix},
\end{gather*}
and
\begin{gather*}
\frac{\partial\mathbf{Z}}{\partial y}=\mathbf{B}\mathbf{Z},\qquad\mathbf{B}:=\frac{1}{2}\ii\sigma_3\zeta +\ii\begin{bmatrix}0 & -\mathcal{U}\\\mathcal{V} & 0\end{bmatrix}.
\end{gather*}
These equations constitute the Lax pair for Painlev\'e-II found by Jimbo and Miwa~\cite{JimboM81}. The existence of a simultaneous fundamental solution matrix of both differential equations implies compatibility of the Lax pair, i.e., the coefficient matrices $\mathbf{A}$ and~$\mathbf{B}$ necessarily satisfy the zero-curvature condition $\mathbf{A}_y-\mathbf{B}_\zeta+[\mathbf{A},\mathbf{B}]=\mathbf{0}$, which by direct computation is equivalent to the following coupled system for the functions $\mathcal{U}=\mathcal{U}(y)$ and $\mathcal{V}=\mathcal{V}(y)$:
\begin{gather}
\frac{\dd^2\mathcal{U}}{\dd y^2}-\frac{1}{3}y\mathcal{U}-2\mathcal{U}^2\mathcal{V}=0\qquad\text{and}\qquad
\frac{\dd^2\mathcal{V}}{\dd y^2}-\frac{1}{3}y\mathcal{V}-2\mathcal{U}\mathcal{V}^2=0.\label{eq:PII-system}
\end{gather}
Combining the diagonal part of \eqref{eq:Z-zeta-Z-inverse-zeta-m-1} with the off-diagonal part of \eqref{eq:Z-y-Z-inverse-zeta-m-1} and using the identity $W^1_{11}(y)+W^1_{22}(y)=0$ following from $\det(\mathbf{Z}(\zeta;y))=1$ yields
\begin{gather}
\frac{1}{3}\ii p = \mathcal{U}'(y)\mathcal{V}(y)-\mathcal{U}(y)\mathcal{V}'(y).\label{eq:PII-COM}
\end{gather}
Since the left-hand side is independent of $y$, so must be the right-hand side, and indeed it is straightforward to check that $\mathcal{U}'(y)\mathcal{V}(y)-\mathcal{U}(y)\mathcal{V}'(y)$ is a constant of motion for the system \eqref{eq:PII-system}. It follows immediately that the logarithmic derivatives
\begin{gather}
\mathcal{P}(y):=\frac{\mathcal{U}'(y)}{\mathcal{U}(y)}\qquad\text{and}\qquad \mathcal{Q}(y):=\frac{\mathcal{V}'(y)}{\mathcal{V}(y)}\label{eq:PQ-define}
\end{gather}
satisfy uncoupled inhomogeneous Painlev\'e-II equations:
\begin{gather}
\frac{\dd^2\mathcal{P}}{\dd y^2}+\frac{2}{3}y\mathcal{P}-2\mathcal{P}^3+\frac{2}{3}\ii p-\frac{1}{3}=0\qquad\text{and}\qquad
\frac{\dd^2\mathcal{Q}}{\dd y^2}+\frac{2}{3}y\mathcal{Q}-2\mathcal{Q}^3-\frac{2}{3}\ii p-\frac{1}{3}=0.\label{eq:two-PIIs}
\end{gather}
These two equations are simply rescaled forms of the standard Painlev\'e-II equation \eqref{eq:PII-intro}. Indeed, the function $u(x):=-\big(\tfrac{3}{2}\big)^{1/3}\mathcal{P}\big({-}\big(\tfrac{3}{2}\big)^{1/3}x\big)$ is a solution of~\eqref{eq:PII-intro} with parameter $\alpha=\tfrac{1}{2}-\ii p$. Likewise, the function $u(x):=-\big(\tfrac{3}{2}\big)^{1/3}\mathcal{Q}\big({-}\big(\tfrac{3}{2}\big)^{1/3}x\big)$ satisfies \eqref{eq:PII-intro} with $\alpha=\tfrac{1}{2}+\ii p$.

\begin{Remark}Let $\mathbf{W}(\zeta;y,p,\tau)$ denote the solution of Riemann--Hilbert Problem~\ref{rhp:PII-generalized} with parameters $(p,\tau)$ related by $\tau^2=\ee^{2\pi p}-1$. Observing that $(p,-\tau)$ also satisfies the same relation, it is a direct matter to check that $\mathbf{W}(\zeta;y,p,-\tau)=(\ii\sigma_3)\mathbf{W}(\zeta;y,p,\tau)(\ii\sigma_3)^{-1}$. At the level of the functions~$\mathcal{U}$, $\mathcal{V}$, $\mathcal{P}$, and $\mathcal{Q}$, this symmetry implies that $\mathcal{U}(y;p,-\tau)=-\mathcal{U}(y;p,\tau)$, $\mathcal{V}(y;p,-\tau)=-\mathcal{V}(y;p,\tau)$, $\mathcal{P}(y;p,-\tau)=\mathcal{P}(y;p,\tau)$, and $\mathcal{Q}(y;p,-\tau)=\mathcal{Q}(y;p,\tau)$. The latter two identities indicate that, for the purposes of studying the solutions of~\eqref{eq:two-PIIs}, or equivalently of~\eqref{eq:PII-intro} for $\alpha=\tfrac{1}{2}\pm\ii p$, the sign of the second parameter $\tau$ in Riemann--Hilbert Problem~\ref{rhp:PII-generalized} is arbitrary and can be chosen for convenience.
\end{Remark}

\subsection{Asymptotic behavior for large $y$ and solution identification} When $p\in\mathbb{C}\setminus\ii\mathbb{Z}$, Riemann--Hilbert Problem~\ref{rhp:PII-generalized} determines a specific solution of each of the two Painlev\'e-II equations \eqref{eq:two-PIIs} as well as specific corresponding logarithmic ``potentials'' $\mathcal{U}$ and $\mathcal{V}$ (which would involve additional integration constants to determine from $\mathcal{P}$ and $\mathcal{Q}$). The monograph \cite{FokasIKN06} contains an exhaustive description of all solutions of the standard form Painlev\'e-II equation \eqref{eq:PII-intro} in which their properties are associated with the monodromy data of a different Riemann--Hilbert problem, namely that corresponding to the alternate Lax pair discovered by Flaschka and Newell \cite{FlaschkaN80}. Since to our knowledge the literature does not yet contain a complete description of the relation between the monodromy data for the Lax pairs of Flaschka--Newell and Jimbo--Miwa, in order to identify the specific solutions arising we need to compute the large-$y$ asymptotic behavior directly from Riemann--Hilbert Problem~\ref{rhp:PII-generalized} and then compare with known asymptotics of solutions catalogued in \cite{FokasIKN06}.

Let us write $y\in\mathbb{C}$, $y\neq 0$, uniquely in the form $y=M\ee^{\ii\theta}$, $M=|y|$. Letting $\zeta = M^{1/2}\xi$ and setting
\begin{gather*}
\mathbf{Y}\big(\xi;M,\ee^{\ii\theta}\big):=M^{\ii p\sigma_3/2}\mathbf{W}\big(M^{1/2}\xi;M\ee^{\ii\theta}\big),
\end{gather*}
the conditions of Riemann--Hilbert Problem~\ref{rhp:PII-generalized} imply an equivalent rescaled Riemann--Hilbert problem for $\mathbf{Y}\big(\xi;M,\ee^{\ii\theta}\big)$:
\begin{rhp}[rescaled Painlev\'e-II parametrix]\label{rhp:PII-rescaled} Let $p\in\mathbb{C}\setminus\ii\mathbb{Z}$ and $\tau=\pm\sqrt{\ee^{2\pi p}-1}$ be given parameters. Given also a number $\ee^{\ii\theta}$ on the unit circle and $M>0$, seek a~$2\times 2$ matrix-valued function $\mathbf{Y}\big(\xi;M,\ee^{\ii\theta}\big)$ with the following properties.
\begin{itemize}\itemsep=0pt
\item[]\textbf{Analyticity:} $\mathbf{Y}\big(\xi;M,\ee^{\ii\theta}\big)$ is analytic for $\xi$ in the five sectors shown in Fig.~{\rm \ref{fig:PII-jumps-generalized}}, now interpreted in the $\xi$-plane. It takes continuous boundary values on the excluded rays and at the origin from each sector.
\item[]\textbf{Jump conditions:} $\mathbf{Y}_+\big(\xi;M,\ee^{\ii\theta}\big) =\mathbf{Y}_-\big(\xi;M,\ee^{\ii\theta}\big)\mathbf{V}\big(\xi;M,\ee^{\ii\theta}\big)$, where $\mathbf{V}\big(\xi;M,\ee^{\ii\theta}\big)$ is the same jump matrix as shown in Fig.~{\rm \ref{fig:PII-jumps-generalized}} except that the exponent $\zeta^3+y\zeta$ is everywhere replaced with $M^{3/2}\varphi\big(\xi;\ee^{\ii\theta}\big)$, $\varphi\big(\xi;\ee^{\ii\theta}\big):=\xi^3+\ee^{\ii\theta}\xi$.
\item[]\textbf{Normalization:} $\mathbf{Y}\big(\xi;M,\ee^{\ii\theta}\big)\xi^{\ii p\sigma_3}\to\mathbb{I}$ as $\xi\to\infty$ uniformly in all directions.
\end{itemize}
\end{rhp}

\subsubsection[The case $\ee^{\ii\theta}=-1$. Asymptotic behavior as $y\to -\infty$]{The case $\boldsymbol{\ee^{\ii\theta}=-1}$. Asymptotic behavior as $\boldsymbol{y\to -\infty}$}\label{sec:minus-one}

\begin{figure}[t]\centering
\includegraphics{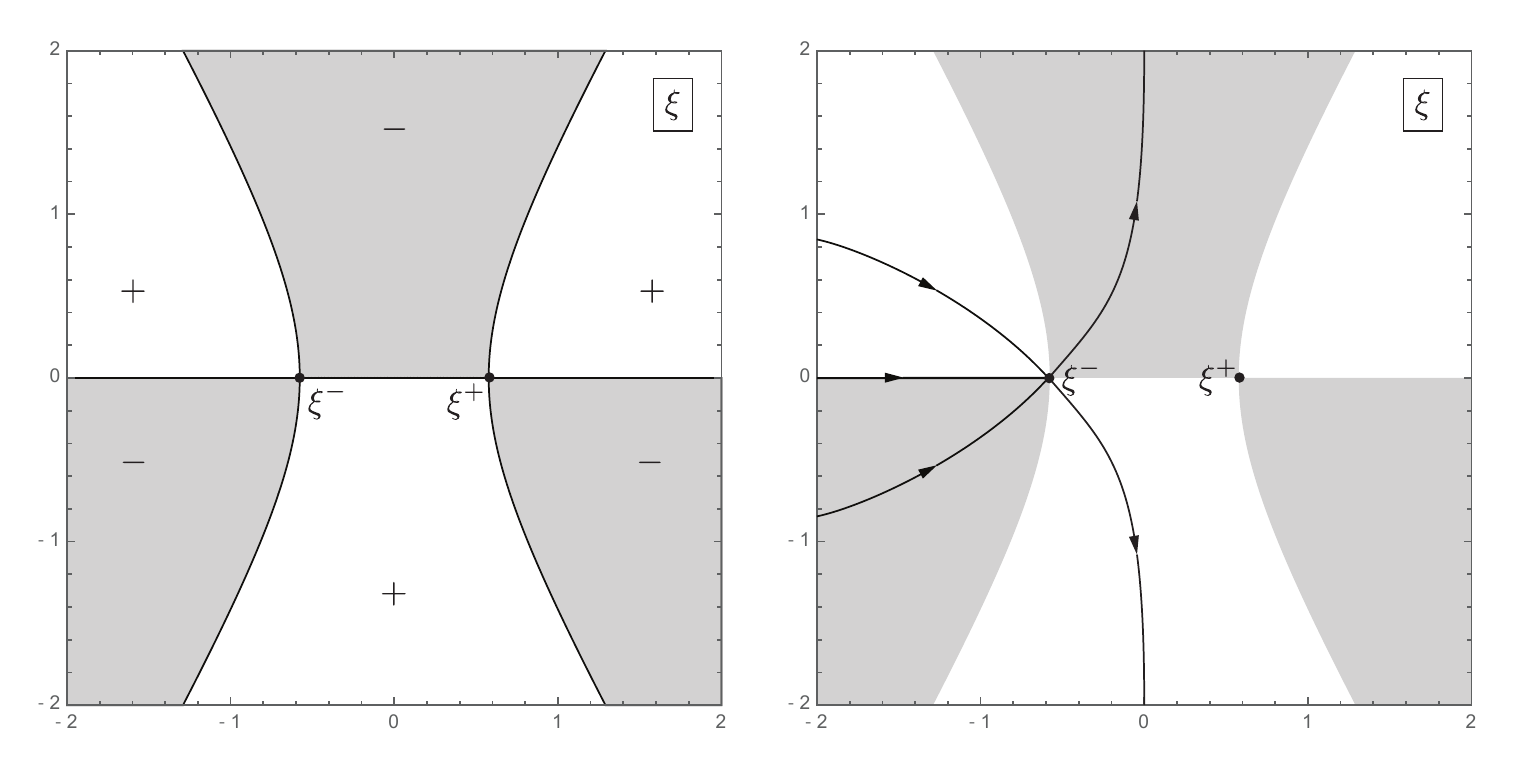}
\caption{Left: the sign chart for $\mathrm{Im}(\varphi(\xi;-1))$. Right: the deformed jump contour $\Sigma_\mathbf{Y}$ for $\mathbf{Y}(\xi;M,-1)$.}\label{fig:PII-Minus}
\end{figure}

Note that $\varphi(\xi;-1)$ has critical points $\xi=\xi^\pm:=\pm 1/\sqrt{3}$. We observe that the sign table for $\mathrm{Im}(\varphi(\xi;-1))$ has the structure plotted in the left-hand panel of Fig.~\ref{fig:PII-Minus}.
We take the jump contour for $\mathbf{Y}(\xi;M,-1)$ to be deformed as shown in the right-hand panel of Fig.~\ref{fig:PII-Minus}, so that in particular the self-intersection point now coincides with the left-most critical point $\xi=\xi^-$. As an outer parametrix for $\mathbf{Y}(\xi;M,-1)$ we take the matrix function
\begin{gather}
\dot{\mathbf{Y}}^{\mathrm{out}}(\xi):=(\xi-\xi^-)^{-\ii p\sigma_3},\label{eq:Y-out-PII-Minus}
\end{gather}
which satisfies exactly the jump condition of $\mathbf{Y}(\xi;M,-1)$ on the ray $\xi<\xi^-$, as well as the normalization condition on $\mathbf{Y}(\xi;M,-1)$. Let $D(\xi^-)$ denote a disk centered at $\xi=\xi^-$ of radius less than $2/\sqrt{3}$ (so that it excludes the other critical point $\xi^+$). For $\xi\in D(\xi^-)$, we introduce a~conformal mapping $\xi\mapsto Z(\xi)$ that satisfies the equation
\begin{gather}
\varphi(\xi^-;-1)-\varphi(\xi;-1)=Z^2,\label{eq:PII-Minus-conformal}
\end{gather}
which has a unique analytic solution $Z=Z(\xi)$ for which $Z(\xi^-)=0$ and $Z'(\xi^-)>0$. Assuming that when $\xi\in D(\xi^-)$ the jump contour for $\mathbf{Y}(\xi;M,-1)$ is taken to consist of five straight-line segments joined at the origin in the $Z$-plane with $\arg(Z)=\pm\tfrac{1}{4}\pi$, $\arg(Z)=\pm\tfrac{3}{4}\pi$, and $\arg(-Z)=0$, we observe that the matrix $\mathbf{P}\big(M^{3/4}Z(\xi);p,\tau\big)\ee^{-\ii M^{3/2}\varphi(\xi^-;-1)\sigma_3/2}$, where $\mathbf{P}(\lambda;p,\tau)$ is the solution of the parabolic cylinder Riemann--Hilbert problem solved in the appendix, is an exact local solution of Riemann--Hilbert Problem~\ref{rhp:PII-rescaled} near $\xi^-$. Writing the outer parametrix in the form
\begin{gather*}
\dot{\mathbf{Y}}^\mathrm{out}(\xi)=M^{3\ii p\sigma_3/4}\mathbf{H}(\xi)(M^{3/4}Z(\xi))^{-\ii p\sigma_3},\qquad \mathbf{H}(\xi):=\left(\frac{\xi-\xi^-}{Z(\xi)}\right)^{-\ii p\sigma_3},
\end{gather*}
where we note that $\mathbf{H}(\xi)$ is analytic within $D(\xi^-)$ and is independent of $M$, we define the inner parametrix by the formula
\begin{gather*}
\dot{\mathbf{Y}}^\mathrm{in}(\xi;M):=M^{3\ii p\sigma_3/4}\ee^{\ii M^{3/2}\varphi(\xi^-;-1)\sigma_3/2}\mathbf{H}(\xi)\mathbf{P}\big(M^{3/4}Z(\xi);p,\tau\big)\\
\hphantom{\dot{\mathbf{Y}}^\mathrm{in}(\xi;M):=}{}\times \ee^{-\ii M^{3/2}\varphi(\xi^-;-1)\sigma_3/2},\qquad \xi\in D(\xi^-).
\end{gather*}
We combine the outer and inner parametrices into a \emph{global parametrix} $\dot{\mathbf{Y}}(\xi;M)$ defined by
\begin{gather*}
\dot{\mathbf{Y}}(\xi;M):=\begin{cases}
\dot{\mathbf{Y}}^\mathrm{in}(\xi;M),& \xi\in D(\xi^-),\\
\dot{\mathbf{Y}}^\mathrm{out}(\xi),& \xi\in\mathbb{C}\setminus\overline{D(\xi^-)}.
\end{cases}
\end{gather*}

Consider the error matrix defined by
\begin{gather}
\mathbf{F}(\xi;M):=M^{-3\ii p\sigma_3/4}\ee^{-\ii M^{3/2}\varphi(\xi^-;-1)\sigma_3/2}\mathbf{Y}(\xi;M,-1)\dot{\mathbf{Y}}(\xi;M)^{-1}\nonumber\\
\hphantom{\mathbf{F}(\xi;M):=}{}\times \ee^{\ii M^{3/2}\varphi(\xi^-;-1)\sigma_3/2}M^{3\ii p\sigma_3/4}.\label{eq:F-def-PII}
\end{gather}
Because $\mathbf{Y}$ and its parametrix $\dot{\mathbf{Y}}$ satisfy exactly the same jump conditions both on the real axis to the left of $\xi^-$ and also on all jump contour arcs within the disk $D(\xi^-)$, it can be shown that $\mathbf{F}(\xi;M)$ admits analytic continuation to these contours and hence can be regarded as a~function analytic in the complex $\xi$-plane except on the four non-real jump contours shown in the right-hand panel of Fig.~\ref{fig:PII-Minus} restricted to the exterior of the disk, and the disk boundary $\partial D(\xi^-)$ where the inner and outer parametrices fail to match exactly. The jump matrix for $\mathbf{F}(\xi;M)$ on the resulting jump contour just described is an exponentially small perturbation of the identity matrix except when $\xi\in\partial D(\xi^-)$ due to the placement of the contour relative to the sign chart of $\mathrm{Im}(\varphi(\xi;-1))$ as shown in Fig.~\ref{fig:PII-Minus}. Taking the circle $\partial D(\xi^-)$ to have clockwise orientation, the jump condition for $\mathbf{F}(\xi;M)$ expresses the mismatch between outer and inner parametrices in the form $\mathbf{F}_+(\xi;M)=\mathbf{F}_-(\xi;M)\mathbf{V}^\mathbf{F}(\xi;M)$, where
\begin{gather}
\mathbf{V}^\mathbf{F}(\xi;M) =M^{-3\ii p\sigma_3/4}\ee^{-\ii M^{3/2}\varphi(\xi^-;-1)\sigma_3/2}\dot{\mathbf{Y}}^\mathrm{in}(\xi;M)\dot{\mathbf{Y}}^\mathrm{out}(\xi)^{-1}\ee^{\ii M^{3/2}\varphi(\xi^-;-1)\sigma_3/2}M^{3\ii p\sigma_3/4}\nonumber\\
\hphantom{\mathbf{V}^\mathbf{F}(\xi;M)}{} =\mathbf{H}(\xi) \mathbf{P}(\lambda;p,\tau)\lambda^{\ii p\sigma_3}\mathbf{H}(\xi)^{-1},\qquad\xi\in\partial D(\xi^-),
\label{eq:F-jump-circle-PII-Minus}
\end{gather}
where $\lambda=M^{3/4}Z(\xi)$. Since the conjugating factors are bounded as $M\to\infty$ while $\lambda$ is bounded below by a multiple of $M^{3/4}$ when $\xi\in D(\xi^-)$, it follows that $\mathbf{V}^\mathbf{F}(\xi;M)$ differs from the identity uniformly on the jump contour for $\mathbf{F}(\xi;M)$ by $O\big(M^{-3/4}\big)$. The boundary value $\mathbf{F}_-(\xi;M)$ necessarily satisfies the integral equation
\begin{gather}
\mathbf{F}_-(\xi;M)-\mathbb{I}=\mathcal{C}_-^{\Sigma_\mathbf{F}}\big(\mathbf{V}^\mathbf{F}(\cdot;M)-\mathbb{I}\big)(\xi) \nonumber\\
\hphantom{\mathbf{F}_-(\xi;M)-\mathbb{I}=}{} +\mathcal{C}_-^{\Sigma_\mathbf{F}}\big((\mathbf{F}_-(\cdot;M)-\mathbb{I})\big(\mathbf{V}^\mathbf{F}(\cdot;M)-\mathbb{I}\big)\big)(\xi),\qquad \xi\in\Sigma_\mathbf{F},\label{eq:singular-integral-eqn}
\end{gather}
which is to be solved for $\mathbf{F}_-(\cdot;M)-\mathbb{I}\in L^2(\Sigma_\mathbf{F})$. Here, for a contour $\Sigma$, $\mathcal{C}^\Sigma_-$ denotes the right boundary value Cauchy operator defined by
\begin{gather*}
\mathcal{C}^\Sigma_-(\mathbf{F}(\cdot))(\xi):=\frac{1}{2\pi\ii}\int_\Sigma\frac{\mathbf{F}(z)\,\dd z}{z-\xi_-},\qquad \xi\in\Sigma\setminus\{\text{self-intersection points}\},
\end{gather*}
where the subscript ``$-$'' indicates that a nontangential limit is taken toward $\xi$ from the right side of $\Sigma$ by local orientation. It is well-known to be a bounded operator on $L^2(\Sigma)$ for contours such as $\Sigma^\mathbf{F}$, with a norm depending only on geometrical details of the contour. Now due to the exponential decay of $\mathbf{V}^\mathbf{F}(\xi;M)$ as $\xi\to\infty$ along the four unbounded rays of $\Sigma_\mathbf{F}$, it is easy to see that not only $\|\mathbf{V}^\mathbf{F}(\cdot;M)-\mathbb{I}\|_\infty=O\big(M^{-3/4}\big)$ (uniform estimate on $\Sigma_\mathbf{F}$) but also $\big\|\mathbf{V}^\mathbf{F}(\cdot;M)-\mathbb{I}\big\|_2=O\big(M^{-3/4}\big)$ ($L^2$ estimate on $\Sigma_\mathbf{F}$). It follows easily that the integral equation is uniquely solvable in $L^2(\Sigma_\mathbf{F})$ with solution satisfying $\|\mathbf{F}_-(\cdot;M)-\mathbb{I}\|_2=O\big(M^{-3/4}\big)$ as $M\to\infty$. From the integral representation
\begin{gather*}
\mathbf{F}(\xi;M)=\mathbb{I}+\frac{1}{2\pi\ii}\int_{\Sigma_\mathbf{F}}\frac{\mathbf{V}^\mathbf{F}(w;M)-\mathbb{I}}{w-\xi}\,\dd w \\
\hphantom{\mathbf{F}(\xi;M)=}{} + \frac{1}{2\pi\ii}\int_{\Sigma_\mathbf{F}}\frac{(\mathbf{F}_-(w;M)-\mathbb{I})\big(\mathbf{V}^\mathbf{F}(w;M)-\mathbb{I}\big)}{w-\xi}\,\dd w,\qquad \xi\in\mathbb{C}\setminus\Sigma_\mathbf{F},
\end{gather*}
the exponential decay of $\mathbf{V}^\mathbf{F}(\xi;M)-\mathbb{I}$ as $\xi\to\infty$ in $\Sigma_\mathbf{F}$ guarantees that
\begin{gather}
\mathbf{F}(\xi;M)=\mathbb{I}+\mathbf{F}^1(M)\xi^{-1}+\mathbf{F}^2(M)\xi^{-2}+O\big(\xi^{-3}\big),\qquad \xi\to\infty, \label{eq:F-PII-Minus-Expansion}
\end{gather}
where
\begin{gather*}
\mathbf{F}^k(M):=-\frac{1}{2\pi\ii}\int_{\Sigma_\mathbf{F}}\big(\mathbf{V}^\mathbf{F}(w;M)-\mathbb{I}\big)w^{k-1}\,\dd w \\
\hphantom{\mathbf{F}^k(M):=}{}-\frac{1}{2\pi\ii}\int_{\Sigma_\mathbf{F}}(\mathbf{F}_-(w;M)-\mathbb{I})\big(\mathbf{V}^\mathbf{F}(w;M) -\mathbb{I}\big)w^{k-1}\,\dd w,\qquad k=1,2.
\end{gather*}
Since $\xi\big(\mathbf{V}^\mathbf{F}(\xi;M)-\mathbb{I}\big)$ is also in $L^2(\Sigma_\mathbf{F})$ with norm proportional to $M^{-3/4}$, it follows by Cauchy--Schwarz that
\begin{gather}
\mathbf{F}^k(M)=-\frac{1}{2\pi\ii}\int_{\Sigma_\mathbf{F}}\big(\mathbf{V}^\mathbf{F}(w;M)-\mathbb{I}\big)w^{k-1}\,\dd w + O\big(M^{-3/2}\big),\nonumber\\
k=1,2,\qquad M\to\infty. \label{eq:F-PII-Minus-moments}
\end{gather}
The contribution to the first term from integration over $\Sigma_\mathbf{F}\setminus\partial D(\xi^-)$ is exponentially small as $M\to\infty$, and therefore we may take the integration over just $\partial D(\xi^-)$ with no change in the order of the error estimate. Combining \eqref{eq:F-PII-Minus-Expansion} with~\eqref{eq:F-def-PII} and the identity $\dot{\mathbf{Y}}(\xi;M)=\dot{\mathbf{Y}}^\mathrm{out}(\xi)$ holding for sufficiently large $|\xi|$, recalling the definition~\eqref{eq:Y-out-PII-Minus} together with the relation $\mathbf{Y}(\xi;M,-1)=M^{\ii p\sigma_3/2}\mathbf{W}\big(M^{1/2}\xi;-M\big)$ and the first expansion in~\eqref{eq:large-zeta-PII}, we see that for sufficiently large negative~$y$,
\begin{gather*}
\mathcal{V}(y)=W^1_{21}(y)=M^{-\ii p/2}\sqrt{M}\ee^{-\ii M^{3/2}\varphi(\xi^-;-1)}F_{21}^1(M), \nonumber\\
\mathcal{Q}(y)=\frac{\mathcal{V}'(y)}{\mathcal{V}(y)}=\ii W^1_{11}(y)-\ii\frac{W^2_{21}(y)}{W^1_{21}(y)}=
\ii \sqrt{M}\left(F^1_{11}(M)-\frac{F_{21}^2(M)}{F_{21}^1(M)}\right), %\label{eq:VQ-Minus}
\end{gather*}
where $M=-y$ for $\ee^{\ii\theta}=-1$. In the case of the formula for $\mathcal{Q}(y)$ we also used~\eqref{eq:Z-y-Z-inverse-zeta-m-1} to eliminate~$\mathcal{V}'(y)$. To calculate $\mathcal{V}(y)$, we therefore substitute~\eqref{eq:F-jump-circle-PII-Minus} into \eqref{eq:F-PII-Minus-moments} for $k=1$, replacing the integration contour by $\partial D(\xi^-)$ with clockwise orientation, and use also the expansion \eqref{eq:U-matrix-asymp-generalized} from the appendix to obtain:
\begin{gather*}
\mathcal{V}(y)=\frac{s(p,\tau)}{2\ii}M^{-\ii p/2-1/4}\ee^{-\ii M^{3/2}\varphi(\xi^-;-1)}\nonumber\\
\hphantom{\mathcal{V}(y)=}{}\times \left[\frac{1}{2\pi\ii}\int_{\partial D(\xi^-)}\left(\frac{w-\xi^-}{Z(w)}\right)^{2\ii p}\frac{\dd w}{Z(w)} + O\big(M^{-5/4}\big)\right],\!\!\qquad M\to +\infty,\qquad\!\! y<0,\!\!\!\! %\label{eq:V-Minus}
\end{gather*}
where $s(p,\tau)\neq 0$ is given by \eqref{eq:beta-def}. Since the integrand is analytic within the circle of integration aside from a simple pole at $w=\xi^-$, we compute it by residues and obtain
\begin{gather}
\mathcal{V}(y)=-\frac{s(p,\tau)}{2\ii}M^{-\ii p/2-1/4}\ee^{-\ii M^{3/2}\varphi(\xi^-;-1)}\nonumber\\
\hphantom{\mathcal{V}(y)=}{} \times \big[Z'(\xi^-)^{-2\ii p-1} + O\big(M^{-5/4}\big)\big],\qquad M\to +\infty,\qquad y<0. \label{eq:V-asymp-M-Minus}
\end{gather}
Note that $\varphi(\xi^-;-1)=2\cdot 3^{-3/2}$ and taking two derivatives of~\eqref{eq:PII-Minus-conformal} with $Z(\xi^-)=0$ and $Z'(\xi^-)>0$ gives $Z'(\xi^-)=3^{1/4}$. Therefore, also using~\eqref{eq:beta-def} to eliminate $s(p,\tau)$, we arrive at the asymptotic formula~\eqref{eq:V-asymp-y-Minus}, which is valid for all $p\in\mathbb{C}\setminus\ii\mathbb{Z}$.

Now the dominant contribution to $F_{11}^1(M)$ is $O\big(M^{-3/2}\big)$ as $M\to\infty$. On the other hand, the dominant contribution to $F_{21}^2(M)$ is calculated exactly as above, by residues. Since the integrand now has an additional factor of $w$, it is easy to see that $F_{21}^2(M)/F_{21}^1(M)=\xi^- + O\big(M^{-3/2}\big)$. We therefore conclude that
\begin{gather*}
\mathcal{Q}(y)=-\ii\xi^-\sqrt{-y}+O\big(y^{-1}\big)=\ii\sqrt{-\frac{y}{3}}+O\big(y^{-1}\big),\qquad y\to -\infty.%\label{eq:Q-Minus}
\end{gather*}

\subsubsection[Generalization to complex $y$]{Generalization to complex $\boldsymbol{y}$}
\begin{figure}[t]\centering
\includegraphics{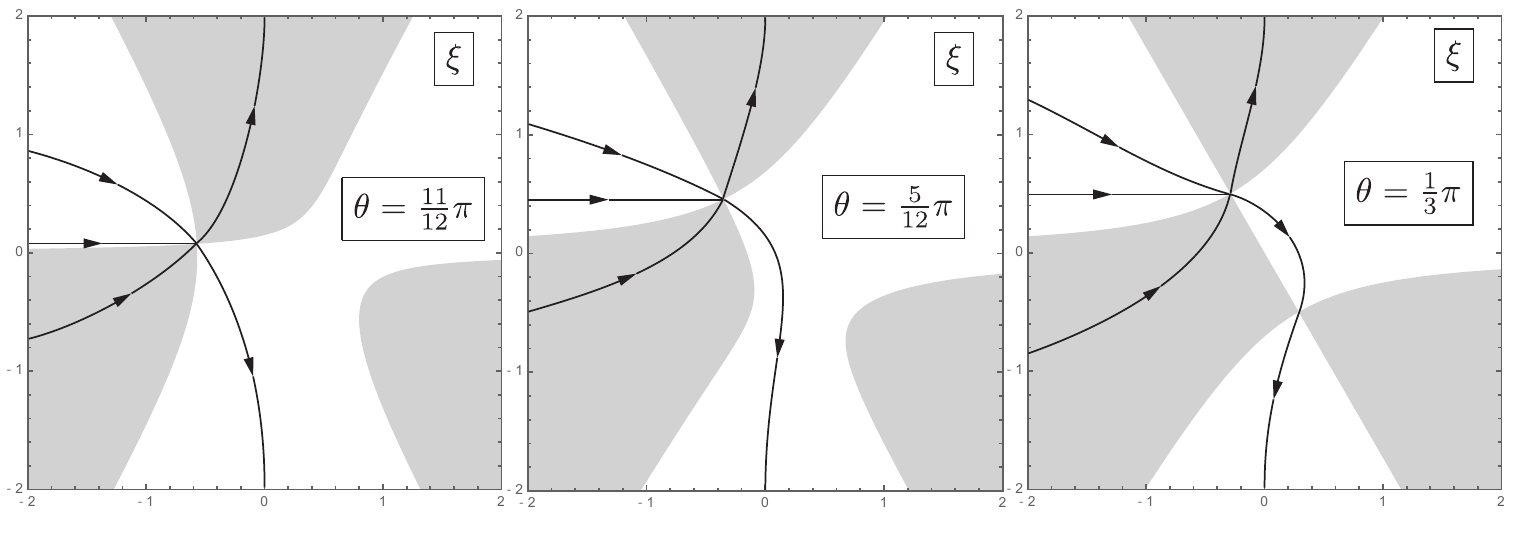}
\caption{Placement of the jump contour $\Sigma_\mathbf{Y}$ relative to the sign chart for $\mathrm{Im}(\varphi(\xi;\ee^{\ii\theta}))=\mathrm{Im}(\varphi(\xi^-;\ee^{\ii\theta}))$ for $\theta=\tfrac{11}{12}\pi$ (left panel), $\theta=\tfrac{5}{12}\pi$ (central panel), and $\theta=\tfrac{1}{3}\pi$ (right panel).}\label{fig:PII-UHP}
\end{figure}

\begin{figure}[t]\centering
\includegraphics{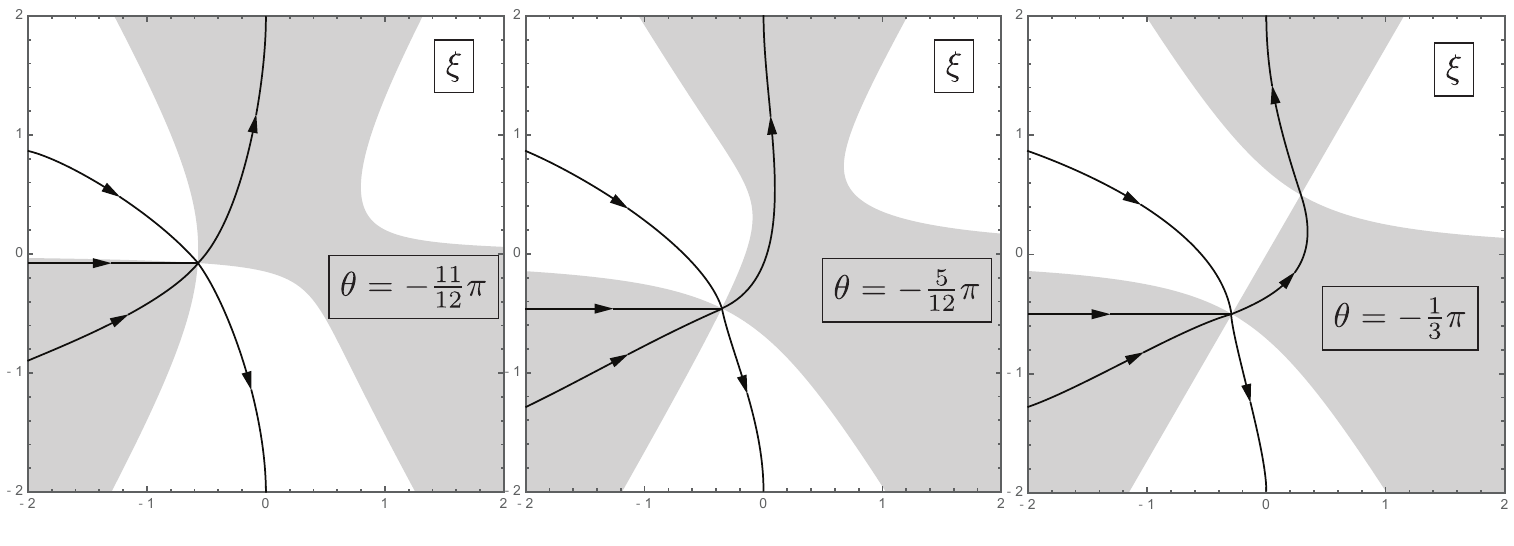}
\caption{As in Fig.~\ref{fig:PII-UHP}, but for $\theta=-\tfrac{11}{12}\pi$ (left panel), $\theta=-\tfrac{5}{12}\pi$ (central panel), and $\theta=-\tfrac{1}{3}\pi$ (right panel).}\label{fig:PII-LHP}
\end{figure}

As soon as $\ee^{\ii\theta}\neq -1$, the two critical points of $\varphi\big(\xi;\ee^{\ii\theta}\big)$ are in general no longer on the same level of $\mathrm{Im}\big(\varphi\big(\xi;\ee^{\ii\theta}\big)\big)$. We denote by $\xi^-=-\big({-}\tfrac{1}{3}\ee^{\ii\theta}\big)^{1/2}$ the critical point that agrees with $-1/\sqrt{3}$ when $\ee^{\ii\theta}=-1$. The level set $\mathrm{Im}\big(\varphi\big(\xi;\ee^{\ii\theta}\big)\big)=\mathrm{Im}\big(\varphi\big(\xi^-;\ee^{\ii\theta}\big)\big)$ undergoes a bifurcation in the neighborhood of the other critical point $\xi=-\xi^-$, when $\ee^{\ii\theta}=-1$ and then again when $\ee^{\ii\theta}=\ee^{\pm\ii\pi/3}$. The bifurcation at $\ee^{\ii\theta}=-1$ is harmless from the point of view of placing the jump contour so as to achieve exponential decay of the jump matrix to the identity relative to the neighborhood of the point $\xi^-$. However, the bifurcations at $\ee^{\ii\theta}=\ee^{\pm\ii\pi/3}$ are genuine obstructions to the basic approach valid for $\ee^{\ii\theta}=-1$. The development of the problematic bifurcation as $\ee^{\ii\theta}$ ranges over the upper (lower) half semicircle from $-1$ toward $\ee^{\ii\pi/3}$ (toward $\ee^{-\ii\pi/3}$) is illustrated in Fig.~\ref{fig:PII-UHP} (in Fig.~\ref{fig:PII-LHP}).

These figures show that the same basic method as applies to the special case of $\ee^{\ii\theta}=-1$ also applies over the whole range $\arg(-y)=\arg\big({-}\ee^{\ii\theta}\big)\in \big({-}\tfrac{2}{3}\pi,\tfrac{2}{3}\pi\big)$. \emph{Mutatis mutandis}, the construction of the parametrix explained in Section~\ref{sec:minus-one} is the same, as is the analysis of the error $\mathbf{F}(\xi;M)$, and we arrive again at the asymptotic formula~\eqref{eq:V-asymp-M-Minus} for $\mathcal{V}(y)$ now valid for $|\arg(-y)|<\tfrac{2}{3}\pi$, except that $\varphi(\xi^-;-1)$ has to be replaced with $\varphi\big(\xi^-;\ee^{\ii\theta}\big)$ in the exponent and that the conformal mapping~$Z(\xi)$ has a generalized interpretation. In general, $Z(\xi)$ is the conformal mapping defined near $\xi^-$ by $\varphi\big(\xi^-;\ee^{\ii\theta}\big) -\varphi\big(\xi;\ee^{\ii\theta}\big)=Z(\xi)^2$ that is chosen to depend continuously on $\ee^{\ii\theta}$. Therefore, $Z'(\xi^-)=3^{1/4}\big({-}\ee^{\ii\theta}\big)^{1/4}$, and also $\varphi\big(\xi^-;\ee^{\ii\theta}\big)=2\cdot 3^{-3/2}\big({-}\ee^{\ii\theta}\big)^{3/2}$, so we see that \eqref{eq:V-asymp-y-Minus} holds true as $y\to\infty$ with $|\arg(-y)|<\tfrac{2}{3}\pi$. Similarly,
\begin{gather}
\mathcal{Q}(y)=\ii\left(-\frac{y}{3}\right)^{1/2}+O\big(|y|^{-1}\big),\qquad y\to\infty,\qquad |\arg(-y)|<\frac{2}{3}\pi.\label{eq:Q-y-wide-asymp}
\end{gather}
In both cases, we use the principal branch to interpret the power functions: $(-y)^P=\ee^{P\log(-y)}$ with $|\mathrm{Im}(\log(-y))|<\pi$.

\subsubsection{Solution identification}
Recall that the function $u(x):=-\big(\tfrac{3}{2}\big)^{1/3}\mathcal{Q}\big({-}\big(\tfrac{3}{2}\big)^{1/3}x\big)$ is a solution of the Painlev\'e-II equation in standard form \eqref{eq:PII-intro} with parameter $\alpha=\tfrac{1}{2}+\ii p$. From \eqref{eq:Q-y-wide-asymp}, we therefore have shown that $u(x)$ matches the asymptotic description of the solution $u_\mathrm{TT}^-\big(x;\tfrac{1}{2}+\ii p\big)$ given in~\eqref{eq:u-x-wide-asymp}, and covering the sector $|\arg(x)|<\tfrac{2}{3}\pi$. According to the analysis of the Riemann--Hilbert problem arising from the Flaschka--Newell Lax pair of Painlev\'e-II described in \cite[Chapter~11]{FokasIKN06}, for each $\alpha\in\mathbb{C}\setminus\big(\mathbb{Z}+\tfrac{1}{2}\big)$ there is \emph{exactly one} solution of \eqref{eq:PII-intro} consistent with the asymptotic formula~\eqref{eq:u-x-wide-asymp} (for each choice of the sign $\pm$). Thus $u(x)=u_\mathrm{TT}^-\big(x;\tfrac{1}{2}+\ii p\big)$. Having square-root asymptotics in such a large sector of the complex $x$-plane, the solutions $u^\pm_\mathrm{TT}(x;\alpha)$ are two of the six increasing tritronqu\'ee solutions of~\eqref{eq:PII-intro} (four others are obtained by $\pm\tfrac{2}{3}\pi$ rotation symmetry in the complex $x$-plane). Each of these six is determined by its leading asymptotic behavior in a sector of maximal opening angle $\tfrac{4}{3}\pi-\epsilon$. It is shown in~\cite{FokasIKN06} that the leading asymptotic~\eqref{eq:u-x-wide-asymp} admits correction in the form of a full asymptotic series in decreasing powers of~$x$ differing by~$\tfrac{3}{2}$ in consecutive terms, the coefficients of which can be determined by formal substitution into the differential equation~\eqref{eq:PII-intro}. Thus in particular it holds that
\begin{gather*}
u(x)=-\ii\left(\frac{x}{2}\right)^{1/2}-\frac{\alpha}{2x}+O\big(|x|^{-5/2}\big),\qquad x\to\infty,\qquad |\arg(x)|<\frac{2}{3}\pi.
\end{gather*}
Using $\alpha=\tfrac{1}{2}+\ii p$, this implies also that \eqref{eq:Q-y-wide-asymp} can be improved to
\begin{gather}
\mathcal{Q}(y)=\ii\left(-\frac{y}{3}\right)^{\tfrac{1}{2}} -\left(\frac{1}{4}+\ii \frac{p}{2}\right)\frac{1}{y}+O\big(|y|^{-5/2}\big),\qquad y\to\infty,\qquad |\arg(-y)|<\frac{2}{3}\pi.\label{eq:Q-asymp-sector}
\end{gather}
The solution $u(x)$ characterized by the asymptotic behavior \eqref{eq:u-x-wide-asymp} is globally meromorphic and has poles near $x=\infty$ in the complementary sector $|\arg(-x)|<\tfrac{1}{3}\pi$; equivalently $\mathcal{Q}(y)$ has poles near $y=\infty$ for $|\arg(y)|<\tfrac{1}{3}\pi$. Indeed, for typical values of $\ee^{\ii\theta}$ with $|\theta|<\frac{1}{3}\pi$, the procedure of asymptotic analysis of Riemann--Hilbert Problem~\ref{rhp:PII-rescaled} requires the introduction of a two-cut $g$-function with one cut shrinking while approaching each of the two critical points of $\varphi\big(\xi;\ee^{\ii\theta}\big)$ as $|\theta|\uparrow \tfrac{1}{3}\pi$ (see the right-hand panels of Figs.~\ref{fig:PII-UHP} and~\ref{fig:PII-LHP}). Thus the asymptotic behavior is given in terms of a certain elliptic function of $M$ with modulus depending on $\theta$. Another degeneration occurs precisely when $\theta=0$ as the two cuts merge at the origin and thus become a single cut. This observation implies that the elliptic asymptotics give way to algebraic/trigonometric asymptotics in the limit $y\to +\infty$, and it forms the basis for the proof of Theorem~\ref{thm:centerline} which we turn to next.

\section{Proof of Theorem~\ref{thm:centerline}}\label{sec:centerline}
To study $\mathcal{V}(y)$ and $\mathcal{Q}(y)$ in the opposite limit $y\to +\infty$, we return to the consideration of $\mathbf{Y}\big(\xi;M,\ee^{\ii\theta}\big)$ as $M\to +\infty$, now in the case $\ee^{\ii\theta}=+1$. In this situation there is no distinction between $y$ and its modulus $M$, so we just replace $M$ with $y>0$ in this section.
\subsection[Introduction of $g$-function]{Introduction of $\boldsymbol{g}$-function}
\begin{figure}[t]\centering
\includegraphics{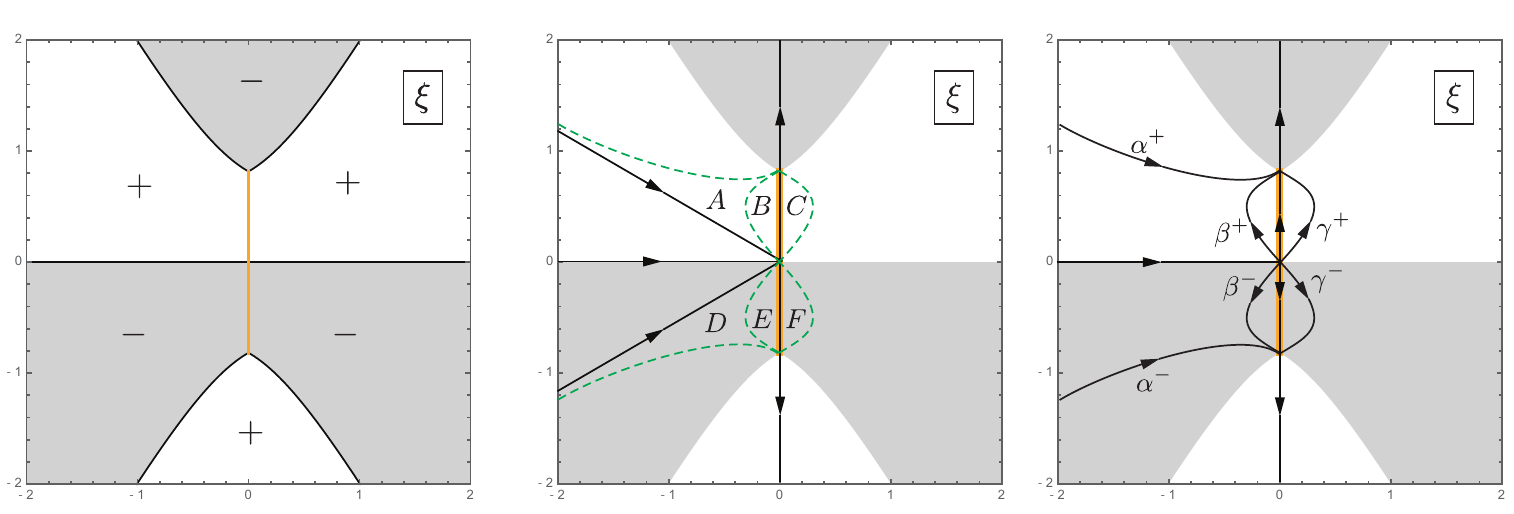}
\caption{Left: the sign chart for $\mathrm{Im}(h(\xi))$. Center: the original jump contour $\Sigma_\mathbf{Y}$ for $\mathbf{Y}(\xi;y,1)$ and the regions $A$, $B$, $C$, $D$, $E$, and $F$. Right: the jump contour $\Sigma_\mathbf{Z}$ for $\mathbf{Z}(\xi;y)$. The branch cut $\Sigma=\Sigma^+\cup\Sigma^-$ for $g$ and $h$ is shown in orange.}\label{fig:PII-Plus}
\end{figure}

The critical points of $\varphi(\xi;1)$ form a conjugate pair, and to deal with this it turns out to be necessary to introduce a so-called $g$-function. Define
\begin{gather}
\nu:=\sqrt{\frac{2}{3}},\label{eq:nu-define}
\end{gather}
let $\Sigma$ denote the straight line segment connecting the points $\pm \ii\nu$, and then set
\begin{gather}
g'(\xi):=\frac{3}{2}\xi R(\xi)-\frac{3}{2}\xi^2-\frac{1}{2}, \qquad\xi\in\mathbb{C}\setminus\Sigma,\label{eq:gprime}
\end{gather}
where $R(\xi)$ is the function analytic for $\xi\in\mathbb{C}\setminus\Sigma$ determined from the conditions $R(\xi)^2=\xi^2+\nu^2$ while $R(\xi)=\xi+O\big(\xi^{-1}\big)$ as $\xi\to\infty$. It is easy to confirm directly that $g'(\xi)=O\big(\xi^{-2}\big)$ as $\xi\to\infty$ so $g(\xi)$ is well defined by the integral
\begin{gather*}
g(\xi):=\int_\infty^\xi g'(s)\,\dd s,\qquad \xi\in\mathbb{C}\setminus\Sigma,
\end{gather*}
where the path of integration is arbitrary in the domain of analyticity of the integrand. It is easy to obtain the asymptotic formula
\begin{gather}
g(\xi)=\frac{1}{12\xi}+O\big(\xi^{-3}\big),\qquad\xi\to\infty.\label{eq:g-asymp}
\end{gather}
Consider the function $h(\xi):=g(\xi)+\tfrac{1}{2}\varphi(\xi;1)$. The left-hand panel of Fig.~\ref{fig:PII-Plus} shows the sign chart for $\mathrm{Im}(h(\xi))$, a function that is continuous across the branch cut $\Sigma$, which in turn is part of the zero level set $\mathrm{Im}(h(\xi))=0$.
Moreover, it is clear from the definition of $h$ that the sum of the boundary values taken by $h'$ on $\Sigma$ is $h'_+(\xi)+h'_-(\xi)=0$. Therefore, $h_+(\xi)+h_-(\xi)$, $\xi\in\Sigma$, is a~real constant. Evaluating this constant for $\xi=0\in\Sigma$ using the fact that $h$ is an odd function of~$\xi$ shows that $h_+(\xi)+h_-(\xi)=0$ holds identically for $\xi\in\Sigma$. Since $h$ is continuous at $\xi=\pm\ii\nu$, this implies that $h(\pm\ii\nu)=0$. However, the limiting values taken at $\xi=0$ are opposite and nonzero:
\begin{gather}
h(0\pm):=\mathop{\lim_{\xi\to 0}}_{\pm\xi>0}h(\xi)=g(0\pm)=\pm\int_0^{+\infty}\left[\frac{3}{2}\xi^2+\frac{1}{2}-\frac{3}{2}\xi\sqrt{\xi^2+\nu^2}\right]\,\dd\xi = \pm\frac{1}{3}\nu.\label{eq:h-zero}
\end{gather}

\begin{Remark}The formula \eqref{eq:gprime} can be motivated by the desire to enforce the condition that $h_+'(\xi)+h_-'(\xi)=0$ or equivalently that $g_+'(\xi)+g_-'(\xi)+\varphi'(\xi;1)=0$ holds on $\Sigma$. Indeed, the latter equation can be re-arranged to read $\big(g'_+(\xi)+\tfrac{3}{2}\xi^2+\tfrac{1}{2}\big)=-\big(g'_-(\xi)+\tfrac{3}{2}\xi^2+\tfrac{1}{2}\big)$, so $\big(g'(\xi)+\tfrac{3}{2}\xi^2+\tfrac{1}{2}\big)^2$ should have no jump across $\Sigma$ and hence can be sought as an entire function. Imposing the integrability condition that $g'(\xi)=O\big(\xi^{-2}\big)$ as $\xi\to\infty$ determines this entire function up to a constant as the polynomial $\tfrac{9}{4}\xi^4+\tfrac{3}{2}\xi^2+c$. The formula \eqref{eq:gprime} then corresponds to the choice $c=0$. This is the only choice of $c$ for which the polynomial has two simple roots and one double root.
\end{Remark}

The central panel of Fig.~\ref{fig:PII-Plus} shows the original jump contour for $\mathbf{Y}(\xi;y,1)$ together with six regions denoted $A$, $B$, $C$, $D$, $E$, and $F$. We use the function $g(\xi)$ to give a piecewise-analytic definition of a new matrix unknown as follows:
\begin{gather*}
\mathbf{Z}(\xi;y):=\mathbf{Y}(\xi;y,1)\begin{bmatrix}1 & -\tau\ee^{-2\pi p}\ee^{\ii y^{3/2}\varphi(\xi;1)}\\0 & 1\end{bmatrix}\ee^{-\ii y^{3/2}g(\xi)\sigma_3},\qquad\xi\in A,\\
\mathbf{Z}(\xi;y):=\mathbf{Y}(\xi;y,1)\begin{bmatrix}1 & -\tau^{-1}\ee^{\ii y^{3/2}\varphi(\xi;1)}\\0 & 1\end{bmatrix}\ee^{-\ii y^{3/2}g(\xi)\sigma_3},\qquad \xi\in B,\\
\mathbf{Z}(\xi;y):=\mathbf{Y}(\xi;y,1)\begin{bmatrix}1 & \tau^{-1}\ee^{\ii y^{3/2}\varphi(\xi;1)}\\0 & 1\end{bmatrix}\ee^{-\ii y^{3/2}g(\xi)\sigma_3},\qquad\xi\in C, \\
\mathbf{Z}(\xi;y):=\mathbf{Y}(\xi;y,1)\begin{bmatrix}1 & 0\\\tau\ee^{-2\pi p}\ee^{-\ii y^{3/2}\varphi(\xi;1)} & 1\end{bmatrix}\ee^{-\ii y^{3/2}g(\xi)\sigma_3},\qquad\xi\in D,\\
\mathbf{Z}(\xi;y):=\mathbf{Y}(\xi;y,1)\begin{bmatrix}1 & 0\\\tau^{-1}\ee^{-\ii y^{3/2}\varphi(\xi;1)} & 1\end{bmatrix}\ee^{-\ii y^{3/2}g(\xi)\sigma_3},\qquad\xi\in E,\\
\mathbf{Z}(\xi;y):=\mathbf{Y}(\xi;y,1)\begin{bmatrix}1 & 0\\-\tau^{-1}\ee^{-\ii y^{3/2}\varphi(\xi;1)}&1\end{bmatrix}\ee^{-\ii y^{3/2}g(\xi)\sigma_3},\qquad\xi\in F,
\end{gather*}
and for $\xi\in\mathbb{C}\setminus(\overline{A}\cup\overline{B}\cup\overline{C}\cup\overline{D}\cup\overline{E}\cup\overline{F})$ we set $\mathbf{Z}(\xi;y):=\mathbf{Y}(\xi;y,1)\ee^{-\ii y^{3/2}g(\xi)\sigma_3}$. The resulting jump contour $\Sigma_\mathbf{Z}$ for $\mathbf{Z}(\xi;y)$ is illustrated in the right-hand panel of Fig.~\ref{fig:PII-Plus}. $\mathbf{Z}(\xi;y)$ is analytic for $\xi\in\mathbb{C}\setminus\Sigma_\mathbf{Z}$, and across the arcs of $\Sigma_\mathbf{Z}$ the jump condition $\mathbf{Z}_+(\xi;y)=\mathbf{Z}_-(\xi;y)\mathbf{V}^\mathbf{Z}(\xi;y)$ holds, where the jump matrix $\mathbf{V}^\mathbf{Z}(\xi;y)$ is defined as follows:
\begin{gather}\allowdisplaybreaks
\mathbf{V}^\mathbf{Z}(\xi;y):=\begin{bmatrix}1 & \tau\ee^{-2\pi p}\ee^{2\ii y^{3/2}h(\xi)}\\0 & 1\end{bmatrix},\qquad\xi\in\alpha^+,\label{eq:alpha-plus}
\\
\mathbf{V}^\mathbf{Z}(\xi;y):=\begin{bmatrix}1 & \tau^{-1}\ee^{-2\pi p}\ee^{2\ii y^{3/2}h(\xi)}\\0 & 1\end{bmatrix},\qquad\xi\in\beta^+,\label{eq:beta-plus}
\\
\mathbf{V}^\mathbf{Z}(\xi;y):=\begin{bmatrix}1 & \tau^{-1}\ee^{2\ii y^{3/2}h(\xi)}\\0 & 1\end{bmatrix},\qquad\xi\in\gamma^+,\nonumber
\\
\mathbf{V}^\mathbf{Z}(\xi;y):=\begin{bmatrix}0 & -\tau^{-1}\\ \tau & 0\end{bmatrix},\qquad\xi\in\Sigma^+:=\{\xi\in\Sigma,\,\mathrm{Im}(\xi)>0\},\label{eq:PII-Plus-jump-Sigma-up}
\\
\mathbf{V}^\mathbf{Z}(\xi;y):=\begin{bmatrix}1 & 0\\\tau\ee^{-2\ii y^{3/2}h(\xi)}&1\end{bmatrix},\qquad \xi\in\ii\mathbb{R}_+\setminus\Sigma^+,\nonumber
\\
\mathbf{V}^\mathbf{Z}(\xi;y):=\ee^{2\pi p\sigma_3},\qquad\xi<0,\label{eq:PII-Plus-jump-diagonal}
\\
\mathbf{V}^\mathbf{Z}(\xi;y):=\begin{bmatrix}1&0\\\tau\ee^{-2\pi p}\ee^{-2\ii y^{3/2}h(\xi)} & 1\end{bmatrix},\qquad\xi\in \alpha^-,\nonumber
\\
\mathbf{V}^\mathbf{Z}(\xi;y):=\begin{bmatrix}1&0\\\tau^{-1}\ee^{-2\pi p}\ee^{-2\ii y^{3/2}h(\xi)} & 1 \end{bmatrix},\qquad\xi\in\beta^-, \nonumber
\\
\mathbf{V}^\mathbf{Z}(\xi;y):=\begin{bmatrix}1&0\\\tau^{-1}\ee^{-2\ii y^{3/2}h(\xi)}&1\end{bmatrix},\qquad\xi\in\gamma^-,\nonumber
\\
\mathbf{V}^\mathbf{Z}(\xi;y):=\begin{bmatrix}0 & \tau\\-\tau^{-1} & 0\end{bmatrix},\qquad \xi\in\Sigma^-:=\{\xi\in\Sigma,\,\mathrm{Im}(\xi)<0\},\label{eq:PII-Plus-jump-Sigma-down}
\\
\mathbf{V}^\mathbf{Z}(\xi;y):=\begin{bmatrix}1 & \tau\ee^{2\ii y^{3/2}h(\xi)}\\0 & 1\end{bmatrix},\qquad \xi\in\ii\mathbb{R}_-\setminus\Sigma^-.\label{eq:PII-Plus-jump-last}
\end{gather}
Due to the placement of the jump contour $\Sigma_\mathbf{Z}$ relative to the sign chart for $\mathrm{Im}(h(\xi))$ illustrated in Fig.~\ref{fig:PII-Plus}, the jump matrix $\mathbf{V}^\mathbf{Z}(\xi;y)$ converges exponentially fast to the identity pointwise on~$\Sigma_\mathbf{Z}$ as $y\to +\infty$, with the exception of the two halves of~$\Sigma$ and the negative real line, see \eqref{eq:PII-Plus-jump-Sigma-up}, \eqref{eq:PII-Plus-jump-diagonal}, and~\eqref{eq:PII-Plus-jump-Sigma-down}.

\subsection{Parametrix construction}
\subsubsection{Outer parametrix}
We construct an outer parametrix $\dot{\mathbf{Z}}^\mathrm{out}(\xi)$ to satisfy the latter jump conditions as follows.
Since $\xi^{-\ii p\sigma_3}$ satisfies exactly the normalization condition $\mathbf{Z}(\xi;y)\xi^{\ii p\sigma_3}\to\mathbb{I}$ as $\xi\to\infty$, as well as the jump condition for $\mathbf{Z}$ across the negative real axis, we seek the outer parametrix in the form $\dot{\mathbf{Z}}^\mathrm{out}(\xi):=
\mathbf{J}(\xi)\xi^{-\ii p\sigma_3}$ and we require that $\mathbf{J}(\xi)\to\mathbb{I}$ as $\xi\to\infty$ and that $\mathbf{J}(\xi)$ be analytic except on $\Sigma$, where
\begin{gather}
\mathbf{J}_+(\xi)=\mathbf{J}_-(\xi)\begin{bmatrix}0 & -\tau^{-1}\xi^{-2\ii p}\\\tau\xi^{2\ii p} & 0\end{bmatrix},\qquad \xi\in\Sigma^+, \label{eq:PII-G-jump}
\end{gather}
and
\begin{gather}
\mathbf{J}_+(\xi)=\mathbf{J}_-(\xi)\begin{bmatrix}0 & \tau\xi^{-2\ii p}\\-\tau^{-1}\xi^{2\ii p} & 0\end{bmatrix},\qquad\xi\in\Sigma^-. \label{eq:PII-G-jump-LHP}
\end{gather}
These choices guarantee that $\dot{\mathbf{Z}}^\mathrm{out}(\xi)$ exactly satisfies the jump conditions of $\mathbf{Z}$ described by the jump matrix in \eqref{eq:PII-Plus-jump-Sigma-up}, \eqref{eq:PII-Plus-jump-diagonal}, and \eqref{eq:PII-Plus-jump-Sigma-down}. Given $p\in\mathbb{C}\setminus\ii\mathbb{Z}$, observe that there exists a unique number $\ell\in\mathbb{C}$ such that $\tau=\ee^{\ell}$ solves the equation $\tau^2=\ee^{2\pi p}-1$ and such that the following inequalities hold:
\begin{gather}
-1<\mathrm{Re}\left(2\ii p + \frac{2\ell}{\ii \pi}\right)\le 1.\label{eq:ell-inequalities}
\end{gather}
Indeed, the relation $\tau^2=\ee^{2\pi p}-1$ determines $\ell$ modulo $\ii\pi\mathbb{Z}$, and hence a unique definite value of $\ell$ can be chosen so that $\ee^{2\ell}=\tau^2=\ee^{2\pi p}-1$ and the inequalities \eqref{eq:ell-inequalities} hold. We are thus breaking the symmetry $\tau\mapsto -\tau$ mentioned in the remark at the end of Section~\ref{sec:DE} in order to facilitate subsequent asymptotic analysis. We note that in the application of Riemann--Hilbert Problem~\ref{rhp:PII-generalized} to rogue waves in~\cite{BilmanLM18}, it is necessary to have $\tau=1$ for $p=\ln(2)/(2\pi)$; this amounts to choosing $\ell=0$ which is consistent with the inequalities \eqref{eq:ell-inequalities}, so we can be assured that the analysis that follows applies to that special case.

Now, given $p$ and $\ell$ as above, let $j(\xi)$ be defined by
\begin{gather*}
j(\xi):=\frac{pR(\xi)}{\pi}\int_{\Sigma}\frac{\log(s)\,\dd s}{R_+(s)(s-\xi)} + \ii\eta \\
\hphantom{j(\xi):=}{} + \frac{\ell R(\xi)}{2\pi\ii}\left[\int_{\Sigma^+}\frac{\dd s}{R_+(s)(s-\xi)}-\int_{\Sigma^-}\frac{\dd s}{R_+(s)(s-\xi)}\right],\qquad\xi\in\mathbb{C}\setminus\Sigma,
\end{gather*}
where the constant $\eta$ is defined by
\begin{gather*}
\eta:=\frac{p}{\ii\pi}\int_\Sigma\frac{\log(s)\,\dd s}{R_+(s)}.
\end{gather*}
Here, $\log(s)$ refers to the principal branch, and $R_+(s)$ refers to the boundary value from the left by orientation of the integration contour (the direction of which is irrelevant for making $j$ and $\eta$ well-defined given the factor $R_+(s)$ in the integrand). Thus, taking the integration in the upward direction, and parametrizing by $s=\ii t$, we have $R_+(\ii t)=-\sqrt{\nu^2-t^2}$, and hence
\begin{gather}
\eta = -\frac{2p}{\pi}\int_0^{\sqrt{\tfrac{2}{3}}}\frac{\ln(t)\,\dd t}{\sqrt{\tfrac{2}{3}-t^2}} =\frac{p}{2}\ln(6).\label{eq:eta-value}
\end{gather}
The function $j$ is analytic for $\xi\in\mathbb{C}\setminus\Sigma$ and the value of the constant $\eta$ is determined so that $j(\xi)\to 0$ as $\xi\to\infty$. In fact, it is not difficult to establish that
\begin{gather}
j(\xi)=-\ii\nu \left(p-\frac{\ell}{\pi}\right)\xi^{-1}-\frac{1}{6}\ii p\xi^{-2} + O\big(\xi^{-3}\big),\qquad\xi\to\infty.
\label{eq:j-large-xi}
\end{gather}
The boundary values $j_\pm(\xi)$ taken on $\Sigma^+$ and $\Sigma^-$ are continuous except at the origin, but including the other endpoints $\xi=\pm\ii\nu$. Since $R(\xi)$ changes sign across $\Sigma=\Sigma^+\cup\Sigma^-$, the sum of the boundary values is easy to compute directly:
\begin{gather*}
j_+(\xi)+j_-(\xi)=\begin{cases}
2\ii p\log(\xi)+\ell+2\ii\eta,& \xi\in\Sigma^+,\quad\mathrm{Im}(\log(\xi))=\tfrac{1}{2}\pi,\\
2\ii p\log(\xi) -\ell + 2\ii\eta,& \xi\in\Sigma^-,\quad\mathrm{Im}(\log(\xi))=-\tfrac{1}{2}\pi.
\end{cases}
\end{gather*}
Let $D(0)$ denote a disk of sufficiently small radius (less than $\tfrac{1}{2}\nu=1/\sqrt{6}$ will do) centered at the origin. It is straightforward to show that
\begin{gather}
j(\xi)=\left(2\ii p+\frac{\ell}{\ii\pi}\right)\log(\xi) + k(\xi),\qquad \xi\in D(0),\qquad\mathrm{Re}(\xi)>0,
\label{eq:j-log}
\end{gather}
where $\log(\xi)$ again denotes the principal branch, and where $k(\xi)$ is a function analytic in the full disk $D(0)$. Letting $R_1(\xi)$ denote the analytic continuation of $R(\xi)$ from the right half of $D(0)$ to all of $D(0)$, a formula for $k(\xi)$ that admits direct evaluation for any $\xi\in D(0)$ reads
\begin{gather*}
k(\xi)=\ii\eta + \frac{p}{\pi}R_1(\xi)\left[\int_{-M}^{\ii\nu}\frac{\log(s)\,\dd s}{R(s)(s-\xi)}+\int_{-\ii\nu}^{-M}\frac{\log(s)\,\dd s}{R(s)(s-\xi)}\right]\nonumber\\
\hphantom{k(\xi)=}{}+\frac{\ell}{2\pi\ii}R_1(\xi)\left[\int_{-M}^{\ii\nu}\frac{\dd s}{R(s)(s-\xi)} - \int_{-\ii\nu}^{-M}\frac{\dd s}{R(s)(s-\xi)}\right]\nonumber\\
\hphantom{k(\xi)=}{}-\left(2\ii p+\frac{\ell}{\ii\pi}\right)\log(\xi+M) + \left(2\ii p+\frac{\ell}{\ii\pi}\right)R_1(\xi)\int_{-M}^0\left[\frac{1}{R_1(s)}-\frac{1}{R_1(\xi)}\right]\frac{\dd s}{s-\xi},
%\label{eq:k-formula}
\end{gather*}
where $M>0$ is any sufficiently large constant, and where the contour integrals on the first two lines are taken over straight line segments. It is an exercise to show from this formula that
\begin{gather}
k(0)=\ii\eta -\ii\left(p-\frac{\ell}{\pi}\right)\ln(2\nu).
\label{eq:k-zero}
\end{gather}
Now set $\mathbf{K}(\xi):=\mathbf{J}(\xi)\ee^{-j(\xi)\sigma_3}$. Then $\mathbf{J}(\xi)\to \mathbb{I}$ as $\xi\to\infty$ implies that also $\mathbf{K}(\xi)\to\mathbb{I}$ as $\xi\to\infty$, and $\mathbf{J}$ analytic in $\mathbb{C}\setminus\Sigma$ implies that the same is true for $\mathbf{K}$. According to \eqref{eq:PII-G-jump}, $\mathbf{K}(\xi)$ satisfies the jump condition
\begin{gather*}
\mathbf{K}_+(\xi)=\mathbf{K}_-(\xi)\begin{bmatrix}0 & -\ee^{2\ii\eta}\\\ee^{-2\ii\eta} & 0\end{bmatrix},\qquad \xi\in\Sigma\;\;\text{(here taken with upward orientation).}
%\label{eq:PII-H-jump}
\end{gather*}
We take for $\mathbf{K}(\xi)$ the following solution:
\begin{gather}
\mathbf{K}(\xi):=\ee^{\ii\eta\sigma_3}\mathbf{S}\left(\frac{\xi-\ii\nu}{\xi+\ii\nu}\right)^{\sigma_3/4}\mathbf{S}^{-1}\ee^{-\ii\eta\sigma_3},\qquad \mathbf{S}:=\frac{1}{\sqrt{2}}\begin{bmatrix}1 & -\ii\\-\ii & 1\end{bmatrix}.
\label{eq:K-plus-def}
\end{gather}
This function takes continuous boundary values along $\Sigma$ except at the endpoints, where it exhibits negative fourth-root singularities. Moreover, its boundary values from the left and right half-planes admit analytic continuation into the full neighborhood $D(0)$. The outer parametrix is thus defined by the formula
\begin{gather}
\dot{\mathbf{Z}}^\mathrm{out}(\xi):=\mathbf{K}(\xi)\ee^{j(\xi)\sigma_3}\xi^{-\ii p\sigma_3},\qquad\xi\in\mathbb{C}\setminus(\Sigma\cup \mathbb{R}_-).\label{eq:outer-plus}
\end{gather}

\subsubsection[Inner parametrices near $\xi=\pm\ii\nu$]{Inner parametrices near $\boldsymbol{\xi=\pm\ii\nu}$}
Let $D(\ii\nu)$ denote a disk of sufficiently small radius (less than $\tfrac{1}{2}\nu$) centered at $\xi=\ii\nu$. We define a conformal mapping $\xi\mapsto W_{\ii\nu}(\xi)$ of $D(\ii\nu)$ onto a neighborhood of the origin by the equation
\begin{gather*}
2\ii h(\xi)=W_{\ii\nu}(\xi)^{3/2}
\end{gather*}
in which we understand that both sides of this equation are positive real for $\xi$ on the imaginary axis above $\xi=\ii\nu$. To define the conformal map, we choose the solution $W_{\ii\nu}(\xi)$ that is also positive for such $\xi$ before analytically continuing the resulting solution to $D(\ii\nu)$. That this procedure succeeds is a consequence of the formula $h'(\xi)=\tfrac{3}{2}\xi R(\xi)$ and the fact that $h(\ii\nu)=0$ since $h_+(\xi)+h_-(\xi)=0$ holds identically on $\Sigma$ and $h$ is continuous at $\xi=\ii\nu$.

Suppose that within the disk $D(\ii\nu)$, the contour arcs $\alpha^+$ and $\beta^+$ are merged into a single arc carrying the product of their (commuting) jump matrices (see \eqref{eq:alpha-plus}--\eqref{eq:beta-plus}), and that this arc lies along the ray $\arg(W_{\ii\nu}(\xi))=\tfrac{2}{3}\pi$. Likewise, we take $\gamma^+$ to lie within $D(\ii\nu)$ along the ray $\arg(W_{\ii\nu}(\xi))=-\tfrac{2}{3}\pi$, while $W_{\ii\nu}$ automatically maps the imaginary axis near $\xi=\ii\nu$ to the real axis so that $\Sigma^+$ lies along $\arg(-W_{\ii\nu}(\xi))=0$ and the imaginary axis above $\xi=\ii\nu$ is taken to the positive real axis. Then, the identity $\tau^2=\ee^{2\pi p}-1$ implies that the jump conditions satisfied by the product $\mathbf{Z}(\xi;y)\tau^{-\sigma_3/2}(\ii\sigma_1)$ (here $\tau^{1/2}$ denotes any fixed square root of $\tau\neq 0$) read as follows for $\xi\in D(\ii\nu)$, where $\lambda:=yW_{\ii\nu}(\xi)$:
\begin{gather}
\mathbf{Z}_+(\xi;y)\tau^{-\sigma_3/2}(\ii\sigma_1)=\mathbf{Z}_-(\xi;y)\tau^{-\sigma_3/2}(\ii\sigma_1)\begin{bmatrix}1 & \ee^{-\lambda^{3/2}}\\0 & 1\end{bmatrix},\qquad \arg(\lambda)=0,\label{eq:Airy-1}
\\
\mathbf{Z}_+(\xi;y)\tau^{-\sigma_3/2}(\ii\sigma_1)=\mathbf{Z}_-(\xi;y)\tau^{-\sigma_3/2}(\ii\sigma_1)\begin{bmatrix}1 & 0\\\ee^{\lambda^{3/2}} & 1\end{bmatrix},\qquad \arg(\lambda)=\pm\frac{2}{3}\pi,\nonumber
\end{gather}
and
\begin{gather}
\mathbf{Z}_+(\xi;y)\tau^{-\sigma_3/2}(\ii\sigma_1)=\mathbf{Z}_-(\xi;y)\tau^{-\sigma_3/2}(\ii\sigma_1)\begin{bmatrix}0 & 1\\-1 & 0\end{bmatrix},\qquad \arg(-\lambda)=0.\label{eq:Airy-3}
\end{gather}
To properly interpret these jump conditions, note that all contours carry the orientation induced from that shown in the right-hand panel of Fig.~\ref{fig:PII-Plus} by the conformal mapping $\xi\mapsto W_{\ii\nu}(\xi)$, which means that the rays in the $\lambda$-plane are oriented in the direction of increasing real part.

From the definition \eqref{eq:outer-plus} of the outer parametrix $\dot{\mathbf{Z}}^\mathrm{out}(\xi)$ that there exists a matrix-valued function $\mathbf{H}_{\ii\nu}(\xi)$ analytic in $D(\ii\nu)$ and having unit determinant, such that
\begin{gather*}
\dot{\mathbf{Z}}^\mathrm{out}(\xi)\tau^{-\sigma_3/2}(\ii\sigma_1)=\mathbf{H}_{\ii\nu}(\xi)W_{\ii\nu}(\xi)^{\sigma_3/4}\mathbf{S}=\mathbf{H}_{\ii\nu}(\xi)y^{-\sigma_3/4}\lambda^{\sigma_3/4}\mathbf{S},\qquad\xi\in D(\ii\nu),
\end{gather*}
where $\mathbf{S}$ denotes the same eigenvector matrix defined in \eqref{eq:K-plus-def}. Because we would like to build an inner parametrix that matches well onto the outer parametrix when $\xi\in\partial D(\ii\nu)$, we take guidance from the final two factors $\lambda^{\sigma_3/4}\mathbf{S}$ and define a matrix $\mathbf{A}(\lambda)$ as the solution of the following Riemann--Hilbert Problem.
\begin{rhp}[Airy parametrix]\label{rhp:Airy}Seek a $2\times 2$ matrix function $\mathbf{A}(\lambda)$ defined for $0<|\arg(\lambda)|<\tfrac{2}{3}\pi$ and $\tfrac{2}{3}\pi<|\arg(\lambda)|<\pi$ with the following properties:
\begin{itemize}\itemsep=0pt
\item[]\textbf{Analyticity:} $\mathbf{A}(\lambda)$ is analytic in the four sectors indicated above, and takes continuous boundary values from each sector including at the origin.
\item[]\textbf{Jump conditions:} The boundary values $\mathbf{A}_\pm(\lambda)$ are related by exactly the same jump conditions satisfied by $\mathbf{Z}(\xi;y)\tau^{-\sigma_3/2}(\ii\sigma_1)$ $($see \eqref{eq:Airy-1}--\eqref{eq:Airy-3}$)$ but with the indicated directions extended to infinite rays in the $\lambda$-plane oriented in the direction of increasing real part.
\item[]\textbf{Normalization:} $\mathbf{A}(\lambda)\mathbf{S}^{-1}\lambda^{-\sigma_3/4}\to\mathbb{I}$ as $\lambda\to\infty$.
\end{itemize}
\end{rhp}
This problem is well-known to have a unique solution constructed from Airy functions. The solution has the additional property that the normalization condition is strengthened to
\begin{gather}
\mathbf{A}(\lambda)\mathbf{S}^{-1}\lambda^{-\sigma_3/4}=\mathbb{I} +\begin{bmatrix}O\big(\lambda^{-3}\big) & O(\lambda^{-1})\\O\big(\lambda^{-2}\big) & O\big(\lambda^{-3}\big)\end{bmatrix},\qquad\lambda\to\infty.
\label{eq:Airy-better-norm}
\end{gather}
We then define an inner parametrix for $\mathbf{Z}(\xi;y)$ as follows:
\begin{gather*}
\dot{\mathbf{Z}}^{\mathrm{in}}_{\ii\nu}(\xi;y):=\mathbf{H}_{\ii\nu}(\xi)y^{-\sigma_3/4}\mathbf{A}(yW_{\ii\nu}(\xi))(-\ii\sigma_1)\tau^{\sigma_3/2},\qquad\xi\in D(\ii\nu).
\end{gather*}
Since $W_{\ii\nu}(\xi)$ is bounded away from zero and hence $\lambda$ is proportional to $y$ when $\xi\in\partial D(\ii\nu)$, we then get from \eqref{eq:Airy-better-norm} that
\begin{gather}
\dot{\mathbf{Z}}^\mathrm{in}_{\ii\nu}(\xi;y)\dot{\mathbf{Z}}^\mathrm{out}(\xi)^{-1} =
\mathbf{H}_{\ii\nu}(\xi)y^{-\sigma_3/4}\mathbf{A}(yW_{\ii\nu}(\xi))\mathbf{S}^{-1}(yW_{\ii\nu}(\xi))^{-\sigma_3/4}y^{\sigma_3/4}\mathbf{H}_{\ii\nu}(\xi)^{-1}\nonumber\\
\hphantom{\dot{\mathbf{Z}}^\mathrm{in}_{\ii\nu}(\xi;y)\dot{\mathbf{Z}}^\mathrm{out}(\xi)^{-1}}{}
=\mathbb{I} + \mathbf{H}_{\ii\nu}(\xi)y^{-\sigma_3/4}\begin{bmatrix}O\big(y^{-3}\big) & O\big(y^{-1}\big)\\ O\big(y^{-2}\big) & O\big(y^{-3}\big)\end{bmatrix}y^{\sigma_3/4}\mathbf{H}_{\ii\nu}(\xi)^{-1}\nonumber\\
\hphantom{\dot{\mathbf{Z}}^\mathrm{in}_{\ii\nu}(\xi;y)\dot{\mathbf{Z}}^\mathrm{out}(\xi)^{-1}}{}
=\mathbb{I}+\mathbf{H}_{\ii\nu}(\xi)\begin{bmatrix}O\big(y^{-3}\big) & O\big(y^{-3/2}\big)\\O\big(y^{-3/2}\big) & O\big(y^{-3}\big)\end{bmatrix}\mathbf{H}_{\ii\nu}(\xi)^{-1}\nonumber\\
\hphantom{\dot{\mathbf{Z}}^\mathrm{in}_{\ii\nu}(\xi;y)\dot{\mathbf{Z}}^\mathrm{out}(\xi)^{-1}}{}
=\mathbb{I}+O\big(y^{-3/2}\big),\qquad\xi\in\partial D(\ii\nu),\label{eq:Airy-mismatch}
\end{gather}
where the estimate on the last line holds in the uniform sense on the circle $\partial D(\ii\nu)$ in the limit $y\to +\infty$.

A very similar construction gives an explicit parametrix $\dot{\mathbf{Z}}^\mathrm{in}_{-\ii\nu}(\xi;y)$ defined in a small disk $D(-\ii\nu)$ centered at $\xi=-\ii\nu$ in terms of Airy functions, satisfying exactly the jump conditions of~$\mathbf{Z}$ within the disk and for which the estimate $\dot{\mathbf{Z}}^\mathrm{in}_{-\ii\nu}(\xi;y)\dot{\mathbf{Z}}^\mathrm{out}(\xi)^{-1}=\mathbb{I}+O\big(y^{-3/2}\big)$ holds uniformly on $\partial D(-\ii\nu)$. We will require no further details of these inner parametrices.

\subsubsection[Inner parametrix near $\xi=0$]{Inner parametrix near $\boldsymbol{\xi=0}$}
Recall the disk $D(0)$ of sufficiently small radius centered at the origin. We define a conformal mapping in terms of $h(\xi)$ within $D(0)$ as follows. Here a slightly different approach is required because $h$ is locally two opposite analytic functions in the left and right half-disks, each of which has a simple critical point at the origin. Thus, we will define a mapping $W_0(\xi)$ on $D(0)$ by two different conditions:
\begin{gather}
h(\xi)-h(0\pm)=\pm\frac{1}{2}W_0(\xi)^2,\qquad \xi\in D(0),\qquad \pm\mathrm{Re}(\xi)>0.\label{eq:W0-conf-map}
\end{gather}
Here, $h(0\pm)$ denotes the limiting value of $h(\xi)$ at $\xi=0$ from the domain $\pm\mathrm{Re}(\xi)>0$; see \eqref{eq:h-zero}. These conditions determine $W_0(\xi)$ up to an overall sign as a conformal mapping of the disk $D(0)$ with the property that $W_0(0)=0$. We fix the sign by insisting that $W_0'(0)>0$. In fact, by taking two derivatives in \eqref{eq:W0-conf-map} and setting $\xi=0$ we see that
\begin{gather}
W_0'(0)=\frac{1}{\sqrt{\nu}} = \left(\frac{3}{2}\right)^{1/4}>0.
\label{eq:W0-prime}
\end{gather}
We deform the jump contour within $D(0)$ to ensure both that $\arg(W_0(\gamma^\pm))=\pm\tfrac{1}{4}\pi$ and that $\arg(W_0(\beta^\pm))=\pm\tfrac{3}{4}\pi$. The conformal mapping is real, and therefore it automatically holds that $\arg(-\xi)=0$ corresponds to $\arg(-W_0(\xi))=0$.

Consider the matrix $\mathbf{X}(\xi;y)$ explicitly related to $\mathbf{Z}(\xi;y)$ for $\xi\in D(0)$ as follows:
\begin{gather}
\mathbf{X}(\xi;y):=\begin{cases}\mathbf{Z}(\xi;y)\begin{bmatrix}0 & \tau^{-1}\\-\tau & 0\end{bmatrix}
\ee^{\ii y^{3/2}h(0+)\sigma_3}(\ii\sigma_1),& \tfrac{1}{2}\pi<\arg(\xi)<\pi,\vspace{1mm}\\
\mathbf{Z}(\xi;y)\begin{bmatrix}0 & \tau\\-\tau^{-1} & 0\end{bmatrix}\ee^{\ii y^{3/2}h(0+)\sigma_3}(\ii\sigma_1),& -\pi<\arg(\xi)<-\tfrac{1}{2}\pi,\vspace{1mm}\\
\mathbf{Z}(\xi;y)\ee^{\ii y^{3/2}h(0+)\sigma_3}(\ii\sigma_1),& |\arg(\xi)|<\tfrac{1}{2}\pi.
\end{cases}\label{eq:X-from-Z}
\end{gather}
From this definition, it follows that
\begin{gather}
\mathbf{X}_+(\xi;y)=\mathbf{X}_-(\xi;y),\qquad\xi\in\Sigma^\pm,\label{eq:X-jump-0}
\end{gather}
so that $\mathbf{X}(\xi;y)$ may be considered to be analytic on $\Sigma^\pm$. The remaining jump conditions sa\-tisfied by $\mathbf{X}(\xi;y)$ are easily expressed in terms of the rescaled conformal mapping with a new independent variable $\lambda:=y^{3/4}W_0(\xi)$. They read as follows:
\begin{gather}
\mathbf{X}_+(\xi;y)=\mathbf{X}_-(\xi;y)\begin{bmatrix}1 & 0\\\tau^{-1}\ee^{\ii\lambda^2} & 1\end{bmatrix},\qquad\arg(W_0(\xi))=\frac{1}{4}\pi,\label{eq:X-jump-1}
\\
\mathbf{X}_+(\xi;y)=\mathbf{X}_-(\xi;y)\begin{bmatrix}1 & \tau^{-1}\ee^{-\ii\lambda^2}\\0 & 1\end{bmatrix},\qquad\arg(W_0(\xi))=-\frac{1}{4}\pi,\nonumber
\\
\mathbf{X}_+(\xi;y)=\mathbf{X}_-(\xi;y)\begin{bmatrix}1 & \tau\ee^{-2\pi p}\ee^{-\ii\lambda^2}\\0 & 1\end{bmatrix},\qquad\arg(W_0(\xi))=\frac{3}{4}\pi,\nonumber
\\
\mathbf{X}_+(\xi;y)=\mathbf{X}_-(\xi;y)\begin{bmatrix}1 & 0\\\tau\ee^{-2\pi p}\ee^{\ii\lambda^2} & 1\end{bmatrix},\qquad\arg(W_0(\xi))=-\frac{3}{4}\pi,\nonumber
\end{gather}
and finally,
\begin{gather}
\mathbf{X}_+(\xi;y)=\mathbf{X}_-(\xi;y)\ee^{2\pi p\sigma_3}\tau^{-2\sigma_3},\qquad\arg(-W_0(\xi))=0, \label{eq:X-jump-5}
\end{gather}
where all contours are oriented in the direction of increasing real part (so $\beta^\pm$ have been re-orien\-ted). To confirm these jump conditions, it is necessary to use the fact that \smash{$h(0-)+h(0+)=0$}.

The jump conditions \eqref{eq:X-jump-0}--\eqref{eq:X-jump-5} match those of the parabolic cylinder function parametrix described in the appendix, provided we replace the parameters $(p,\tau)$ with $(q,\mu)$, where (recall that $\ell$ is a specific number for which $\tau=\ee^{\ell}$)
\begin{gather}
q:=p-\frac{\ell}{\pi}\qquad\text{and}\qquad\mu:=\frac{1}{\tau}. \label{eq:p1-tau1}
\end{gather}
It is easy to check that the condition $\tau^2=\ee^{2\pi p}-1$ guarantees that also $\mu^2=\ee^{2\pi q}-1$. It then follows from the definition \eqref{eq:outer-plus} of the outer parametrix $\dot{\mathbf{Z}}^\mathrm{out}(\xi)$ and from~\eqref{eq:j-log} that there exists a matrix function $\mathbf{H}_0(\xi)$, independent of $y$ and having unit determinant, that is analytic in~$D(0)$ such that
\begin{gather*}
\dot{\mathbf{Z}}^\mathrm{out}(\xi)\ee^{\ii y^{3/2}h(0+)\sigma_3}(\ii\sigma_1)= \mathbf{H}_0(\xi)\ee^{-\ii y^{3/2}h(0+)\sigma_3}W_0(\xi)^{-\ii q\sigma_3}\\
\hphantom{\dot{\mathbf{Z}}^\mathrm{out}(\xi)\ee^{\ii y^{3/2}h(0+)\sigma_3}(\ii\sigma_1)}{} =\mathbf{H}_0(\xi)\ee^{-\ii y^{3/2}h(0+)\sigma_3}y^{3\ii q\sigma_3/4}\lambda^{-\ii q\sigma_3},\qquad\xi\in D(0),\qquad\mathrm{Re}(\xi)>0.
\end{gather*}
Note that, in particular,
\begin{gather}
\mathbf{H}_0(0)=\ee^{\ii\eta\sigma_3}\mathbf{S}\ee^{-\ii\pi\sigma_3/4}\mathbf{S}^{-1}\ee^{-\ii\eta\sigma_3}\ee^{k(0)\sigma_3}W_0'(0)^{-\ii q\sigma_3}(\ii\sigma_1).\label{eq:H00}
\end{gather}
Similarly (comparing with \eqref{eq:X-from-Z}),
\begin{gather*}
\dot{\mathbf{Z}}^\mathrm{out}(\xi)\begin{bmatrix}0 & \tau^{-1}\\-\tau & 0\end{bmatrix}\ee^{\ii y^{3/2}h(0+)\sigma_3}(\ii\sigma_1)=\mathbf{H}_0(\xi)\ee^{-\ii y^{3/2}h(0+)\sigma_3}y^{3\ii q\sigma_3/4}\lambda^{-\ii q\sigma_3},\\
\qquad \xi\in D(0),\qquad\frac{1}{2}\pi<\arg(\xi)<\pi,\qquad\text{and}\\
\dot{\mathbf{Z}}^\mathrm{out}(\xi)\begin{bmatrix}0 & \tau\\-\tau^{-1} & 0\end{bmatrix}\ee^{\ii y^{3/2}h(0+)\sigma_3}(\ii\sigma_1)=\mathbf{H}_0(\xi)\ee^{-\ii y^{3/2}h(0+)\sigma_3}y^{3\ii q\sigma_3/4}\lambda^{-\ii q\sigma_3},\\
\qquad \xi\in D(0),\qquad-\pi<\arg(\xi)<-\frac{1}{2}\pi.
\end{gather*}
Taking into account the final factor of $\lambda^{-\ii q\sigma_3}$ to obtain a good match on the boundary $\partial D(0)$, we are led to construct an inner parametrix near the origin by the following formul\ae\
\begin{gather*}
\dot{\mathbf{Z}}^\mathrm{in}_0(\xi;y):=\mathbf{H}_0(\xi)\ee^{-\ii y^{3/2}h(0+)\sigma_3}y^{3\ii q\sigma_3/4}\mathbf{P}(y^{3/4}W_0(\xi);q,\mu)(-\ii\sigma_1)\ee^{-\ii y^{3/2}h(0+)\sigma_3},\\
\hphantom{\dot{\mathbf{Z}}^\mathrm{in}_0(\xi;y):=}{} \xi\in D(0),\qquad\mathrm{Re}(\xi)>0,
\\
\dot{\mathbf{Z}}^\mathrm{in}_0(\xi;y):=
\mathbf{H}_0(\xi)\ee^{-\ii y^{3/2}h(0+)\sigma_3}y^{3\ii q\sigma_3/4}\mathbf{P}(y^{3/4}W_0(\xi);q,\mu)(-\ii\sigma_1)\ee^{-\ii y^{3/2}h(0+)\sigma_3}\begin{bmatrix}0 & -\tau^{-1}\\\tau&0\end{bmatrix},\\
\hphantom{\dot{\mathbf{Z}}^\mathrm{in}_0(\xi;y):=}{} \xi\in D(0),\qquad\frac{1}{2}\pi<\arg(\xi)<\pi,\qquad\text{and}
\\
\dot{\mathbf{Z}}^\mathrm{in}_0(\xi;y):=
\mathbf{H}_0(\xi)\ee^{-\ii y^{3/2}h(0+)\sigma_3}y^{3\ii q\sigma_3/4}\mathbf{P}(y^{3/4}W_0(\xi);q,\mu)(-\ii\sigma_1)\ee^{-\ii y^{3/2}h(0+)\sigma_3}\begin{bmatrix}0 & -\tau\\\tau^{-1}&0\end{bmatrix},\\
\hphantom{\dot{\mathbf{Z}}^\mathrm{in}_0(\xi;y):=}{} \xi\in D(0),\qquad-\pi<\arg(\xi)<-\frac{1}{2}\pi.
\end{gather*}
Here $\mathbf{P}(\lambda;\cdot,\cdot)$ is the solution of Riemann--Hilbert Problem~\ref{rhp:PC-generalized} given in the appendix. We emphasize that $\dot{\mathbf{Z}}_0^\mathrm{in}(\xi;y)$ is an exact local solution of the jump conditions for $\mathbf{Z}(\xi;y)$ within the disk $D(0)$. It satisfies the mismatch condition
\begin{gather}
\dot{\mathbf{Z}}^\mathrm{in}_0(\xi;y)\dot{\mathbf{Z}}^\mathrm{out}(\xi)^{-1}=
\mathbf{H}_0(\xi)\ee^{-\ii y^{3/2}h(0+)\sigma_3}y^{3\ii q\sigma_3/4}\mathbf{P}(\lambda;q,\mu)\lambda^{\ii q\sigma_3}y^{-3\ii q\sigma_3/4}\nonumber\\
\hphantom{\dot{\mathbf{Z}}^\mathrm{in}_0(\xi;y)\dot{\mathbf{Z}}^\mathrm{out}(\xi)^{-1}=}{}\times \ee^{\ii y^{3/2}h(0+)\sigma_3}\mathbf{H}_0(\xi)^{-1},\qquad
\xi\in\partial D(0),\label{eq:PC-Plus-Mismatch}
\end{gather} where $\lambda=y^{3/4}W_0(\xi)$.

\subsection{Error analysis}
The \emph{global parametrix} $\dot{\mathbf{Z}}(\xi;y)$ for $\mathbf{Z}(\xi;y)$ is defined as
\begin{gather*}
\dot{\mathbf{Z}}(\xi;y):=\begin{cases}
\dot{\mathbf{Z}}^\mathrm{in}_{\pm\ii\nu}(\xi;y),&\xi\in D(\pm\ii\nu),\\
\dot{\mathbf{Z}}^\mathrm{in}_0(\xi;y),&\xi\in D(0),\\
\dot{\mathbf{Z}}^\mathrm{out}(\xi),&\xi\in\mathbb{C}\setminus\big(\overline{D(\ii\nu)}\cup\overline{D(0)}\cup\overline{D(-\ii\nu)}\big).
\end{cases}%\label{eq:Plus-global}
\end{gather*}
In a preliminary attempt to gauge the accuracy of the approximation of $\mathbf{Z}(\xi;y)$ by $\dot{\mathbf{Z}}(\xi;y)$, we define the error matrix as $\mathbf{E}(\xi;y):=\mathbf{Z}(\xi;y)\dot{\mathbf{Z}}(\xi;y)^{-1}$ wherever both factors are defined (note that by definition $\det(\dot{\mathbf{Z}}(\xi;y))=1$ holds everywhere). Because the inner parametrices are in each case exact local solutions of the jump conditions for $\mathbf{Z}(\xi;y)$, it is easy to apply Morera's theorem to deduce that $\mathbf{E}(\xi;y)$ can be viewed as a function analytic within all three disks: $D(\pm\ii\nu)$ and~$D(0)$. Similarly, because the outer parametrix exactly satisfies the jump condition for $\mathbf{Z}(\xi;y)$ on the part of the negative real axis outside of $D(0)$ and the parts of $\Sigma^\pm$ lying outside of all three disks, $\mathbf{E}(\xi;y)$ may be considered to be analytic in a neighborhood of these. Therefore, the only points of non-analyticity for $\mathbf{E}(\xi;y)$ correspond to jump discontinuities (i) across the portions of $\alpha^\pm$, $\beta^\pm$, $\gamma^\pm$ and $\ii\mathbb{R}\setminus\Sigma$ lying outside all three disks (where $\mathbf{Z}(\xi;y)$ is discontinuous while $\dot{\mathbf{Z}}(\xi;y)=\dot{\mathbf{Z}}^\mathrm{out}(\xi)$ is analytic) and (ii) across the boundaries of the three disks (where the parametrix $\dot{\mathbf{Z}}(\xi;y)$ is discontinuous while $\mathbf{Z}(\xi;y)$ is analytic).

For jumps of type (i), we directly compute that $\mathbf{E}_+(\xi;y)=\mathbf{E}_-(\xi;y)\dot{\mathbf{Z}}^\mathrm{out}(\xi)\mathbf{V}^\mathbf{Z}(\xi;y)\dot{\mathbf{Z}}^\mathrm{out}(\xi)^{-1}$, where $\mathbf{V}^\mathbf{Z}(\xi;y)$ is defined on the relevant arcs by \eqref{eq:alpha-plus}--\eqref{eq:PII-Plus-jump-last}. The placement of the jump contour~$\Sigma_\mathbf{Z}$ relative to the sign chart for $\mathrm{Im}(h(\xi))$ as shown in Fig.~\ref{fig:PII-Plus}, the fact that for jumps of type~(i), $\xi$ is bounded away from the points $\xi=0,\pm\ii\nu$, and the fact that the conjugating factors are uniformly bounded for such $\xi$ and independent of $y$, show that $\mathbf{E}_+(\xi;y)=\mathbf{E}_-(\xi;y)(\mathbb{I}+\text{exponentially small in $y$ and~$\xi$})$.

For jumps of type (ii), if we take the disk boundaries to be oriented in the clockwise direction, it is easy to see that $\mathbf{E}_+(\xi;y)= \mathbf{E}_-(\xi;y)\dot{\mathbf{Z}}^\mathrm{in}_{\xi_0}(\xi;y)\dot{\mathbf{Z}}^\mathrm{out}(\xi)^{-1}$ holds when $\xi\in\partial D(\xi_0)$, for all three cases $\xi_0=0,\pm\ii\nu$. For $\xi_0=\pm\ii\nu$ we have from \eqref{eq:Airy-mismatch} and its analogue on $\partial D(-\ii\nu)$ that $\mathbf{E}_+(\xi;y)=\mathbf{E}_-(\xi;y)\big(\mathbb{I}+O\big(y^{-3/2}\big)\big)$ holds uniformly on $\partial D(\pm\ii\nu)$. For $\xi_0=0$, we combine \eqref{eq:PC-Plus-Mismatch} with formul\ae\ \eqref{eq:alpha-def}--\eqref{eq:U-matrix-asymp-generalized} from the appendix and the fact that $W_0(\xi)$ is bounded away from zero on $\partial D(0)$ to conclude that
\begin{gather}
\mathbf{E}_+(\xi;y)=\mathbf{E}_-(\xi;y)\mathbf{H}_0(\xi)\left(\mathbb{I}+\frac{1}{W_0(\xi)}\begin{bmatrix}0 & c_{12}(y)\\
-c_{21}(y) & 0\end{bmatrix}\right. \label{eq:E-Jump-Plus-circle}\\ \left. {} + \begin{bmatrix}O\big(y^{-3/2}\big) & O(y^{-9/4 + 3\mathrm{Re}(\ii q)/2}) \\
O(y^{-9/4-3\mathrm{Re}(\ii q)/2}) & O\big(y^{-3/2}\big)\end{bmatrix}\right)\mathbf{H}_0(\xi)^{-1},\qquad\xi\in\partial D(0),\qquad y\to +\infty,\nonumber
\end{gather}
where we have used $\big|\ee^{\pm 2\ii y^{3/2}h(0+)}\big|=1$ (following from \eqref{eq:h-zero}) and where
\begin{gather}\begin{split}
& c_{12}(y):=\frac{1}{2\ii}\ee^{-2\ii y^{3/2}h(0+)}y^{-3/4+3\ii q/2}r(q,\mu)=\frac{\tau\sqrt{\pi}\ee^{-\ii\pi/4}\ee^{\pi q/2}2^{\ii q}}{\Gamma(\ii q)}\ee^{-\ii(2y/3)^{3/2}}y^{-3/4+3\ii q/2},\\
& c_{21}(y):=\frac{1}{2\ii}\ee^{2\ii y^{3/2}h(0+)}y^{-3/4-3\ii q/2}s(q,\mu)=\frac{q\Gamma(\ii q)\ee^{\ii\pi/4}\ee^{-\pi q/2}2^{-\ii q}}{2\tau\sqrt{\pi}}\ee^{\ii(2y/3)^{3/2}}\!y^{-3/4-3\ii q/2}.\end{split}\!\!\!\!\!\!
\label{eq:c12-c21}
\end{gather}
Now the inequalities \eqref{eq:ell-inequalities} and the definition \eqref{eq:p1-tau1} of $q$ show that $-\tfrac{1}{2}<\mathrm{Re}(\ii q)\le\tfrac{1}{2}$, which implies that all matrix elements in the error term on the second line of \eqref{eq:E-Jump-Plus-circle} are $O\big(y^{-3/2}\big)$ as $y\to +\infty$. However, the terms on the first line are larger; indeed, $c_{12}(y)$ does not decay to zero at all as $y\to+\infty$ if $\mathrm{Re}(\ii q)=\tfrac{1}{2}$. So it will be necessary to take that term into account explicitly in order to arrive at a small-norm problem, and it will also be convenient to remove the effect of $c_{21}(y)$, since for $-\tfrac{1}{2}<\mathrm{Re}(\ii q)\le\tfrac{1}{2}$ it can decay arbitrarily slowly as $y\to+\infty$ for $\mathrm{Re}(\ii q)$ near the lower end of the allowed range.

\subsubsection{First parametrix for the error}
Since $c_{12}(y)$ does not decay when $\mathrm{Re}(\ii q)=\tfrac{1}{2}$, to cover the whole range $-\tfrac{1}{2}<\mathrm{Re}(\ii q)\le\tfrac{1}{2}$ we must construct a parametrix for $\mathbf{E}(\xi;y)$ itself, which we will denote by $\dot{\mathbf{E}}(\xi;y)$. We want it to have the following properties. It should be analytic for $\xi\in\mathbb{C}\setminus\partial D(0)$, taking continuous boundary values on $\partial D(0)$ from the interior and exterior. It should converge to the identity as $\xi\to\infty$. Finally, it should satisfy the jump condition
\begin{gather*}
\dot{\mathbf{E}}_+(\xi;y)=\dot{\mathbf{E}}_-(\xi;y)\mathbf{H}_0(\xi)\left(\mathbb{I}+\frac{c_{12}(y)}{W_0(\xi)}\sigma_+\right)\mathbf{H}_0(\xi)^{-1},\qquad\xi\in\partial D(0),\qquad \sigma_+:=\begin{bmatrix}0&1\\0 & 0\end{bmatrix},
\end{gather*}
in which the circle is taken with clockwise orientation. Given that $W_0(\xi)$ has a simple pole at the origin, a reasonable ansatz for $\dot{\mathbf{E}}(\xi;y)$ is that its exterior boundary value is a function with a simple pole at $\xi=0$ and that satisfies the normalization condition at infinity:
\begin{gather}
\dot{\mathbf{E}}_+(\xi;y)=\mathbb{I}+\xi^{-1}\mathbf{R}(y).\label{eq:Edot-plus-ansatz}
\end{gather}
We then attempt to determine the residue $\mathbf{R}(y)$ by insisting that the corresponding boundary value $\dot{\mathbf{E}}_-(\xi;y)$ obtained from the jump condition admit analytic continuation to the full interior of the disk $D(0)$. Thus,
\begin{gather*}
\dot{\mathbf{E}}_-(\xi;y)=\dot{\mathbf{E}}_+(\xi;y)\mathbf{H}_0(\xi)\left(\mathbb{I}+\frac{c_{12}(y)}{W_0(\xi)}\sigma_+\right)^{-1}\mathbf{H}_0(\xi)^{-1} \\
\hphantom{\dot{\mathbf{E}}_-(\xi;y)}{} = \left(\mathbb{I}+\xi^{-1}\mathbf{R}(y)\right)\mathbf{H}_0(\xi)\left(\mathbb{I}-\frac{c_{12}(y)}{W_0(\xi)}\sigma_+\right)\mathbf{H}_0(\xi)^{-1}
\end{gather*}
needs to be analytic at $\xi=0$, where we have used $\sigma_+^2=\mathbf{0}$ and \eqref{eq:Edot-plus-ansatz} to obtain the second line. Noting that $\det(\mathbf{H}_0(0))=1$, the Laurent expansion of $\mathbf{H}_0(0)^{-1}\dot{\mathbf{E}}_-(\xi;y)\mathbf{H}_0(0)$ about $\xi=0$ is therefore
\begin{gather*}
\mathbf{H}_0(0)^{-1}\dot{\mathbf{E}}_-(\xi;y)\mathbf{H}_0(0)=-\frac{c_{12}(y)}{W_0'(0)}\mathbf{H}_0(0)^{-1}\mathbf{R}(y)\mathbf{H}_0(0)\sigma_+\xi^{-2} \\
\qquad{}+ \left(-\frac{c_{12}(y)}{W_0'(0)}\sigma_+-\frac{c_{12}(y)}{W_0'(0)}\mathbf{H}_0(0)^{-1}\mathbf{R}(y)\mathbf{H}_0'(0)\sigma_+\right.\\
 \qquad{} {}+\frac{c_{12}(y)}{W_0'(0)}\mathbf{H}_0(0)^{-1}\mathbf{R}(y)\mathbf{H}_0(0)\sigma_+\mathbf{H}_0(0)^{-1}\mathbf{H}_0'(0) \\
\left. \qquad{}+\mathbf{H}_0(0)^{-1}\mathbf{R}(y)\mathbf{H}_0(0)\left(\mathbb{I}+\frac{c_{12}(y)W_0''(0)}{2W_0'(0)^2}\sigma_+\right)\right)\xi^{-1}+O(1),\qquad\xi\to 0.
\end{gather*}
To remove the coefficient of $\xi^{-2}$, we will again use $\sigma_+^2=\mathbf{0}$ to express $\mathbf{R}(y)\mathbf{H}_0(0)$ in the rank-$1$ form
\begin{gather}
\mathbf{R}(y)\mathbf{H}_0(0) = \mathbf{a}(y)\mathbf{e}_2^\top,\qquad\mathbf{e}_2:=\begin{bmatrix}0\\1\end{bmatrix} \quad\implies\quad\mathbf{R}(y)\mathbf{H}_0(0)\sigma_+=\mathbf{0},\label{eq:Edot-plus-residue}
\end{gather}
and $\mathbf{a}(y)$ is a vector unknown yet to be determined. Using this form, two terms in the coefficient of $\xi^{-1}$ are automatically cancelled, and the condition that this coefficient also vanish becomes
\begin{gather*}
-\frac{c_{12}(y)}{W_0'(0)}\sigma_+-\frac{c_{12}(y)}{W_0'(0)}\mathbf{H}_0(0)^{-1}\mathbf{a}(y)\mathbf{e}_2^\top\mathbf{H}_0(0)^{-1}\mathbf{H}_0'(0)\sigma_++\mathbf{H}_0(0)^{-1}\mathbf{a}(y)\mathbf{e}_2^\top=\mathbf{0}.
\end{gather*}
The first column is an identity, and after multiplying on the left by $\mathbf{H}_0(0)$, the second column reads
\begin{gather*}
-\frac{c_{12}(y)}{W_0'(0)}\mathbf{H}_0(0)\mathbf{e}_1-\frac{c_{12}(y)}{W_0'(0)}\mathbf{a}(y)\mathbf{e}_2^\top
\mathbf{H}_0(0)^{-1}\mathbf{H}_0'(0)\mathbf{e}_1+\mathbf{a}(y)=\mathbf{0},\qquad\mathbf{e}_1:=\begin{bmatrix}1\\0\end{bmatrix}.
\end{gather*}
Since $\mathbf{e}_2^\top\mathbf{H}_0(0)^{-1}\mathbf{H}_0'(0)\mathbf{e}_1$ is a scalar, the solution is explicitly given by
\begin{gather}
\mathbf{a}(y):=\widetilde{c}_{12}(y)\mathbf{H}_0(0)\mathbf{e}_1,\qquad\widetilde{c}_{12}(y):=\frac{c_{12}(y)}{W_0'(0)-c_{12}(y)\mathbf{e}_2^\top\mathbf{H}_0(0)^{-1}\mathbf{H}_0'(0)\mathbf{e}_1},
\label{eq:Edot-plus-vector}
\end{gather}
provided the denominator of the coefficient $\widetilde{c}_{12}(y)$ is nonzero. Since $W_0'(0)>0$ and $\mathbf{H}_0(\xi)$ are independent of $y$, while $c_{12}(y)=O\big(y^{-3/4+3\mathrm{Re}(\ii q)/2}\big)$ as $y\to +\infty$, it is clear that $\widetilde{c}_{12}(y)$ will be well-defined and will also satisfy $\widetilde{c}_{12}(y)=O\big(y^{-3/4+3\mathrm{Re}(\ii q)/2}\big)$ if $\mathrm{Re}(\ii q)<\tfrac{1}{2}$. This is a decaying estimate as $y\to+\infty$. Only in the case $\mathrm{Re}(\ii q)=\tfrac{1}{2}$ is it even possible for $\dot{\mathbf{E}}(\xi;y)$ to fail to exist for sufficiently large positive~$y$. In the latter case, by excluding neighborhoods of certain points, we may still assume that $\widetilde{c}_{12}(y)=O(1)$ as $y\to +\infty$. Thus it is certainly true that $\dot{\mathbf{E}}(\xi;y)=O(1)$ as $y\to +\infty$ with the above caveat if $\mathrm{Re}(\ii q)=\tfrac{1}{2}$. Since $\det(\dot{\mathbf{E}}(\xi;y))=1$ the same holds for the inverse matrix.

Using the rational expression \eqref{eq:Edot-plus-ansatz} together with \eqref{eq:Edot-plus-residue} and \eqref{eq:Edot-plus-vector} and the identity \smash{$\mathbf{e}_1\mathbf{e}_2^\top\!=\!\sigma_+$} then shows that
\begin{gather}
\dot{\mathbf{E}}(\xi;y)=\mathbb{I}+\widetilde{c}_{12}(y)\mathbf{H}_0(0)\sigma_+\mathbf{H}_0(0)^{-1}\xi^{-1},\qquad\xi\in\mathbb{C}\setminus\overline{D(0)}. \label{eq:Edot-outside}
\end{gather}

\subsubsection{Second parametrix for the error}
We now compare the initial error matrix $\mathbf{E}(\xi;y)$ with its parametrix $\dot{\mathbf{E}}(\xi;y)$ by setting $\mathbf{F}(\xi;y):=\mathbf{E}(\xi;y)\dot{\mathbf{E}}(\xi;y)^{-1}$. Clearly, $\mathbf{F}(\xi;y)$ is analytic in the same domain as is $\mathbf{E}(\xi;y)$, and takes its boundary values on the jump contour in the same continuous fashion. Moreover, since we require $\mathbf{E}(\xi;y)\to\mathbb{I}$ as $\xi\to\infty$, we also have the corresponding condition on $\mathbf{F}(\xi;y)$ because $\dot{\mathbf{E}}(\xi;y)\to\mathbb{I}$ as $\xi\to\infty$. For all points $\xi$ on the jump contour for $\mathbf{E}$ omitting the circle $\partial D(0)$, we use the uniform estimate on $\dot{\mathbf{E}}(\xi;y)$ and its inverse to conclude that the previously obtained estimates of the deviation of the jump matrix for $\mathbf{E}(\xi;y)$ from the identity persist. Hence
\begin{gather}
\mathbf{F}_+(\xi;y)=\mathbf{F}_-(\xi;y)\big(\mathbb{I}+O\big(y^{-3/2}\big)\big),\qquad\xi\in\Sigma_\mathbf{F}\setminus\partial D(0),\label{eq:F-jump-outside}
\end{gather}
where the estimate is in the $L^\infty$ sense, but the error term also decays exponentially to zero in both $\xi\to\infty$ and $y>0$ along unbounded portions of $\Sigma_\mathbf{F}$. On the circle $\partial D(0)$ where $\mathbf{E}_+(\xi;y)=\mathbf{E}_-(\xi;y)\mathbf{V}^\mathbf{E}(\xi;y)$ and $\dot{\mathbf{E}}_+(\xi;y)=\dot{\mathbf{E}}_-(\xi;y)\mathbf{V}^{\dot{\mathbf{E}}}(\xi;y)$, we have $\mathbf{F}_+(\xi;y)=\mathbf{F}_-(\xi;y)\mathbf{V}^\mathbf{F}(\xi;y)$, where
\begin{gather*}
\mathbf{V}^\mathbf{F}(\xi;y)=\dot{\mathbf{E}}_+(\xi;y)\mathbf{V}^{\dot{\mathbf{E}}}(\xi;y)^{-1}\mathbf{V}^\mathbf{E}(\xi;y)\dot{\mathbf{E}}_+(\xi;y)^{-1},\qquad\xi\in\partial D(0).
\end{gather*}
Now, using $c_{12}(y)=O(1)$ along with $\sigma_+^2=\mathbf{0}$ and the fact that over the whole range $-\tfrac{1}{2}<\mathrm{Re}(\ii q)\le\tfrac{1}{2}$ we have $c_{12}(y)c_{21}(y)=O\big(y^{-3/2}\big)$ as $y\to +\infty$, we get
\begin{gather*}
\mathbf{V}^{\dot{\mathbf{E}}}(\xi;y)^{-1}\mathbf{V}^\mathbf{E}(\xi;y)=
\mathbf{H}_0(\xi)\left(\mathbb{I}-\frac{c_{12}(y)}{W_0(\xi)}\sigma_+\right)\\
\hphantom{\mathbf{V}^{\dot{\mathbf{E}}}(\xi;y)^{-1}\mathbf{V}^\mathbf{E}(\xi;y)=}{}\times
\left(\mathbb{I}+\frac{c_{12}(y)}{W_0(\xi)}\sigma_+-\frac{c_{21}(y)}{W_0(\xi)}\sigma_- + O\big(y^{-3/2}\big)\right)\mathbf{H}_0(\xi)^{-1}\\
\hphantom{\mathbf{V}^{\dot{\mathbf{E}}}(\xi;y)^{-1}\mathbf{V}^\mathbf{E}(\xi;y)}{} =
\mathbf{H}_0(\xi)\left(\mathbb{I}-\frac{c_{21}(y)}{W_0(\xi)}\sigma_-+O\big(y^{-3/2}\big)\right)\mathbf{H}_0(\xi)^{-1},\qquad\xi\in\partial D(0),
\end{gather*}
where $\sigma_-:=\sigma_+^\top$.
Therefore, using \eqref{eq:Edot-outside} and defining $\mathbf{B}(\xi):=\mathbf{H}_0(0)^{-1}\mathbf{H}_0(\xi)$,
\begin{gather*}
\mathbf{H}_0(0)^{-1}\mathbf{V}^\mathbf{F}(\xi;y)\mathbf{H}_0(0)\\
\qquad{} =\left(\mathbb{I}+\widetilde{c}_{12}(y)\sigma_+\xi^{-1}\right)\mathbf{B}(\xi)\left(\mathbb{I}-\frac{c_{21}(y)}{W_0(\xi)}\sigma_-+O\big(y^{-3/2}\big)\right)\mathbf{B}(\xi)^{-1}
\left(\mathbb{I}-\widetilde{c}_{12}(y)\sigma_+\xi^{-1}\right)\\
\qquad{} =\mathbb{I}-\frac{c_{21}(y)}{W_0(\xi)}\mathbf{B}(\xi)\sigma_-\mathbf{B}(\xi)^{-1} + O\big(y^{-3/2}\big),\qquad\xi\in\partial D(0),
\end{gather*}
where again we used $\sigma_+^2=\mathbf{0}$ and took into account that $\widetilde{c}_{12}(y)c_{21}(y)=O\big(y^{-3/2}\big)$ and that $\widetilde{c}_{12}(y)=O(1)$. Therefore conjugating again by $\mathbf{H}_0(0)$,
\begin{gather}
\mathbf{V}^\mathbf{F}(\xi;y)=\mathbb{I}-\frac{c_{21}(y)}{W_0(\xi)}\mathbf{H}_0(\xi)\sigma_-\mathbf{H}_0(\xi)^{-1}+O\big(y^{-3/2}\big),\qquad\xi\in\partial D(0),\qquad y\to +\infty.\label{eq:F-jump-circle}
\end{gather}
This form suggests that we might finally arrive at a small-norm target problem with jump matrices uniformly deviating from the identity by $O\big(y^{-3/2}\big)$ if we choose to remove the term proportional to $c_{21}(y)$ by a second parametrix for the error, which we denote by $\dot{\mathbf{F}}(\xi;y)$. Technically, since $c_{21}(y)=O\big(y^{-3/4-3\mathrm{Re}(\ii q)/2}\big)$ which is $o(1)$ as $y\to +\infty$ over the whole range of values $-\tfrac{1}{2}<\mathrm{Re}(\ii q)\le\tfrac{1}{2}$, $\mathbf{F}(\xi;y)$ already satisfies the conditions of a small-norm problem, but it is convenient to first construct $\dot{\mathbf{F}}(\xi;y)$ and then compare to $\mathbf{F}(\xi;y)$ before obtaining estimates from small-norm theory.

Therefore, we define $\dot{\mathbf{F}}(\xi;y)$ as the matrix analytic for $\xi\in\mathbb{C}\setminus\partial D(0)$, normalized to the identity as $\xi\to\infty$, and whose continuous boundary values on $\partial D(0)$ are related by the jump condition
\begin{gather}
\dot{\mathbf{F}}_+(\xi;y)=\dot{\mathbf{F}}_-(\xi;y)\left(\mathbb{I}-\frac{c_{21}(y)}{W_0(\xi)}\mathbf{H}_0(\xi)\sigma_-\mathbf{H}_0(\xi)^{-1}\right),\qquad\xi\in\partial D(0). \label{eq:dot-F-jump}
\end{gather}
Of course, the solution of this problem is very similar to that for $\dot{\mathbf{E}}(\xi;y)$, and we obtain the result that
\begin{gather}
\dot{\mathbf{F}}(\xi;y)=\mathbb{I}-\widetilde{c}_{21}(y)\mathbf{H}_0(0)\sigma_-\mathbf{H}_0(0)^{-1}\xi^{-1},\nonumber\\
\widetilde{c}_{21}(y):=\frac{c_{21}(y)}{W_0'(0)+c_{21}(y)\mathbf{e}_1^\top\mathbf{H}_0(0)^{-1}\mathbf{H}_0'(0)\mathbf{e}_2},\qquad\xi\in\mathbb{C}\setminus\overline{D(0)},\label{eq:Fdot-outside}
\end{gather}
and one can check directly that $\dot{\mathbf{F}}_-(\xi;y)$, obtained from evaluating the above expression on $\partial D(0)$ as $\dot{\mathbf{F}}_+(\xi;y)$ and applying the jump condition~\eqref{eq:dot-F-jump}, admits analytic continuation to the interior of~$D(0)$. Unlike $\widetilde{c}_{12}(y)$ for which existence may require conditions on arbitrarily large $y>0$ if $\mathrm{Re}(\ii q)=\tfrac{1}{2}$, $\widetilde{c}_{21}(y)$ always exists as long as $y>0$ is sufficiently large and $-\tfrac{1}{2}<\mathrm{Re}(\ii q)\le\tfrac{1}{2}$.

\subsubsection{Improved error analysis}
We finally set
\begin{gather*}
\mathbf{G}(\xi;y):=\mathbf{F}(\xi;y)\dot{\mathbf{F}}(\xi;y)^{-1}.
\end{gather*}
Carrying forward the conditions on $y\to +\infty$ sufficient to guarantee that $\widetilde{c}_{12}(y)=O(1)$ as $y\to +\infty$ (no conditions needed unless $\mathrm{Re}(\ii q)=\tfrac{1}{2}$ in which case it may be necessary to exclude certain intervals), we can show that $\mathbf{G}(\xi;y)$ satisfies the conditions of a small-norm Riemann--Hilbert problem. Clearly, $\mathbf{G}(\xi;y)$ is analytic in the domain $\xi\in\mathbb{C}\setminus\Sigma_\mathbf{G}$, where $\Sigma_\mathbf{G}=\Sigma_\mathbf{F}=\Sigma_\mathbf{E}$, and $\mathbf{G}(\xi;y)\to\mathbb{I}$ as $\xi\to\infty$ as this normalization holds for both factors in the definition. It remains to estimate the difference of the jump matrix for $\mathbf{G}$ from the identity. First note that since $\dot{\mathbf{F}}(\xi;y)$ and $\dot{\mathbf{F}}(\xi;y)^{-1}$ are bounded uniformly for $y>0$ sufficiently large, it is easy to see that the estimate~\eqref{eq:F-jump-outside} implies a similar estimate for~$\mathbf{G}(\xi;y)$:
\begin{gather}
\mathbf{G}_+(\xi;y)=\mathbf{G}_-(\xi;y)\big(\mathbb{I}+O\big(y^{-3/2}\big)\big),\qquad\xi\in\Sigma_\mathbf{G}\setminus\partial D(0) \label{eq:G-jump}
\end{gather}
again with the estimate holding uniformly and also exhibiting exponential decay in both $y>0$ and $|\xi|$ on unbounded arcs of $\Sigma_\mathbf{G}$. Because the jump condition \eqref{eq:dot-F-jump} for $\dot{\mathbf{F}}(\xi;y)$ on $\partial D(0)$ takes into account everything except the $O\big(y^{-3/2}\big)$ error term in the jump matrix \eqref{eq:F-jump-circle} for $\mathbf{F}(\xi;y)$, it is easy to show that the estimate \eqref{eq:G-jump} extends also to $\xi\in\partial D(0)$. Through standard analysis of a system of singular integral equations equivalent to the Riemann--Hilbert conditions for $\mathbf{G}(\xi;y)$ (completely analogous to \eqref{eq:singular-integral-eqn}) it follows that $\mathbf{G}(\xi;y)$ exists for $y>0$ sufficiently large (again with possible exclusions if $\mathrm{Re}(\ii q)=\tfrac{1}{2}$) and satisfies
\begin{gather}
\mathbf{G}(\xi;y)=\mathbb{I}+\mathbf{G}^1(y)\xi^{-1}+\mathbf{G}^2(y)\xi^{-2}+O\big(\xi^{-3}\big),\qquad\xi\to\infty, \label{eq:G-expand}
\end{gather}
where the limit $\xi\to\infty$ may be taken in any direction, and in which the moments $\mathbf{G}^k(y)$, $k=1,2$, satisfy the estimate $\mathbf{G}^k(y)=O\big(y^{-3/2}\big)$ as $y\to+\infty$.

Now $\mathbf{Z}(\xi;y)=\mathbf{E}(\xi;y)\dot{\mathbf{Z}}(\xi;y)=\mathbf{F}(\xi;y)\dot{\mathbf{E}}(\xi;y)\dot{\mathbf{Z}}(\xi;y)=\mathbf{G}(\xi;y)\dot{\mathbf{F}}(\xi;y)\dot{\mathbf{E}}(\xi;y)\dot{\mathbf{Z}}(\xi;y)$. Provided that $\xi=\zeta y^{-1/2}$ is sufficiently large, we may take the rational forms \eqref{eq:Edot-outside} and \eqref{eq:Fdot-outside} for $\dot{\mathbf{E}}(\xi;y)$ and $\dot{\mathbf{F}}(\xi;y)$ respectively, and also we will have $\dot{\mathbf{Z}}(\xi;y)=\dot{\mathbf{Z}}^\mathrm{out}(\xi)$. Therefore, for such $\xi$, we have the exact formula
\begin{gather*}
\mathbf{W}(\zeta;y)\zeta^{\ii p\sigma_3}=y^{-\ii p\sigma_3/2}\mathbf{Y}\big(y^{-1/2}\zeta;y\big)\zeta^{\ii p\sigma_3}\\
\hphantom{\mathbf{W}(\zeta;y)\zeta^{\ii p\sigma_3}}{} =y^{-\ii p\sigma_3/2}\mathbf{Z}\big(y^{-1/2}\zeta;y\big)\big(y^{-1/2}\zeta\big)^{\ii p\sigma_3}\ee^{\ii y^{3/2}g(y^{-1/2}\zeta)\sigma_3}y^{\ii p\sigma_3/2}\\
\hphantom{\mathbf{W}(\zeta;y)\zeta^{\ii p\sigma_3}}{} =y^{-\ii p\sigma_3/2}\mathbf{G}\big(y^{-1/2}\zeta;y\big)\dot{\mathbf{F}}\big(y^{-1/2}\zeta;y\big)\dot{\mathbf{E}}\big(y^{-1/2}\zeta;y\big)\mathbf{K}\big(y^{-1/2}\zeta\big)\\
\hphantom{\mathbf{W}(\zeta;y)\zeta^{\ii p\sigma_3}=}{} \times\ee^{j(y^{-1/2}\zeta)\sigma_3}\ee^{\ii y^{3/2}g(y^{-1/2}\zeta)\sigma_3}y^{\ii p\sigma_3/2}.
\end{gather*}
This expression admits an asymptotic expansion in descending non-negative powers of $\zeta^{-1}$ as $\zeta\to\infty$, i.e., of the form of the expansion in parentheses on the first line of \eqref{eq:large-zeta-PII}. It is straightforward but tedious to calculate explicitly the first two coefficients $\mathbf{W}^1(y)$ and $\mathbf{W}^2(y)$ by combi\-ning~\eqref{eq:G-expand}, the rational expressions \eqref{eq:Edot-outside} and \eqref{eq:Fdot-outside}, the large-$\xi$ expansion of $\mathbf{K}(\xi)$ (see~\eqref{eq:K-plus-def}):
\begin{gather*}
\mathbf{K}(\xi)=\mathbb{I} + \begin{bmatrix}0 & 6^{\ii p-1/2}\\-6^{-\ii p-1/2} & 0\end{bmatrix}
\xi^{-1}-\frac{1}{12}\mathbb{I}\xi^{-2}+O\big(\xi^{-3}\big),\qquad\xi\to\infty,
\end{gather*}
and the expansion (see \eqref{eq:g-asymp} and \eqref{eq:j-large-xi})
\begin{gather*}
\ee^{j(\xi)\sigma_3}\ee^{\ii y^{3/2}g(\xi)\sigma_3}=\mathbb{I}+\left(
-\sqrt{\frac{2}{3}} q
+\frac{1}{12} y^{3/2}\right)\ii\sigma_3\xi^{-1} \\
\qquad{}+\left(-\frac{1}{6}\ii p\sigma_3
-\frac{1}{3}q^2\mathbb{I}
+\frac{1}{12}\sqrt{\frac{2}{3}}qy^{3/2}\mathbb{I}
-\frac{1}{288}y^3\mathbb{I}\right)\xi^{-2}+O\big(\xi^{-3}\big),\qquad\xi\to\infty.
\end{gather*}
In substituting from the formul\ae\ \eqref{eq:Edot-outside} and \eqref{eq:Fdot-outside}, we use the fact
\begin{gather*}
\mathbf{H}_0(0)\sigma_\pm\mathbf{H}_0(0)^{-1}=\frac{1}{2}W_0'(0)^{\pm 2\ii q}\ee^{\pm(2\ii\eta- 2k(0))}\begin{bmatrix}1 & \mp\ee^{2\ii\eta}\\\pm\ee^{-2\ii\eta} & -1\end{bmatrix}\\
\hphantom{\mathbf{H}_0(0)\sigma_\pm\mathbf{H}_0(0)^{-1}}{} =\frac{1}{2}\left(\frac{32}{3}\right)^{\pm \ii q/2}\begin{bmatrix}1 & \mp 6^{\ii p}\\\pm 6^{-\ii p} & -1\end{bmatrix},
\end{gather*}
which follows from \eqref{eq:H00} with the help of \eqref{eq:eta-value}, \eqref{eq:k-zero}, and \eqref{eq:W0-prime}, recalling also \eqref{eq:nu-define}. This leads to the following exact formul\ae\ for $\mathcal{V}(y)$ and $\mathcal{Q}(y)$:
\begin{gather}
\mathcal{V}(y)=W^1_{21}(y)=y^{\ii p+1/2}\left(-6^{-(\ii p+1/2)} + \frac{1}{2}\left(\frac{32}{3}\right)^{\ii q/2}6^{-\ii p}\widetilde{c}_{12}(y)\right.\nonumber\\
\left. \hphantom{\mathcal{V}(y)=W^1_{21}(y)=}{} +\frac{1}{2}\left(\frac{32}{3}\right)^{-\ii q/2}6^{-\ii p}\widetilde{c}_{21}(y) + G^1_{21}(y)\right),\label{eq:V-basic}
\end{gather}
and
\begin{gather*}
\mathcal{Q}(y)=\ii W^1_{11}(y)-\ii\frac{W^2_{21}(y)}{W^1_{21}(y)}\\
\hphantom{\mathcal{Q}(y)}{} =\ii y^{1/2}\left(\frac{1}{2}\left(\frac{32}{3}\right)^{\ii q/2}\widetilde{c}_{12}(y)-\frac{1}{2}\left(\frac{32}{3}\right)^{-\ii q/2}\widetilde{c}_{21}(y)+G_{11}^1(y)\right)-\ii y^{\ii p+1}\frac{\mathcal{N}(y)}{\mathcal{V}(y)},
\end{gather*}
where
\begin{gather}
\mathcal{N}(y):=\frac{1}{2}\left(\frac{32}{3}\right)^{\ii q/2}6^{-\ii p-1/2}\widetilde{c}_{12}(y) -\frac{1}{2}\left(\frac{32}{3}\right)^{-\ii q/2}6^{-\ii p-1/2}\widetilde{c}_{21}(y)+\frac{1}{2}6^{-\ii p}\widetilde{c}_{12}(y)\widetilde{c}_{21}(y)\nonumber\\
\hphantom{\mathcal{N}(y):=}{}+G^2_{21}(y)-6^{-\ii p-1/2}G^1_{22}(y) +\frac{1}{2}\left(\frac{32}{3}\right)^{\ii q/2}\widetilde{c}_{12}(y)\big(G^1_{21}(y)+6^{-\ii p}G_{22}^1(y)\big) \nonumber\\
\hphantom{\mathcal{N}(y):=}{} -\frac{1}{2}\left(\frac{32}{3}\right)^{-\ii q/2}\widetilde{c}_{21}(y)\big(G^1_{21}(y)-6^{-\ii p}G_{22}^1(y)\big).
\label{eq:N-formula}
\end{gather}

Suppose now that $-\tfrac{1}{2}<\mathrm{Re}(\ii q)<\tfrac{1}{2}$, so we are excluding only the endpoint case $\mathrm{Re}(\ii q)=\tfrac{1}{2}$. Then from \eqref{eq:c12-c21} we see that $c_{12}(y)$ and $c_{21}(y)$ both decay algebraically as $y\to +\infty$, so from \eqref{eq:Edot-plus-vector} and \eqref{eq:Fdot-outside} we get (also using \eqref{eq:W0-prime}) that
\begin{gather}
\widetilde{c}_{12}(y)=\left(\frac{2}{3}\right)^{1/4}c_{12}(y)+O\big(c_{12}(y)^2\big)\qquad\text{and}\nonumber\\
\widetilde{c}_{21}(y)=\left(\frac{2}{3}\right)^{1/4}c_{21}(y)+O\big(c_{21}(y)^2\big),\qquad y\to +\infty.
\label{eq:tildes-no-tildes}
\end{gather}
Furthermore, this condition on $q$ guarantees that we are in the small-norm setting for $\mathbf{G}(\xi;y)$, so $\mathbf{G}^k(y)=O\big(y^{-3/2}\big)$ holds for $k=1,2$.
Then, using \eqref{eq:c12-c21} again we see that \eqref{eq:V-basic} becomes simply
\begin{gather}
\mathcal{V}(y)=-\left(\frac{y}{6}\right)^{\ii p+1/2}\big(1+O\big(y^{-3/4+3|\mathrm{Re}(\ii q)|/2}\big)\big)\nonumber\\
\hphantom{\mathcal{V}(y)}{} =-\left(\frac{y}{6}\right)^{\ii p+1/2}(1+o(1)),\qquad y\to +\infty.\label{eq:V-asymp-y-plus}
\end{gather}
This becomes \eqref{eq:V-asymp-y-Plus-specific} when $p=\ln(2)/(2\pi)$, which means $\mathrm{Re}(\ii q)=0$.
Similarly, all but the terms on the first line of \eqref{eq:N-formula} are $O\big(y^{-3/2}\big)$, so using \eqref{eq:c12-c21} and \eqref{eq:tildes-no-tildes},
\begin{gather*}
\mathcal{N}(y)=24^{-1/4}6^{-\ii p-1/2}y^{-3/4}\\
\quad{}\times\left[\frac{\tau\sqrt{\pi}\ee^{\pi q/2}}{\Gamma(\ii q)}\ee^{-\ii \pi/4}\left(\frac{128}{3}\right)^{\ii q/2}\ee^{-\ii(2y/3)^{3/2}}y^{3\ii q/2}\big(1+O\big(y^{-3/4+3\mathrm{Re}(\ii q)/2}\big)\big)\right.\\
\left. \quad{} -\frac{q\Gamma(\ii q)\ee^{-\pi q/2}}{2\tau\sqrt{\pi}}
\ee^{\ii\pi/4}\left(\frac{128}{3}\right)^{-\ii q/2}\ee^{\ii(2y/3)^{3/2}}y^{-3\ii q/2}\big(1+O\big(y^{-3/4-3\mathrm{Re}(\ii q)/2}\big)\big)\right]
+O\big(y^{-3/2}\big)\\
=24^{-1/4}6^{-\ii p-1/2}y^{-3/4}\left[\frac{\tau\sqrt{\pi}\ee^{\pi q/2}}{\Gamma(\ii q)}\ee^{-\ii \pi/4}\left(\frac{128}{3}\right)^{\ii q/2}\ee^{-\ii(2y/3)^{3/2}}y^{3\ii q/2}\right.\\
\left.\quad{} -\frac{q\Gamma(\ii q)\ee^{-\pi q/2}}{2\tau\sqrt{\pi}}
\ee^{\ii\pi/4}\left(\frac{128}{3}\right)^{-\ii q/2}\ee^{\ii(2y/3)^{3/2}}y^{-3\ii q/2}\right]+O\big(y^{-3/2+3|\mathrm{Re}(\ii q)|}\big),\qquad y\to +\infty.
\end{gather*}
The two terms in square brackets are only of equal size in the limit $y\to+\infty$ if $\mathrm{Re}(\ii q)=0$. Thus, keeping only the leading term(s), there are three different asymptotic formul\ae\ for $\mathcal{Q}(y)$ in the limit $y\to +\infty$ assuming that $-\tfrac{1}{2}<\mathrm{Re}(\ii q)<\tfrac{1}{2}$:
\begin{gather*}
\mathcal{Q}(y)=\frac{q\Gamma(\ii q)\ee^{-\pi q/2}}{24^{1/4}\tau\sqrt{\pi}}\ee^{-\ii\pi/4}\left(\frac{128}{3}\right)^{-\ii q/2}\ee^{\ii(2y/3)^{3/2}}y^{-1/4-3\ii q/2}\big(1+O\big(y^{M^-(q)}\big)\big),\\ y\to+\infty,\qquad -\frac{1}{2}<\mathrm{Re}(\ii q)<0,\qquad M^-(q):=\max\left\{-\frac{3}{4}-\frac{3}{2}\mathrm{Re}(\ii q),3\mathrm{Re}(\ii q)\right\}<0,\\
\mathcal{Q}(y)=\frac{2\tau\sqrt{\pi}\ee^{\pi q/2}}{24^{1/4}\Gamma(\ii q)}\ee^{\ii\pi/4}\left(\frac{128}{3}\right)^{\ii q/2}\ee^{-\ii(2y/3)^{3/2}}y^{-1/4+3\ii q/2}\big(1+O\big(y^{M^+(q)}\big)\big),\\ y\to+\infty,\qquad 0<\mathrm{Re}(\ii q)<\frac{1}{2},\qquad M^+(q):=\max\left\{-\frac{3}{4}+\frac{3}{2}\mathrm{Re}(\ii q),-3\mathrm{Re}(\ii q)\right\}<0,
\end{gather*}
and, using the identity \cite[equation~(5.4.3)]{DLMF} and the relations $\mu=\tau^{-1}$ and $\mu^2=\ee^{2\pi q}-1$,
\begin{gather*}
\mathcal{Q}(y)=-\frac{2\tau\sqrt{2q(\ee^{2\pi q}-1)}}{(24y)^{1/4}}\sin(\Theta(y)) + O\big(y^{-1}\big),\\
y\to+\infty,\qquad\!\Theta(y):=-\left(\frac{2}{3}y\right)^{3/2}\!+\frac{3}{2}q\ln(y)-\frac{1}{4}\pi+\frac{1}{2}q\ln\left(\frac{128}{3}\right)-\arg(\Gamma(\ii q)),\!\!\qquad q\in\mathbb{R}.
\end{gather*}
It is straightforward to translate these formul\ae\ into corresponding the asymptotic formul\ae\ for $u^-_\mathrm{TT}(x;\alpha)=-\big(\tfrac{3}{2}\big)^{1/3}\mathcal{Q}\big({-}\big(\tfrac{3}{2}\big)^{1/3}x\big)$ given in the statement of Theorem~\ref{thm:centerline}, where $\alpha=\tfrac{1}{2}+\ii p$. In doing so, one must express $q$ and $\tau$ explicitly in terms of $\alpha$, or equivalently, $p$, which results in the definitions \eqref{eq:q0}--\eqref{eq:q-tau}. This concludes the proof of Theorem~\ref{thm:centerline}.

\section{Proof of Theorem~\ref{thm:integral}}\label{sec:integral}
The computation of total integrals of Painlev\'e-II solutions $u(x;\alpha)$ is fairly straightforward when such solutions are extracted from a Riemann--Hilbert problem associated with the Jimbo--Miwa Lax pair because the ``fundamental'' potentials in this setting are the quantities $\mathcal{U}(y)$ and $\mathcal{V}(y)$ (see \eqref{eq:UV-define}) whose asymptotics are easiest to compute for large $|y|$ and whose logarithmic derivatives yield solutions of Painlev\'e-II.

Since $\mathcal{V}'(y)/\mathcal{V}(y)=\mathcal{Q}(y)$ and since \eqref{eq:Q-asymp-sector} gives asymptotics for $\mathcal{Q}(y)$ as $y\to-\infty$ up to an error term that it absolutely integrable in this limit, it is easy to see that
\begin{gather*}
\mathcal{V}(y)=K\frac{\ee^{-2\ii (-y/3)^{3/2}}}{(-3y)^{1/4+\ii p/2}}\exp\left(\int_{-\infty}^y
\left[\mathcal{Q}(\eta)-\ii\sqrt{-\frac{\eta}{3}}+\left(\frac{1}{4}+\ii\frac{p}{2}\right)\frac{1}{\eta}\right]\,\dd\eta\right)
\end{gather*}
holds for sufficiently negative $y$, where $K$ is an integration constant, and where the integral in the exponent is absolutely convergent. The integration constant $K$ is easily determined from the asymptotic formula \eqref{eq:V-asymp-y-Minus} valid for $\mathcal{V}(y)$ as $y\to -\infty$. Therefore we have the representation
\begin{gather}
\mathcal{V}(y)=\frac{\tau p\Gamma(\ii p)\ee^{-3\pi\ii/4}}{2^{1+\ii p}\sqrt{\pi}\ee^{\pi p/2}} \frac{\ee^{-2\ii(-y/3)^{3/2}}}{(-3y)^{1/4+\ii p/2}}\nonumber\\
\hphantom{\mathcal{V}(y)=}{}\times \exp\left(\int_{-\infty}^y
\left[\mathcal{Q}(\eta)-\ii\sqrt{-\frac{\eta}{3}}+\left(\frac{1}{4}+\ii\frac{p}{2}\right)\frac{1}{\eta}\right]\,\dd\eta\right), \label{eq:V-Q-left}
\end{gather}
which is valid for $y<0$ of sufficiently large magnitude that it lies to the left of all real poles of~$\mathcal{Q}(y)$.

Pick $a<0$ and $b>0$ so that all real poles of $\mathcal{Q}$ lie in the interval $(a,b)$. This is possible under the hypothesis $\mathrm{Re}(\ii q)\neq\tfrac{1}{2}$ according to the defining asymptotic formula \eqref{eq:u-x-wide-asymp} and Theorem~\ref{thm:centerline}. Then using the formula \eqref{eq:V-Q-left} for $y=a$ and again the relation $\mathcal{V}'(y)/\mathcal{V}(y)=\mathcal{Q}(y)$, we get the following formula, now valid for $y>b$:
\begin{gather}
\mathcal{V}(y)=\frac{\tau p\Gamma(\ii p)\ee^{-3\pi\ii/4}}{2^{1+\ii p}\sqrt{\pi}\ee^{\pi p/2}}
\frac{\ee^{-2\ii(-a/3)^{3/2}}}{(-3a)^{1/4+\ii p/2}}\exp\left(\int_{-\infty}^{a}\left[\mathcal{Q}(\eta)-\ii\sqrt{-\frac{\eta}{3}}+\left(\frac{1}{4}+\ii\frac{p}{2}\right)\frac{1}{\eta}\right]\,\dd\eta \right.
\nonumber\\
\left. \hphantom{\mathcal{V}(y)=}{} + \int_{C(a,b)}\mathcal{Q}(\eta)\,\dd\eta+\int_b^y\mathcal{Q}(\eta)\,\dd\eta\right).\label{eq:V-formula-last}
\end{gather}
Here, $C(a,b)$ is any contour in the complex plane from $a$ to $b$ that avoids all poles of $\mathcal{Q}$. Now rewriting \eqref{eq:V-asymp-y-plus} for $y>b$ in the form
\begin{gather*}
\mathcal{V}(y)=-\left(\frac{b}{6}\right)^{1/2+\ii p}\exp\left(\int_b^y\left(\frac{1}{2}+\ii p\right)\frac{\dd\eta}{\eta}\right)(1+o(1)),\\
 y\to +\infty,\qquad-\frac{1}{2}<\mathrm{Re}(\ii q)<\frac{1}{2},
\end{gather*}
we compare with \eqref{eq:V-formula-last} and eliminate $\mathcal{V}(y)$ to find
\begin{gather}
\exp\left(\int_{-\infty}^{a}\left[\mathcal{Q}(\eta)-\ii\sqrt{-\frac{\eta}{3}}+\left(\frac{1}{4}+\ii\frac{p}{2}\right)\frac{1}{\eta}\right]\,\dd\eta \right.\nonumber\\
\left.\qquad\quad{}+ \int_{C(a,b)}\mathcal{Q}(\eta)\,\dd\eta+ \int_b^y\left[\mathcal{Q}(\eta)-\left(\frac{1}{2}+\ii p\right)\frac{1}{\eta}\right]\,\dd\eta\right) \nonumber\\
\qquad{} = \frac{\sqrt{2\pi}\ee^{2\ii(-a/3)^{3/2}}}{\tau p\Gamma(\ii p)}\ee^{-\ii\pi(1/2+\ii p)/2}\left(b\sqrt{-\frac{a}{3}}\right)^{1/2+\ii p}(1+o(1)),\nonumber\\
\qquad \quad {} \ y\to +\infty, \quad -\frac{1}{2}\mathrm{Re}(\ii q)<\frac{1}{2}.\label{eq:Almost-total-integral}
\end{gather}

Taking the limit $y\to +\infty$ in \eqref{eq:Almost-total-integral} yields the total integral identity
\begin{gather*}
\exp\left(\int_{-\infty}^{a}\left[\mathcal{Q}(y)-\ii\sqrt{-\frac{y}{3}}+\left(\frac{1}{4}+\ii\frac{p}{2}\right)\frac{1}{y}\right]\,\dd y \right.\nonumber\\
\left.\qquad \quad {}+ \int_{C(a,b)}\mathcal{Q}(y)\,\dd y+\int_b^{+\infty}\left[\mathcal{Q}(y)-\left(\frac{1}{2}+\ii p\right)\frac{1}{y}\right]\,\dd y\right) \nonumber\\
\qquad{}= \frac{\sqrt{2\pi}\ee^{2\ii(-a/3)^{3/2}}}{\tau p\Gamma(\ii p)}\ee^{-\ii\pi(1/2+\ii p)/2}\left(b\sqrt{-\frac{a}{3}}\right)^{1/2+\ii p},\qquad -\frac{1}{2}<\mathrm{Re}(\ii q)<\frac{1}{2},
%\label{eq:total-integral}
\end{gather*}
in which the integral over the interval $(b,+\infty)$ is not generally absolutely convergent, but it necessarily makes sense as an improper integral. Suppose that we take for $C(a,b)$ a real path with infinitesimal semicircular indentations in the lower half-plane centered at each real pole of~$\mathcal{Q}$. It is well-known that all poles of $\mathcal{Q}$ are simple and have residue $\pm 1$. Therefore,
\begin{gather*}
\int_{C(a,b)}\mathcal{Q}(y)\,\dd y = \mathrm{P.V.}\int_a^b\mathcal{Q}(y)\,\dd y + \ii\pi (N_+-N_-),
\end{gather*}
where ``P.V.'' denotes the Hadamard principal value and $N_\pm$ is the number of real poles of $\mathcal{Q}$ of residue $\pm 1$.

Setting $A:=-(\tfrac{2}{3})^{1/3}b<0$ and $B:=-(\tfrac{2}{3})^{1/3}a>0$ so all real poles of $u^-_\mathrm{TT}(x;\alpha)=-\big(\tfrac{3}{2}\big)^{1/3}\mathcal{Q}(-\big(\tfrac{3}{2}\big)^{1/3}x)$ lie in the interval $(A,B)$, we have established the total integral formu\-la~\eqref{eq:u-total-integral} for the increasing tritronqu\'ee solution $u=u^-_\mathrm{TT}(x;\alpha)$ of \eqref{eq:PII-intro}. This completes the proof of Theorem~\ref{thm:integral}.

\section{Proof of Theorem~\ref{thm:global}}\label{sec:global}
\begin{figure}[t]\centering
\includegraphics{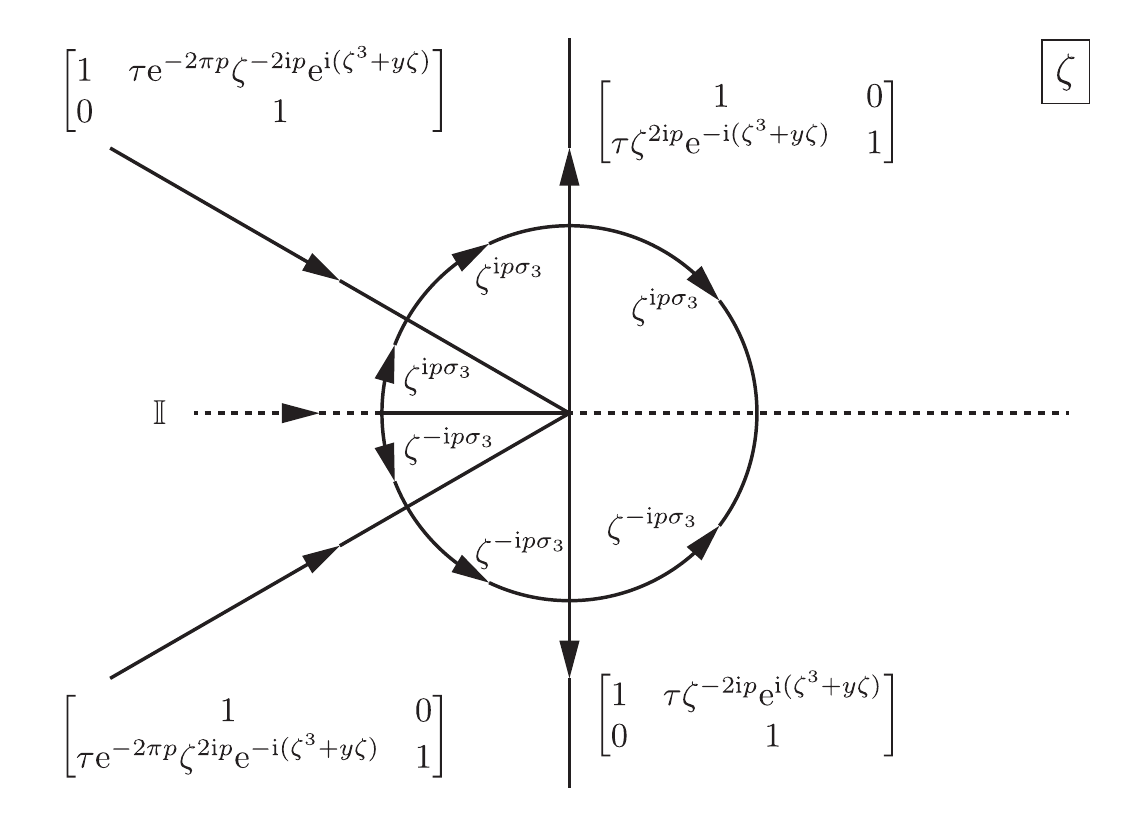}
\caption{The jump contour and jump matrix $\widetilde{\mathbf{V}}^\mathrm{PII}(\zeta;y)$ for $\widetilde{\mathbf{W}}(\zeta;y)$. On the segments of the jump contour within the unit disk, the jump matrix is exactly as shown in Fig.~\ref{fig:PII-jumps-generalized}.}\label{fig:PII-VL}
\end{figure}

When $p>0$, which also implies that $\tau\in\mathbb{R}$, Riemann--Hilbert Problem~\ref{rhp:PII-generalized} has a unique solution for every $y\in\mathbb{R}$ as a consequence of Zhou's vanishing lemma \cite{Zhou89}. Indeed, let us define a modified unknown $\widetilde{\mathbf{W}}$ as $\widetilde{\mathbf{W}}:=\mathbf{W}$ for $|\zeta|<1$ and $\widetilde{\mathbf{W}}:=\mathbf{W}\zeta^{\ii p\sigma_3}$ for $|\zeta|>1$. This transformation removes the jump from the negative real axis for $|\zeta|>1$, modifies the remaining jump matrices for $|\zeta|>1$ by conjugation, and introduces a new diagonal jump across the unit circle. For the latter, it is convenient in the setting of \cite{Zhou89} to take the upper and lower semicircles both to be oriented from left to right. The jump conditions for $\widetilde{\mathbf{W}}(\zeta;y)$ are illustrated in Fig.~\ref{fig:PII-VL}.
Now $\widetilde{\mathbf{W}}(\zeta;y)$ satisfies the conditions of an identity-normalized (as $\zeta\to\infty$) Riemann--Hilbert problem whose jump matrix is easily checked to satisfy the cyclic condition at each self-intersection point that is necessary for consistency. Moreover, when $p>0$ and $y\in\mathbb{R}$, the jump matrix satisfies $\widetilde{\mathbf{V}}^\mathrm{PII}(\zeta;y)=\widetilde{\mathbf{V}}^\mathrm{PII}(\zeta^*;y)^\dagger$ for all non-real $\zeta$ and that $\widetilde{\mathbf{V}}^\mathrm{PII}(\zeta;y)+\widetilde{\mathbf{V}}^\mathrm{PII}(\zeta;y)^\dagger$ is positive-definite for all real $\zeta$ in the jump contour. The vanishing lemma~\cite{Zhou89} therefore guarantees the existence of $\widetilde{\mathbf{W}}(\zeta;y)$, which also yields $\mathbf{W}(\zeta;y)$ by inverting the substitution made for $|\zeta|>1$. Uniqueness of the solution is standard, as is the fact that the solution satisfies $\det(\mathbf{W}(\zeta;y))=1$.

The existence of a solution of Riemann--Hilbert Problem~\ref{rhp:PII-generalized} for every real $y$ when $p>0$ actually implies via analytic Fredholm theory that the functions $\mathcal{U}(y)$ and $\mathcal{V}(y)$ defined by~\eqref{eq:UV-define} are analytic for $y\in\mathbb{R}$. Moreover, these functions are Schwarz reflections of each other: $\mathcal{V}(y)=-\mathcal{U}(y^*)^*$. The relation \eqref{eq:PII-COM} then implies that $\mathcal{U}(y)$ and $\mathcal{V}(y)$ cannot have any real zeros unless $p=0$ (in which case they vanish identically because $\mathbf{W}(\zeta;y)\equiv \mathbb{I}$), because $\mathcal{U}(y)$ and $\mathcal{V}(y)$ would have to vanish simultaneously. The same argument shows that $\mathcal{U}(y)$ and $\mathcal{V}(y)$ cannot have any real critical points for positive real $p$. Therefore, when $p>0$, $\mathcal{P}(y)$ and $\mathcal{Q}(y)$ defined by \eqref{eq:PQ-define} are nonvanishing analytic functions for $y\in\mathbb{R}$, and they are related by the symmetry $\mathcal{Q}(y)=\mathcal{P}(y^*)^*$. This implies that when $p>0$, the increasing tritronqu\'ee solution $u=u^-_\mathrm{TT}\big(x;\tfrac{1}{2}+\ii p\big)$ of \eqref{eq:PII-intro} for $\alpha=\tfrac{1}{2}+\ii p$ is globally analytic and nonvanishing for real $x$. Applying the symmetries \eqref{eq:TT-symmetries} for $x\in\mathbb{R}$ then shows that $u_\mathrm{TT}^-\big(x,-\tfrac{1}{2}+\ii p\big)$, $u_\mathrm{TT}^+\big(x,\tfrac{1}{2}-\ii p\big)$, and $u_\mathrm{TT}^+\big(x,-\tfrac{1}{2}-\ii p\big)$ are global nonvanishing solutions for $x\in\mathbb{R}$ whenever $p>0$. This completes the proof of Theorem~\ref{thm:global}.

\appendix

\setcounter{section}{1}

\section*{Appendix A.~Parabolic cylinder function parametrix}\label{app:PC}\pdfbookmark[1]{Appendix A.~Parabolic cylinder function parametrix}{appendix}

\begin{figure}[t]\centering
\includegraphics{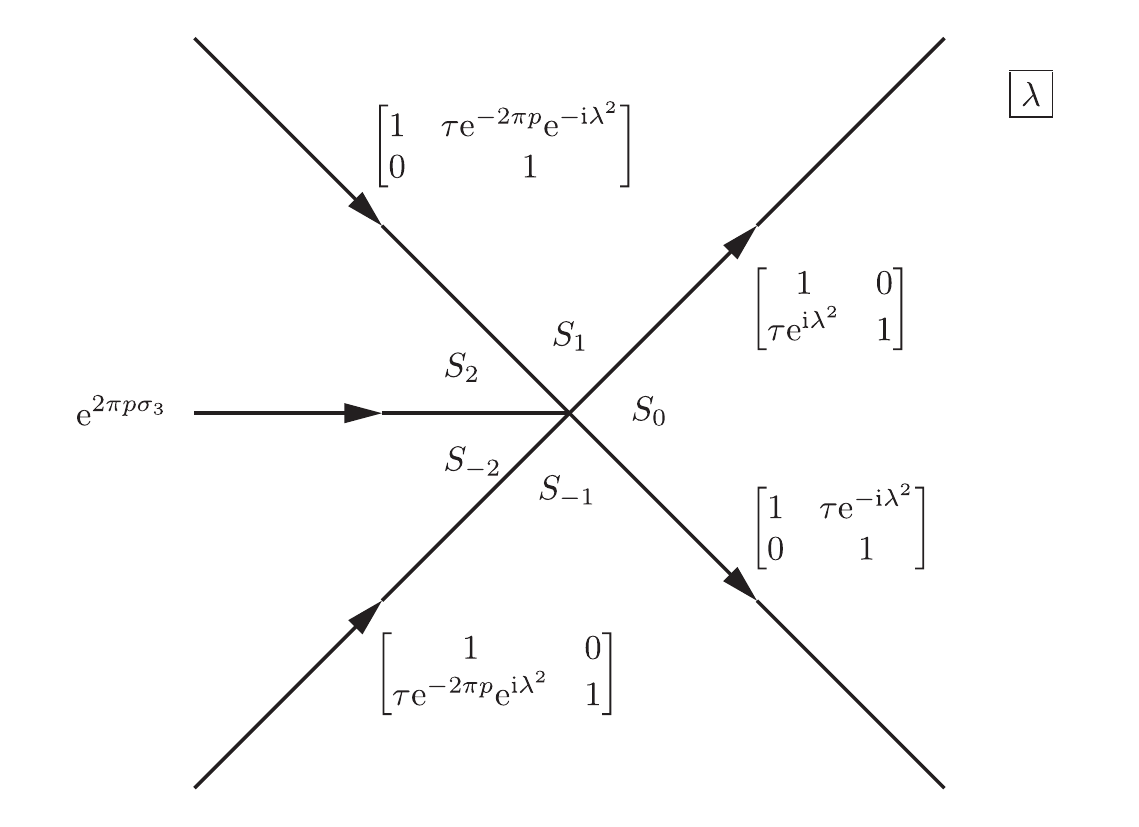}
\caption{The jump contour in the $\lambda$-plane and jump matrix $\mathbf{V}^{\mathrm{PC}}$.}\label{fig:PC-jumps-generalized}
\end{figure}

Consider the contour shown in Fig.~\ref{fig:PC-jumps-generalized} and the indicated jump matrix defined thereon.
\begin{rhp}[parabolic cylinder parametrix]\label{rhp:PC-generalized} Let $p,\tau\in\mathbb{C}$ be related by $\tau^2=\ee^{2\pi p}-1$. Seek a $2\times 2$ matrix-valued function $\mathbf{P}(\lambda)=\mathbf{P}(\lambda;p,\tau)$ with the following properties.
\begin{itemize}\itemsep=0pt
\item[]\textbf{Analyticity:} $\mathbf{P}(\lambda)$ is analytic for $\lambda$ in the five sectors $S_0\colon |\arg(\lambda)|<\tfrac{1}{4}\pi$, $S_1\colon \tfrac{1}{4}\pi<\arg(\lambda)<\tfrac{3}{4}\pi$, $S_{-1}\colon -\frac{3}{4}\pi<\arg(\lambda)<-\tfrac{1}{4}\pi$, $S_2\colon \tfrac{3}{4}\pi<\arg(\lambda)<\pi$, and $S_{-2}\colon -\pi<\arg(\lambda)<-\tfrac{3}{4}\pi$. It takes continuous boundary values on the excluded rays and at the origin from each sector.
\item[]\textbf{Jump conditions:} $\mathbf{P}_+(\lambda)=\mathbf{P}_-(\lambda)\mathbf{V}^\mathrm{PC}(\lambda)$, where $\mathbf{V}^\mathrm{PC}(\lambda)$ is the matrix function defined on the jump contour shown in Fig.~{\rm \ref{fig:PC-jumps-generalized}}.
\item[]\textbf{Normalization:} $\mathbf{P}(\lambda)\lambda^{\ii p\sigma_3}\to\mathbb{I}$ as $\lambda\to\infty$ uniformly in all directions.
\end{itemize}
\end{rhp}
To solve Riemann--Hilbert Problem~\ref{rhp:PC-generalized}, first note that there can be at most one solution, and if it exists it must have unit determinant. This problem has no solution if $p\in\ii\mathbb{Z}\setminus\{0\}$ for the same reason as in the case of Riemann--Hilbert Problem~\ref{rhp:PII-generalized}, and its solution is $\mathbf{P}(\lambda)\equiv\mathbb{I}$ if $p=0$.

Therefore we now assume that $p\not\in\ii\mathbb{Z}$, and consider the matrix $\mathbf{Q}(\lambda):=\mathbf{P}(\lambda)\ee^{-\ii\lambda^2\sigma_3/2}$. This change of variables removes the exponentials $\ee^{\pm\ii\lambda^2}$ from the jump conditions, from which it follows that $\mathbf{Q}'(\lambda)$ and $\mathbf{Q}(\lambda)$ satisfy exactly the same jump conditions. Appealing to invertibility and the specified regular behavior of $\mathbf{P}$ near $\lambda=0$ we see that $\mathbf{Q}'(\lambda)\mathbf{Q}(\lambda)^{-1}$ is an entire function of $\lambda$. If $\mathbf{P}(\lambda)$ has an asymptotic expansion as $\lambda\to\infty$ consistent with the normalization condition and hence of the form $\mathbf{P}(\lambda)=\big(\mathbb{I}+\lambda^{-1}\mathbf{P}^{[1]} + \cdots\big)\lambda^{-\ii p\sigma_3}$, the additional assumption that the series is differentiable term-by-term with respect to $\lambda$ shows by a Liouville argument that $\mathbf{Q}'(\lambda)\mathbf{Q}(\lambda)^{-1}=-\ii\lambda\sigma_3+\ii\big[\sigma_3,\mathbf{P}^{[1]}\big]$; therefore $\mathbf{Q}(\lambda)$ is a matrix solution of the first-order system
\begin{gather}
\frac{\dd\mathbf{Q}}{\dd\lambda}=\begin{bmatrix}-\ii\lambda & r\\s & \ii\lambda\end{bmatrix}\mathbf{Q},\qquad r:=2\ii P^{[1]}_{12},\qquad s:=-2\ii P^{[1]}_{21}.\label{eq:Y-system}
\end{gather}
Eliminating the second row of $\mathbf{Q}$ and rescaling by $t=\sqrt{2}\ee^{-\ii\pi/4}\lambda$ shows that the first row elements of $\mathbf{Q}$ satisfy Weber's equation for parabolic cylinder functions \cite[equation~(12.2.2)]{DLMF}
\begin{gather*}
\frac{\dd^2 Q_{1k}}{\dd t^2}-\left(\frac{1}{4}t^2+a\right)Q_{1k}=0,\qquad a:=\frac{1}{2}(1+\ii rs),\qquad k=1,2.
\end{gather*}
This differential equation has particular solutions $U(a,\pm t)$ and $U(-a,\pm\ii t)$, and the special function $U(\cdot,\cdot)$ is described in detail in \cite[Section~12]{DLMF}. In each of the five sectors, the general solution can be written as a linear combination of a suitable ``numerically satisfactory'' fundamental pair:
\begin{gather*}
Q_{1k}(\lambda)=\begin{cases}
r A^{(0)}_kU(a,t) + r B^{(0)}_kU(-a,\ii t), &\lambda\in S_0,\\
r A^{(1)}_kU(a,t)+r B^{(1)}_kU(-a,-\ii t), &\lambda\in S_1,\\
r A^{(-1)}_kU(a,-t)+r B^{(-1)}_kU(-a,\ii t), &\lambda\in S_{-1},\\
r A^{(2)}_kU(a,-t)+r B^{(2)}_kU(-a,-\ii t), &\lambda\in S_{2},\\
r A^{(-2)}_kU(a,-t)+r B^{(-2)}_kU(-a,-\ii t), &\lambda\in S_{-2}.
\end{cases}
\end{gather*}
The first row of the matrix differential equation \eqref{eq:Y-system} then gives the elements of the second row explicitly in terms of those of the first and their derivatives; however the derivatives can be eliminated using \cite[equations~(12.8.2)--(12.8.3)]{DLMF}, and we therefore obtain
\begin{gather*}
Q_{2k}(\lambda)=\sqrt{2}\ee^{-\ii\pi/4}\begin{cases}
-A^{(0)}_kU(a-1,t) + \ii \big(a-\tfrac{1}{2}\big)B^{(0)}_kU(1-a,\ii t), &\lambda\in S_0,\\
-A^{(1)}_kU(a-1,t) - \ii \big(a-\tfrac{1}{2}\big)B^{(1)}_kU(1-a,-\ii t), &\lambda\in S_1,\\
A^{(-1)}_kU(a-1,-t)+\ii \big(a-\tfrac{1}{2}\big)B^{(-1)}_kU(1-a,\ii t), &\lambda\in S_{-1},\\
A^{(2)}_kU(a-1,-t)-\ii \big(a-\tfrac{1}{2}\big)B^{(2)}_kU(1-a,-\ii t), &\lambda\in S_2,\\
A^{(-2)}_kU(a-1,-t)-\ii \big(a-\tfrac{1}{2}\big)B^{(-2)}_kU(1-a,-\ii t), &\lambda\in S_{-2}.
\end{cases}
\end{gather*}
In each sector, the asymptotic expansion \cite[equation~(12.9.1)]{DLMF}:
\begin{gather}
U\left(a,t\right)\sim \ee^{-\frac{1}{4}t^{2}}t^{-a-\frac{1}{2}}\sum_{j=0}^{\infty}(-1)^{j}\frac{{\big(\frac{1}{2}+a\big)_{2j}}}{j!(2t^{2})^{j}},\qquad t\to\infty,\qquad |\arg(t)|<\frac{3}{4}\pi \label{eq:U-asymp}
\end{gather}
can be used to obtain the asymptotic behavior of the matrix elements of $\mathbf{Q}(\lambda)$. Since the normalization condition on $\mathbf{P}(\lambda)=\mathbf{Q}(\lambda)\ee^{\ii\lambda^2\sigma_3/2} =\mathbf{Q}(\lambda)\ee^{-t^2\sigma_3/4}$ forbids exponential growth, it is necessary that $A_1^{(i)}=0$ and $B_2^{(i)}=0$ for all sector indices $i=0,\pm 1,\pm 2$. From the jump conditions it follows that the second column of $\mathbf{Q}$ must match between sectors $S_0$ and $S_1$, as well as between sectors $S_{-1}$ and $S_{-2}$. This implies that $A_2^{(0)}=A_2^{(1)}$ and that $A_2^{(-1)}=A_2^{(-2)}$. Similarly, the first column of $\mathbf{Q}$ must match between sectors $S_0$ and $S_{-1}$, as well as between sectors $S_1$ and $S_2$. This implies that $B_1^{(0)}=B_1^{(-1)}$ and that $B_1^{(1)}=B_1^{(2)}$. Lastly, the jump condition between sectors $S_2$ and $S_{-2}$ implies that $B_1^{(-2)}=\ee^{-2\pi p}B_1^{(2)}=\ee^{-2\pi p}B_1^{(1)}$ and that $A_2^{(2)}=\ee^{-2\pi p}A_2^{(-2)}=\ee^{-2\pi p}A_2^{(-1)}$. Therefore, if $\mathbf{Q}^{(i)}(\lambda)$ denotes the restriction of $\mathbf{Q}(\lambda)$ to $\lambda\in S_i$, we have so far found that{\allowdisplaybreaks
\begin{gather}
\mathbf{Q}^{(0)}(\lambda)=\begin{bmatrix}r B_1^{(0)}U(-a,\ii t) & r A_2^{(0)}U(a,t)\\
\sqrt{2}\ee^{\ii\pi/4}\big(a-\tfrac{1}{2}\big)B_1^{(0)}U(1-a,\ii t) & \sqrt{2}\ee^{3\ii\pi/4}A_2^{(0)}U(a-1,t)\end{bmatrix},\qquad\lambda\in S_0,\label{eq:Y0}
\\
\mathbf{Q}^{(1)}(\lambda)=\begin{bmatrix}r B_1^{(1)}U(-a,-\ii t) & r A_2^{(0)}U(a,t)\\
\sqrt{2}\ee^{-3\ii\pi/4}\big(a-\tfrac{1}{2}\big)B_1^{(1)}U(1-a,-\ii t) & \sqrt{2}\ee^{3\ii\pi/4}A_2^{(0)}U(a-1,t)
\end{bmatrix},\qquad\lambda\in S_1,\nonumber
\\
\mathbf{Q}^{(-1)}(\lambda)=\begin{bmatrix}r B_1^{(0)}U(-a,\ii t) & r A_2^{(-1)}U(a,-t)\\
\sqrt{2}\ee^{\ii\pi/4}\big(a-\tfrac{1}{2}\big)B_1^{(0)}U(1-a,\ii t) & \sqrt{2}\ee^{-\ii\pi/4}A_2^{(-1)}U(a-1,-t)
\end{bmatrix},\qquad\lambda\in S_{-1},\nonumber
\\
\mathbf{Q}^{(2)}(\lambda)=\begin{bmatrix}r B_1^{(1)}U(-a,-\ii t) & \ee^{-2\pi p}r A_2^{(-1)}U(a,-t)\\
\sqrt{2}\ee^{-3\ii\pi/4}\big(a-\tfrac{1}{2}\big)B_1^{(1)}U(1-a,-\ii t) & \ee^{-2\pi p}\sqrt{2}\ee^{-\ii\pi/4} A_2^{(-1)}U(a-1,-t)\end{bmatrix},\nonumber\\
\hphantom{\mathbf{Q}^{(2)}(\lambda)=}{} \ \lambda\in S_2,\nonumber
\\
\mathbf{Q}^{(-2)}(\lambda)=\begin{bmatrix}\ee^{-2\pi p}r B_1^{(1)} U(-a,-\ii t) & r A_2^{(-1)} U(a,-t)\\
\ee^{-2\pi p}\sqrt{2}\ee^{-3\ii\pi/4}\big(a-\tfrac{1}{2}\big)B_1^{(1)} U(1-a,-\ii t) & \sqrt{2}\ee^{-\ii\pi/4}A_2^{(-1)} U(a-1,-t)\end{bmatrix},\nonumber\\
\hphantom{\mathbf{Q}^{(-2)}(\lambda)=}{} \ \lambda\in S_{-2}.\label{eq:Ym2}
\end{gather}
Now} taking $\lambda\to\infty$ in each sector we find that the matrix $\mathbf{P}(\lambda)=\mathbf{Q}(\lambda)\ee^{\ii\lambda^2\sigma_3/2}$ satisfies the required normalization condition provided that
\begin{gather}
a=\frac{1}{2}-\ii p \qquad\text{which implies that}\quad rs=-2p, \label{eq:a-from-p}
\end{gather}
and that the remaining unknown coefficients are subject to the following:
\begin{alignat}{3}
& r B_1^{(0)}= \ee^{-\pi p/4}\ee^{\ii p\ln(2)/2}, \qquad && A_2^{(0)}=\frac{\ee^{-3\ii \pi/4}}{\sqrt{2}}\ee^{-\pi p/4}\ee^{-\ii p\ln(2)/2}, &\nonumber\\
& r B_1^{(1)}= \ee^{3\pi p/4}\ee^{\ii p\ln(2)/2},\qquad && A_2^{(-1)}=\frac{\ee^{\ii \pi/4}}{\sqrt{2}}\ee^{3\pi p/4}\ee^{-\ii p\ln(2)/2}.& \label{eq:coefficient-identities}
\end{alignat}
Only the values of the constants $r$ and $s$ remain undetermined. It is interesting to observe that if one tries to use the definitions in \eqref{eq:Y-system} to find $r$ and $s$ from the formul\ae\ \eqref{eq:Y0}--\eqref{eq:Ym2} after substituting from \eqref{eq:a-from-p} and \eqref{eq:coefficient-identities} and using the asymptotic expansion \eqref{eq:U-asymp} for the off-diagonal entries, one finds simply the truisms $r=r$ and $s=s$. These constants must therefore be determined from the only other information we have not yet used: the nontrivial jumps captured by the first (resp.\@ second) column of the jump matrix $\mathbf{V}^\mathrm{PC}$ for $\arg(\lambda)=\tfrac{1}{4}\pi,-\tfrac{3}{4}\pi$ (resp.\@ for $\arg(\lambda)=-\tfrac{1}{4}\pi,\tfrac{3}{4}\pi$). In particular, the jump condition for $\arg(\lambda)=\tfrac{1}{4}\pi$ requires that
\begin{gather*}
Q^{(1)}_{11}(\lambda)=Q^{(0)}_{11}(\lambda)+\tau Q^{(0)}_{12}(\lambda),
\end{gather*}
which given all available information is equivalent to the condition
\begin{gather*}
\ee^{\ii\pi(a+1/2)/2}U(-a,-\ii t)+\ee^{-\ii\pi(1+1/2)/2}U(-a,\ii t)\\
\qquad{} =\tau r\frac{\ee^{-\ii\pi/4}}{\sqrt{2}}\ee^{-\pi p/2}\ee^{-\ii p\ln(2)}U(a,t),\qquad a=\frac{1}{2}-\ii p,\qquad \arg(t)=0.
\end{gather*}
But it is easy to check that this matches the connection formula \cite[equation~(12.2.18)]{DLMF} for the parabolic cylinder function $U(\cdot,\cdot)$ provided that
\begin{gather}
r=r(p,\tau):=2\ee^{\ii\pi/4}\sqrt{\pi}\frac{\ee^{\pi p/2}\ee^{\ii p\ln(2)}}{\tau\Gamma(\ii p)}.\label{eq:alpha-def}
\end{gather}
Using \eqref{eq:a-from-p}
one then deduces that
\begin{gather}
s=s(p,\tau)=-\frac{2 p}{r(p,\tau)}=\frac{\ee^{3\ii\pi/4}}{\sqrt{\pi}}\tau p\Gamma(\ii p)\ee^{-\pi p/2}\ee^{-\ii p\ln(2)}.
\label{eq:beta-def}
\end{gather}
At last all free parameters have been determined, and it is straightforward to verify with the use of connection formul\ae\ for $U(\cdot,\cdot)$ that the other jump conditions that we have not yet enforced are automatically satisfied.

\begin{Remark}In fact, the unenforced jump conditions cannot contain any further information, because the cyclic product of the jump matrices for $\mathbf{Q}$ is the identity (as must be the case for consistency at the origin). This matrix identity modulo unit determinant amounts to three conditions, which determine all of the off-diagonal entries of $\mathbf{V}^\mathrm{PC}$ given any one of them.
\end{Remark}

This completes the solution of Riemann--Hilbert Problem~\ref{rhp:PC-generalized}. We pause to point out that the normalization condition on $\mathbf{P}(\lambda)=\mathbf{P}(\lambda;p,\tau)$ indeed holds in every sense we assumed at the start in order to derive the differential equation \eqref{eq:Y-system}, and therefore in particular from the asymptotic expansion \eqref{eq:U-asymp},
\begin{gather}
\mathbf{P}(\lambda;p,\tau)\lambda^{\ii p\sigma_3}=\mathbb{I} + \frac{1}{2\ii\lambda}\begin{bmatrix}0 & r(p,\tau) \\ -s(p,\tau) & 0\end{bmatrix}+\begin{bmatrix}O\big(\lambda^{-2}\big) & O\big(\lambda^{-3}\big)\\O\big(\lambda^{-3}\big) & O\big(\lambda^{-2}\big)\end{bmatrix},\qquad\lambda\to\infty.\label{eq:U-matrix-asymp-generalized}
\end{gather}
In fact, \eqref{eq:U-asymp} shows that the diagonal (resp.\ off-diagonal) elements have explicit asymptotic expansions in descending even (resp.\ odd) powers of $\lambda$. The full expansion of $\mathbf{P}(\lambda)\lambda^{\ii p\sigma_3}$ as $\lambda\to\infty$ is the same in all five sectors.

In the special case that $p>0$, which also implies that $\tau\in\mathbb{R}$, we may use the identity for the modulus of the gamma function on the imaginary axis (see \cite[equation~(5.4.3)]{DLMF}) to deduce that $s=-r^*$.

\subsection*{Acknowledgements} The author's work was supported by the National Science Foundation under grants DMS-1513054 and DMS-1812625. The author thanks Thomas Bothner, Deniz Bilman, and Liming Ling for useful discussions.

\pdfbookmark[1]{References}{ref}
\LastPageEnding

\end{document}